\documentclass{article}

\usepackage{url}
\usepackage{soul}
\usepackage{xurl}
\usepackage{cite}
\usepackage{arxiv}
\usepackage{array}
\usepackage{float}
\usepackage{lipsum}
\usepackage{lineno}
\usepackage{amssymb}
\usepackage{amsmath}
\usepackage{authblk}
\usepackage{siunitx}
\usepackage{amsfonts}
\usepackage{placeins}
\usepackage{geometry}
\usepackage{tabularx}
\usepackage{booktabs}
\usepackage{multirow}
\usepackage{hyperref}
\usepackage{graphicx}
\usepackage{makecell}
\usepackage{nicefrac}
\usepackage{microtype}
\usepackage{tcolorbox}
\usepackage{adjustbox}
\usepackage{dirtytalk}
\usepackage{subcaption}
\usepackage[T1]{fontenc}
\usepackage[utf8]{inputenc}
\usepackage[version=4]{mhchem}

\title{Stochastic Generalized-Order Constitutive Modeling of Viscoelastic Spectra of Polyurea-Graphene Nanocomposites}

\author[1]{Arman Khoshnevis}
\author[2,3]{Demetrios A. Tzelepis}
\author[2]{Valeriy V. Ginzburg}
\author[1,4,*]{Mohsen Zayernouri}

\affil[1]{Department of Mechanical Engineering, Michigan State University, East Lansing, MI 48824, USA}
\affil[2]{Department of Chemical Engineering and Materials Science, Michigan State University, East Lansing, MI 48824, USA}
\affil[3]{Materials Division, US-Army, Ground Vehicle System Center, Warren, MI 48397, USA}
\affil[4]{Department of Statistics and Probability, Michigan State University, East Lansing, MI 48824, USA}
\affil[*]{\small \textit{Corresponding author}}

\begin{document}

\maketitle

\footnotetext[1]{Email: khoshne1@msu.edu}
\footnotetext[2]{Email: tzelepi1@msu.edu}
\footnotetext[3]{Email: ginzbur7@msu.edu}
\footnotetext[4]{Email: zayern@msu.edu}

\begin{abstract}
Polyurea (PUa) elastomers are extensively used in a wide range of applications spanning from biomedical to defense fields due to their enabling mechanical properties. These materials can be further reinforced through the incorporation of nanoparticles to form nanocomposites. This study focuses on an IPDI-based PUa matrix with exfoliated graphene nanoplatelet (xGnP) fillers. We propose a generalized constitutive model by integrating one Fractional Maxwell Model (FMM) and one Fractional Maxwell Gel (FMG) branch in a parallel configuration via introducing a new dimensionless number to bridge between these branches physically and mathematically. Through systematic local-to-global sensitivity analyses, we investigate the behavior of these nanocomposites to facilitate simulation, design, and performance prediction. Consistently, the constructed models share the same most/least influential model parameters. \(\alpha_1\) and \(E_{c_1}\), the power exponent and the characteristic modulus of the first branch, are found to be the most influential model parameters, while \(\tau_{c_2}\) and \(\tau_{c_1}\), the characteristic time-scales of each branch, are recognized as the least influential model parameters. The proposed PU nanocomposite constitutive laws can now make an impact to the design and optimization of coating and shock-absorbing coatings in a range of applications.
\end{abstract}


\section{Introduction}
Polyurethanes (PUs), polyureas (PUas), and poly(urethaneureas) (PUUs) are amongst an interesting division of polymers with a broad scope of industrial applications \cite{petrovic_polyurethane_1991, szycher_polyurethanes_2013, akindoyo_polyurethane_2016, das_brief_2020, sonnenschein_polyurethanes_2021}, including coatings \cite{chattopadhyay_thermal_2005, chattopadhyay_structural_2007, driffield_method_2007}, shock-absorbing coatings \cite{ding2024recent}, adhesives \cite{cognard_handbook_2005, lee_preparation_2008, jia-hu_synthesis_2015, meier-westhues_polyurethanes_2019}, foams \cite{bicerano_flexible_2004, harikrishnan_modeling_2010, allan_thermal_2013, dsouza_polyurethane_2014}, and elastomers \cite{furukawa_microphase-separated_2005, bagdi_quantitative_2012}. PUas are generally considered as a multi-block copolymer consisting of altering soft segment blocks (typically a polyamine) and hard segment blocks (formed by a reaction between a diisocyanate and a chain extender) \cite{koberstein_smallangle_1983, koberstein_smallangle_1983-1, leung_smallangle_1985, christenson_model_1986, koberstein_compression-molded_1992, koberstein_compression-molded_1992-1}. Depending on the hard segment weight fraction (HSWF) and the length of the soft segments, PUas can exhibit a wide spectrum of physical and mechanical properties ranging from hard and brittle to soft and elastomeric, making them suitable for various applications \cite{ginzburg_theoretical_2007}. Due to the strong hydrogen bonding within the hard segments and their thermodynamic incompatibility with the soft segments \cite{benoit_scattering_1988, heintz_spectroscopic_2005}, a microphase separation occurs creating various morphologies such as spherical \cite{christenson_relationship_2005}, cylindrical \cite{garrett_microdomain_2001}, lamellar \cite{kaushiva_uniaxial_2000}, etc. These morphologies are similar to those observed in diblock copolymers and correspond to different values of composition, \(f\), and segregation strength, \(\chi N\) \cite{aou_characterization_2013}.

Data-driven constitutive modeling and model identification, particularly for non-linear viscoelastic models \cite{as2023mechanics, abdolazizi2024viscoelastic}, fractional-order constitutive models \cite{suzuki2021data, dabiri2023fractional}, and soft matters \cite{gonzalez2023model, upadhyay2024physics}, have gained significant popularity among researchers in recent years. These efforts mainly aim to develop physics-informed machine learning or deep learning frameworks based on synthesized and experimental stress-strain datasets in the time domain. Alongside this trend, there has been a continuous effort to conduct numerical simulations of such materials across various applications \cite{khoshnevis2024double, doherty2023stabilised, gruninger2024benchmarking}.

Moreover, performing sensitivity analysis (SA) is a crucial step in developing reliable constitutive models for materials with complex mechanical behavior. As such models often include numerous parameters, standard parameter identification procedures do not provide insights into the relative importance of these parameters on the output. Commonly used sensitivity analysis techniques include normalized partial derivatives \cite{saltelli1999sensitivity} and One-At-a-Time (OAT) \cite{morris1991factorial} methods, which are considered local SA techniques, as well as variance-based Sobol’ indices \cite{sobol1990sensitivity} and the Extended Fourier Amplitude Sensitivity Test (EFAST) \cite{saltelli1999quantitative}, which are considered global SA techniques. Implementation of these approaches is well-documented in the literature for a wide range of materials, including composites \cite{sasikumar2023sensitivity}, biological composites \cite{rao2023frequency}, cementitious composites \cite{zhou2020global}, flexoelectric materials \cite{hamdia2018sensitivity}, salt caverns \cite{khaledi2016sensitivity}, and standard materials used in transmission lines \cite{de2024thermo, kc2024thermo, kc2024multi}.

Next, a brief introduction of fractional-order constitutive modeling is reviewed. Scott-Blair \cite{blair_role_1947} proposed a constitutive equation to describe the stress response of complex materials through a fractional derivative of strain \cite{jaishankar_power-law_2013}
\begin{linenomath}
\begin{equation}
\sigma(t) = \mathbb{V} \frac{d^\alpha \epsilon(t)}{dt^\alpha},
\end{equation}
\end{linenomath}
where \(\frac{d^\alpha }{dt^\alpha}\) is the fractional-order derivative operator and \(\mathbb{V}\) (with SI unit of \(\mathrm{Pa\,s^{\text{\(\alpha\)}}}\)) is a quasi-material property representing the \say{firmness} of a material \cite{blair_subjective_1942}. These two parameters fully characterize the constitutive model. This model can also be described by \(E\), modulus, \(\tau\), time-scale, and \(\alpha\), power exponent; however, only the combination of \((E\tau^{\alpha})\) and \(\alpha\) are measurable experimentally \cite{schiessel_generalized_1995}. The model is symbolized by an element known as the spring-pot, which interpolates between a purely elastic spring \((\alpha = 0)\) and a purely viscous dashpot \((\alpha = 1)\). Therefore, it is effectively captures a wide range of relaxation time-scales inherent to complex rheological materials \cite{jaishankar_power-law_2013}.

\begin{figure}[t]
    \includegraphics[width=6.5 cm]{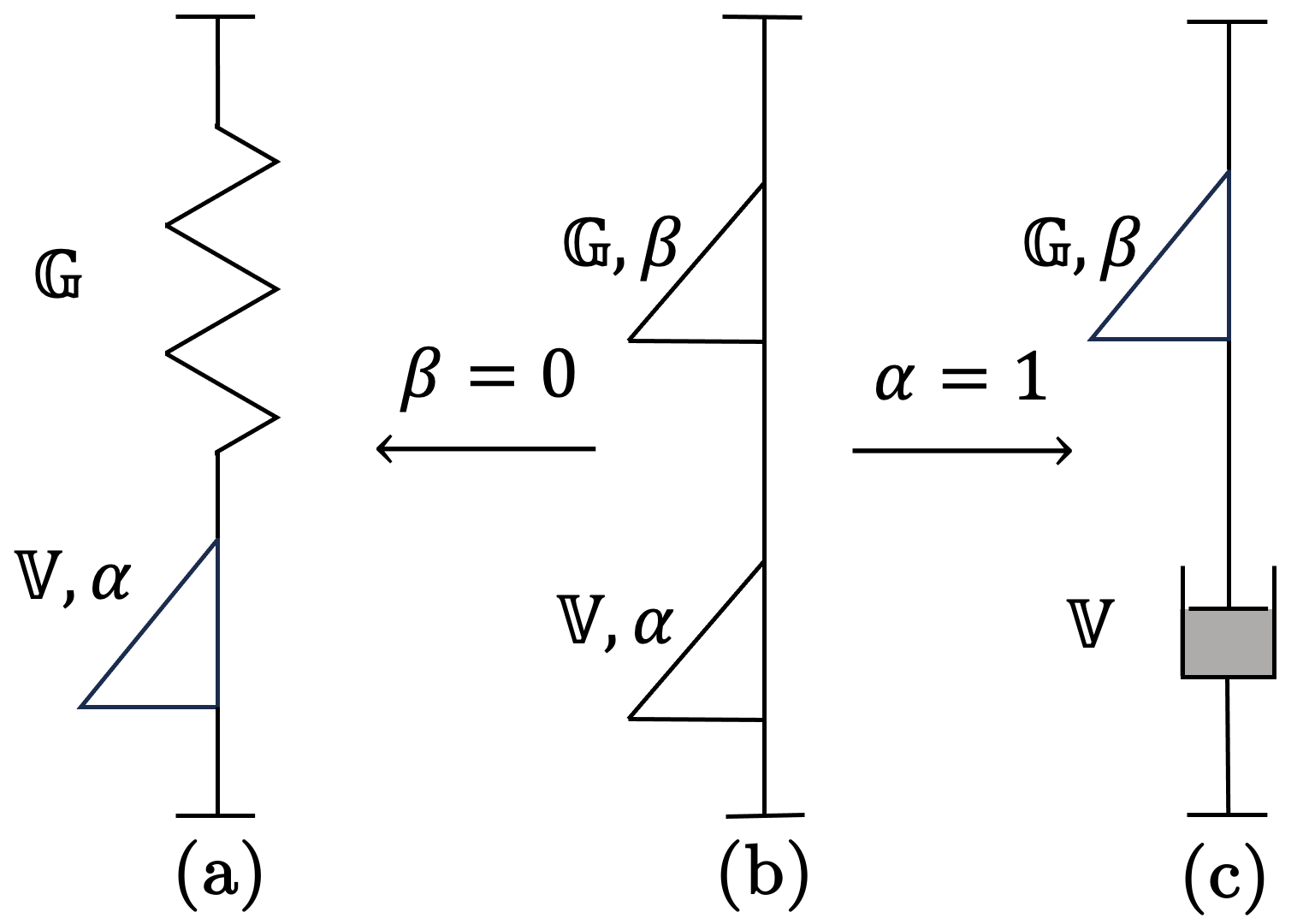}
    \centering
    \caption{Schematic illustration of a (b) Fractional Maxwell Model (FMM) and its two limiting cases: (a) Fractional Maxwell Gel (FMG) and (c) Fractional Maxwell Liquid (FML).} \label{fig-fmmsimplified}
\end{figure}

Configuring spring-pots in different arrangements paves the way of developing more intricate constitutive models. To begin with, placing two spring-pots, characterized by \((\mathbb{V},\alpha)\) and \((\mathbb{G},\beta)\), in series (as depicted in Figure \ref{fig-fmmsimplified}b), and assuming equality of stresses \((\sigma=\sigma_1 = \sigma_2)\) and additivity of strains \((\epsilon = \epsilon_1 + \epsilon_2)\), results in the Fractional Maxwell Model (FMM) \cite{schiessel_generalized_1995, jaishankar_power-law_2013}
\begin{linenomath}
\begin{equation} \label{eq-fmm}
\sigma(t) + \frac{\mathbb{V}}{\mathbb{G}} \frac{d^{\alpha-\beta} \epsilon(t)}{dt^{\alpha-\beta}} = \mathbb{V} \frac{d^\alpha \epsilon(t)}{dt^\alpha},
\end{equation}
\end{linenomath}
where the inequality \(0\leq\beta<\alpha\leq1\) ensures the thermodynamic consistency \cite{lion_thermodynamics_1997} without loss of generality \footnote{A thorough derivation of this constitutive model is elaborated in the literature; for instance Schiessel \cite{schiessel_generalized_1995}.}. By applying the Fourier transform to \eqref{eq-fmm}, the complex modulus can be derived as \cite{rathinaraj_incorporating_2021}:
\begin{linenomath}
\begin{equation}
E^*(\omega) = \frac{\mathbb{V}(i\omega)^{\alpha} \cdot \mathbb{G}(i\omega)^{\beta}}{\mathbb{V}(i\omega)^{\alpha} + \mathbb{G}(i\omega)^{\beta}}.
\end{equation}
\end{linenomath}
This expression can be non-dimensionalized and written as
\begin{linenomath}
\begin{equation}
\frac{E^*(\omega)}{E_c} = \frac{(i\omega\tau_c)^{\alpha}}{1 + (i\omega\tau_c)^{\alpha-\beta}},
\end{equation}
\end{linenomath}
where \(E_c \equiv (\mathbb{G}^{\alpha}/\mathbb{V}^{\beta})^{1/(\alpha-\beta)}\) and \(\tau_c \equiv (\mathbb{V}/\mathbb{G})^{1/(\alpha-\beta)}\) \cite{rathinaraj_incorporating_2021}. Furthermore, the storage and loss moduli can be obtained by splitting the complex modulus into its real and imaginary components, respectively as
\begin{subequations}
\begin{linenomath}
\begin{equation}
\frac{E^{\prime}}{E_c} = \frac{(\omega\tau_c)^\alpha \cos\left(\frac{\pi\alpha}{2}\right) + (\omega\tau_c)^{2\alpha-\beta} \cos\left(\frac{\pi\beta}{2}\right)}{1 + (\omega\tau_c)^{\alpha-\beta} \cos\left(\frac{\pi(\alpha-\beta)}{2}\right) + (\omega\tau_c)^{2(\alpha-\beta)}},
\end{equation}
\end{linenomath}

\begin{linenomath}
\begin{equation}
\frac{E^{\prime\prime}}{E_c} = \frac{(\omega\tau_c)^\alpha \sin\left(\frac{\pi\alpha}{2}\right) + (\omega\tau_c)^{2\alpha-\beta} \sin\left(\frac{\pi\beta}{2}\right)}{1 + (\omega\tau_c)^{\alpha-\beta} \cos\left(\frac{\pi(\alpha-\beta)}{2}\right) + (\omega\tau_c)^{2(\alpha-\beta)}}.
\end{equation}
\end{linenomath}
\end{subequations}

As noted, the spring-pot can be simplified into either a pure spring or a dashpot element, depending on the value of its power exponent. This adaptability gives rise to two limiting cases for the Fractional Maxwell Model (FMM). First, setting the exponent \(\alpha\) to unity results in the Fractional Maxwell Liquid (FML) (Figure \ref{fig-fmmsimplified}c), which is best suited for describing the behavior of complex fluids in the pre-gel state. Second, letting the power exponent \(\beta\) approaches zero leads to the Fractional Maxwell Gel (FMG) model (Figure \ref{fig-fmmsimplified}a), capable of capturing the elastic behavior of materials beyond their gel point \cite{rathinaraj_incorporating_2021}.

In the present study, we extend the fractional modeling approach by representing the soft phase with a general fractional Maxwell model (FMM) branch. This modification aims to parsimoniously capture the viscoelastic behavior of both the neat and nanocomposite systems and reduce modeling errors across 10 decades of frequency. To achieve this goal, we employ a multi-objective Particle Swarm Optimization (PSO) via an in-house developed MATLAB code to fit the modeled linear viscoelastic properties, namely, the storage and loss moduli, to the experimental dataset corresponding to various polyurea nanocomposites. Notably, we introduce a new dimensionless number and incorporate it into the optimization process to establish a physical and mathematical connection between the two fractional branches.

Furthermore, owing to the inherent modeling and experimental errors, we shift our perspective from a frequentist approach to a stochastic one, recognizing that slightly different model parameters may be obtained for each new set of neat or nanocomposite systems with the same specifications. Consequently, we seek to answer the following questions for the FMG-FMG and FMM-FMG models: which model parameter(s), when fixed at any value within their range of variability, reduce the variance of the model output the most? which model parameter(s) can be considered non-influential and treated deterministically? To address these questions, we implement a derivative-based local sensitivity analysis and a variance-based global sensitivity analysis, quantifying the effect of parameters on the model outputs.

The structure of the remainder of the paper is organized as follows: Section \ref{sec.M&M} provides a summarized description of the synthesis process for the nanocomposite systems and their characterization methods, the fractional constitutive modeling approach, the optimization technique employed, and local and global sensitivity analyses procedures adopted. Section \ref{sec.Exp&Mod} presents the fitting of the fractional model to the experimental datasets, elucidates the effect of the HSWF and xGnP weight fractions on model parameters, and discusses factor prioritization and factor fixing based on the local and global sensitivity analysis. Finally, Section \ref{sec.con} concludes the paper with remarks summarizing the key findings and insights.

\section{Materials and Methods} \label{sec.M&M}
\subsection{Polymer and Nanocomposite Synthesis}
The synthesis of three polyurea (PUa)-Neat materials, characterized by the hard segment weight fractions (HSWF) of 20\%, 30\%, and 40\%, as well as their corresponding nanocomposite systems, incorporating exfoliated nano-graphene platelets (xGnP) weight fractions of 0.5, 1.0, and 1.5 wt.\%, has been detailed in our prior studies \cite{tzelepis_experimental_2023, tzelepis_polyureagraphene_2023}. The chemical constituents employed in this synthesis included isophorone diisocyanate (IPDI)-VESTANAT, sourced from Evonik corporation, Piscataway, NJ, USA; Jeffamine T5000 and D2000, provided by Huntsman Corporation, The Woodlands, TX, USA; diethyltoluenediamine (DETDA) (Lonzacure) chain extender, acquired from Lonza, Morristown, NJ, USA; and Toluene, purchased from Fisher Scientific, Hampton, NH, USA. Additionally, the exfoliated nano-graphene (grade R-10) was procured from XG Sciences. All materials were utilized as \say{received}, without further processing. The formulations for the three neat PUa-s with HSWF of 20\%, 30\%, and 40\%, are listed in Table \ref{tab1}. The sample preparation process is summarized as follows:

The polyurea prepolymer (A-side) was formed by introducing \ce{IPDI} into the reactor, followed by the addition of toluene to mitigate potential gelling. Concurrently, a blend of Jeffamine D2000 and T5000 was mixed for \SI{5}{\minute} in a separate beaker, degassed for \SI{10}{\minute}, and then added to the \ce{IPDI}–toluene mixture. For the B-side, a distinct blend of Jeffamine D2000 and Lonzacure \ce{DETDA} was prepared by mixing in another beaker for \SI{5}{\minute}, and subjected to vacuum degassing for approximately 10 min, and then poured into a separate additional funnel. Upon assembling the reactor, a vacuum was applied for \SI{5}{\minute}, followed by the addition of \ce{N2} gas at a flow rate of \SIrange{0.3}{0.4}{\liter\per\second}. The reactor temperature was raised to \SI{80}{\degreeCelsius}, and the A-side amine blend was added drop-wise under mechanical stirring at 120 rpm. This mixture underwent continuous stirring for an additional hour at \SI{80}{\degreeCelsius}. Subsequently, the reactor was cooled to \SI{0}{\degreeCelsius}, and the B-side amine blend was added in a drop-wise manner. The material was then transferred to molds, where it underwent gelation and solvent evaporation at room temperature for \SI{24}{\hour}, followed by further solvent evaporation in an oven at \SI{40}{\degreeCelsius} for \SIrange{12}{24}{\hour}. The final curing of the polyurea occurred at \SI{60}{\degreeCelsius} for \SI{72}{\hour}.

\begin{table}[H]
    \caption{Summary of the constituents, used in the synthesis of the model PUa-Neat. (Adapted with permission from Ref. \cite{tzelepis_experimental_2023}.)\label{tab1}}
    \centering
    \begin{tabular}{lcccc}
        \toprule
        \textbf{Sides} & \textbf{Component} & \textbf{IPDI-2k-20HS} & \textbf{IPDI-2k-30HS} & \textbf{IPDI-2k-40HS}\\
        \midrule
        \multirow{5}{*}{\makecell[l]{Isocyanate \\ Prepolymer \\ (A-Side)}} & IPDI & 30.8 g & 41.6 g & 52.1 g \\
        & T5000 & 14.5 g & 12.1 g & 10.1 g \\
        & D2000 & 57.1 g & 48.5 g & 40.4 g \\
        & Toluene & \SI{82.7}{g} (\SI{95}{mL}) & \SI{165.3}{g} (\SI{190}{mL}) & \SI{208.8}{g} (\SI{240}{mL}) \\
        & \%NCO & 8.7\% & 12.9\% & 15.3\% \\
        \midrule
        \multirow{2}{*}{\makecell[l]{Amine Blend \\ (B-Side)}} & DETDA & 10 g & 19.4 g & 29.3 g \\
        & D2000 & 90 g & 80.6 g & 70.7 g \\
        \bottomrule
    \end{tabular}
\end{table}

The synthesis process for all PUa-GnP nanocomposites mirrored that of the neat systems, with the following additional steps. Exfoliated nano-graphene (grade R-10) underwent heat treatment at 400°C for \SI{1}{\hour}, followed by a furnace cooling. In a \SI{500}{mL} beaker, the required amount of xGnPs was combined with \SI{190}{mL} of toluene and subjected to simultaneous mechanical stirring and sonication. Magnetic stirring at 200 rpm was employed for mechanical agitation, while sonication was performed using a Qsonica sonicator, manufactured by Qsonica L.L.C, Newt.own, CT, USA. The sonication parameters included an amplitude of 20 and a process time of 30 min with a pulse time of \SI{10}{\second} on and \SI{10}{\second} off, ensuring the slurry temperature remained below \SI{32}{\degreeCelsius} throughout the process. The total energy input amounted to \SI{38,610}{J} over a runtime of \SI{\approx 1}{\hour}, with xGnP weights summarized as follows: \SI{1.02}{g} for 0.5 wt.\% xGnP formulations, \SI{2.04}{g} for 1.0 wt.\% xGnP formulations, and \SI{3.06}{g} for 1.5 wt.\% xGnP formulation.

\subsection{PUa–xGnP Experimental Characterization}
Several experimental analyses were conducted to investigate various aspects of both the neat materials and nanocomposites, as detailed in \cite{tzelepis_experimental_2023, tzelepis_polyureagraphene_2023}. This study centers primarily on the master curves derived from Dynamic Mechanical Analysis (DMA), and its procedure is outlined as follows. DMA tests were carried out using a TA Instruments RSA-G2 rheometer (New Castle, DE, USA) in a tensile configuration, with all film samples initially loaded in tension. Temperature sweeps were performed with six repetitions per formulation, spanning a range from \SI{-95}{\degreeCelsius} to a maximum temperature determined by the hard segment content of the polyurea formulation, while maintaining the heating rate at \SI{3}{\degreeCelsius\per\minute}. The reference temperature for each material was set to its glass transition temperature, defined as the peak of the loss modulus (see Table S1). All Time-Temperature Superposition (TTS) shifts were executed using TA Instruments' TRIOS software package, version 5.1.1.46572.

\subsection{Mathematical Modeling} \label{sec.MM}
\begin{figure}[H]
    \centering
    \begin{subfigure}{0.4\textwidth}
        \includegraphics[width=\textwidth]{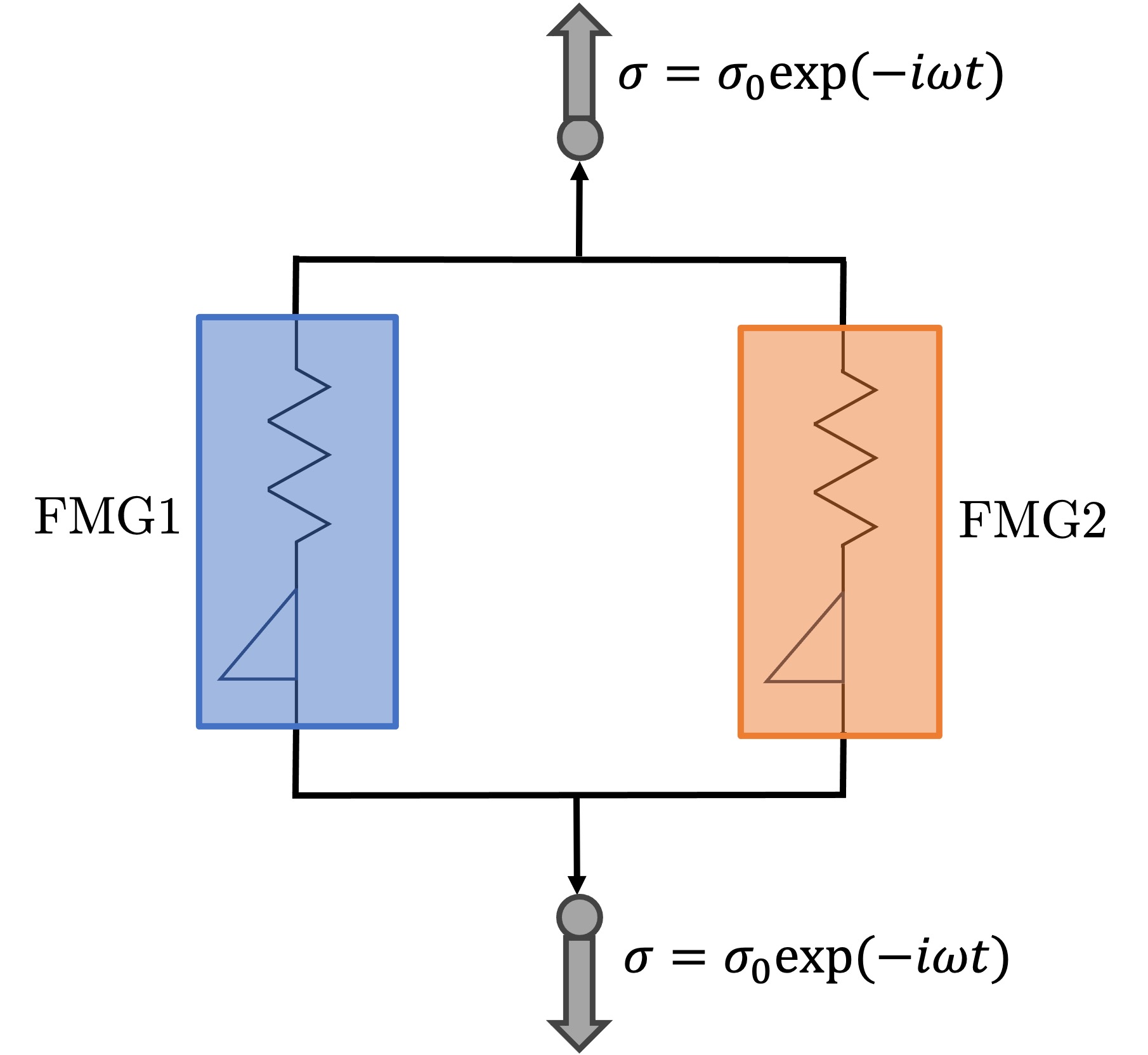}
        \caption{\label{fig-FMG-model}}
    \end{subfigure}
    \begin{subfigure}{0.4\textwidth}
        \includegraphics[width=\textwidth]{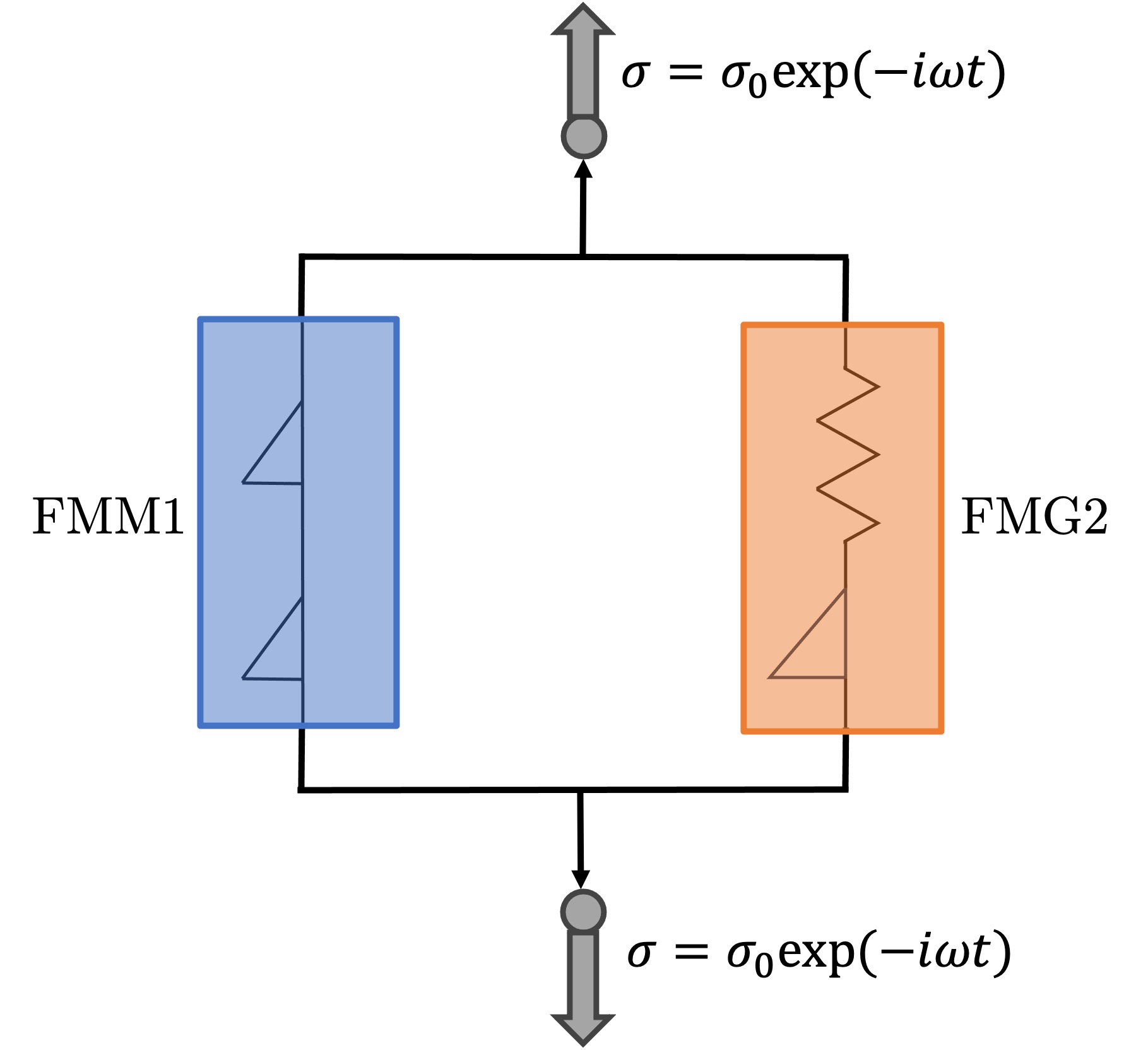}
        \caption{\label{fig-FMM-model}}
    \end{subfigure}
    \caption{Schematic of two constitutive models considered in the present study: (a) FMG-FMG model: The FMG1 branch corresponds to the filler soft phase, and the FMG2 branch corresponds to the percolated hard phase (both consist of a spring-pot and a spring in series). (b) FMM-FMG model: the FMM1 branch (consists of two spring-pots in series) corresponds to the filled soft phase, whereas the FMG2 branch corresponds to the percolated hard phase.} \label{fig-models}
\end{figure}

The complexity of a single-branch FMM model can be further advanced through the parallel configuration of FMM and FMG branches. In our preceding research \cite{tzelepis_experimental_2023, tzelepis_polyureagraphene_2023}, we examined a constitutive model consisting of two parallel FMG branches, one corresponds to the soft-phase matrix and the other one to the percolated hard phase, ultimately forming a FMG-FMG model, as depicted in Figure \ref{fig-FMG-model}. It is  important to highlight that an additional branch is not designated to the incorporated nanoparticles since they are assumed to be uniformly dispersed within the two phases without forming a new phase by themselves due to their low mass fraction. In the present study, we introduced another model by replacing the first FMG branch with an FMM branch, resulting in the FMM-FMG model, illustrated in Figure \ref{fig-FMM-model}. This modification aims to more accurately represent the viscoelastic behavior of the soft phase segment. Consequently, in both models, characteristic moduli (\(E_{c_1}, E_{c_2}\)), relaxation time-scales (\(\tau_{c_1}, \tau_{c_2}\)), and power exponents (\(\alpha_1, \alpha_2\)), along with an additional power exponent (\(\beta_1\)) associated with the second spring-pot element in the FMM branch, constitute the six to seven model parameters essential for the characterization of nanocomposite systems.

DMA data corresponds to multiple temperatures are combined utilizing the time-temperature superposition (TTS) principal and assuming both phases share the same shift factor \(a_T\). Therefore, the equivalent storage and loss moduli of our model can be expressed as
\begin{subequations}
\begin{linenomath}
\begin{equation}  \label{eq-EpSum}
E^{\prime}(x) = \sum_{k=1}^2 E_{c_k} \frac{\left(x\tau_{c_k}\right)^{\alpha_k} \cos\left(\frac{\pi\alpha_k}{2}\right) + \left(x\tau_{c_k}\right)^{2\alpha_k-\beta_k} \cos\left(\frac{\pi\beta_k}{2}\right)}{1 + \left(x\tau_{c_k}\right)^{\alpha_k-\beta_k} \cos\left(\frac{\pi\left(\alpha_k-\beta_k\right)}{2}\right) + \left(x\tau_{c_k}\right)^{2\left(\alpha_k-\beta_k\right)}},
\end{equation}
\end{linenomath}

\begin{linenomath}
\begin{equation}  \label{eq-EppSum}
E^{\prime\prime}(x) = \sum_{k=1}^2 E_{c_k} \frac{(x\tau_{c_k})^{\alpha_k} \sin\left(\frac{\pi\alpha_k}{2}\right) + \left(x\tau_{c_k}\right)^{2\alpha_k-\beta_k} \sin\left(\frac{\pi\beta_k}{2}\right)}{1 + \left(x\tau_{c_k}\right)^{\alpha_k-\beta_k} \cos\left(\frac{\pi\left(\alpha_k-\beta_k\right)}{2}\right) + \left(x\tau_{c_k}\right)^{2\left(\alpha_k-\beta_k\right)}}.
\end{equation}
\end{linenomath}
\end{subequations}
Here, \(x=a_T\omega\), with \(\beta_{1,2}=0\) for the FMG-FMG model, and \(\beta_2=0\) for the FMM-FMG model. It should be noted that the shift factor obeys the TS2 (Two State, Two (Time) Scale (TS2)) function, which was shown to successfully describe the TTS of neat and nanocomposite PUas in the temperature range of \SI{-70}{\degreeCelsius} to \SI{70}{\degreeCelsius} \cite{tzelepis_experimental_2023, tzelepis_polyureagraphene_2023}. This model describes the glass transition as the transition between high-temperature and low-temperature Arrhenius region as follows:
\begin{linenomath}
\begin{equation} \label{eq-at}
\begin{aligned}
\ln(a_T) &\equiv \ln\left(\frac{\tau[T]}{\tau[T_o]}\right) \\
&= \frac{E_1}{RT} + \frac{(E_2 - E_1)}{RT} \left(\frac{1}{1 + \exp\left\{\frac{\Delta S}{R} \left(1 - \frac{T^*}{T}\right)\right\}}\right) -\\
&\quad \frac{E_1}{RT_o} + \frac{(E_2 - E_1)}{RT_o} \left(\frac{1}{1 + \exp\left\{\frac{\Delta S}{R} \left(1 - \frac{T^*}{T_o}\right)\right\}}\right),
\end{aligned}
\end{equation}
\end{linenomath}
where \(E_1\) and \(E_2\) are activation energies (\(\mathrm{J/mol}\)), \(\Delta S/R\) is the dimensionless transition entropy between the solid and liquid states of matter, \(T^*\) is the transition temperature (\(\mathrm{K}\)) (typically, \(T^* \approx T_g\)), and \(T_o\) is the reference temperature of the TTS shifts.

\subsection{Optimization Subject to a Morphological Constraint}
The parameters for both FMG-FMG and FMM-FMG models are fitted using a global particle-swarm optimization (PSO) algorithm \cite{kennedy_particle_1995}. Each optimization run is characterized by a consistent swarm population size of \(N_{pop}=200\) and total of \(N_{it}=6000\) iterations. To account for the stochastic nature of the PSO algorithm, 50 optimization runs are executed, followed by calculating the mean values and standard deviations for each material parameter.

In our PSO, a crucial constraint is integrated into the algorithm to determine the characteristic time of the second FMG branch. This constraint stems from the new dimensionless number (\(\mathcal{N}_P\)) we introduce for the first time for each branch in the models, which is obtained through Buckingham-Pi theorem, defined as:
\begin{equation}
\begin{aligned}
\mathcal{N}_P = \frac{L}{c\tau_c} = \left[\frac{\rho}{E_c}\right]^{1/2}\frac{L}{\tau_c},
\end{aligned}
\end{equation}
in which \(\tau_c\) represents the characteristic time, \(\rho\) is the density, \(L\) is the characteristic length (described below in more detail), \(E_c\) is the Young's modulus, and \(c\) is the speed of sound, defined as \cite{kinsler2000fundamentals},
\begin{equation} \label{eq-speedc}
\begin{aligned}
c \cong \sqrt{\frac{E}{\rho}}.
\end{aligned}
\end{equation}
In writing \eqref{eq-speedc}, we used a one-dimensional instead of a three-dimensional speed of sound equation, but the differences are minor and not relevant for our simple scaling analysis. As shown below, the dimensionless parameter \(\mathcal{N}_P\) is essentially the inverse of a \say{Deborah number}, \(De \cong \omega^* \tau_c\), with the \say{phononic bandgap} frequency \(\omega^* \cong c/L\). We further hypothesize that the parameter \(\mathcal{N}_P\) must be the same for both branches, leading to the following constraint on the characteristic times, 
\begin{equation} \label{eq-constraint}
\begin{aligned}
\tau_{c,2} = \tau_{c,1} \times \left(\frac{E_{c,1}}{E_{c,2}}\right)^{1/2}.
\end{aligned}
\end{equation}

This relationship between the moduli and relaxation times of the two elements can be further justified in the following. As mentioned earlier, the polyurea nanocomposites exhibit time-temperature superposition (TTS) in a broad range of temperatures (\(T \approx \SIrange{-70}{70}{\degreeCelsius}\)), frequencies (\(\omega \approx \SIrange{0.1}{100}{\per\second}\)) and formulations (HSWF \(\approx\) 20 to 40 wt.\% and GNP loading between 0 and 1.5 wt.\%). This was (pleasantly) surprising and somewhat unexpected, and the reasons for this apparent \say{rheological simplicity} are still uncertain. Even so, we can address at least some conditions necessary for this to happen.

\begin{figure}[H]
\includegraphics[width=0.22\linewidth, clip]{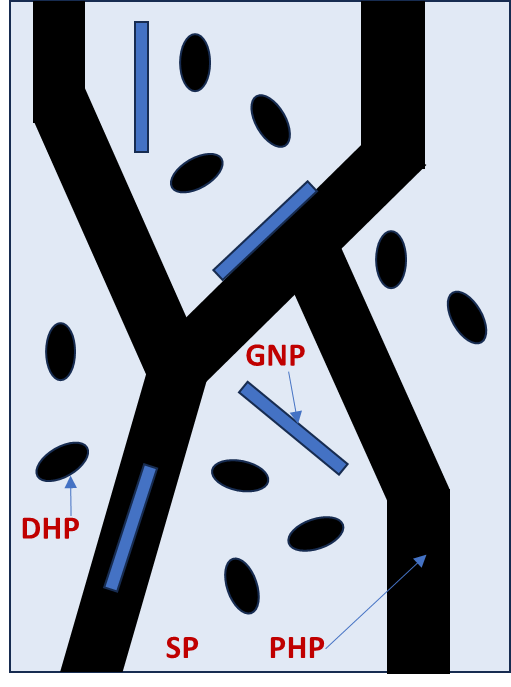}
\centering
\caption{Schematic representation of the polyurea-GnP nanocomposite morphology. Here, SP is soft phase, PHP is percolated hard phase, DHP is dispersed hard phase, and GNP is graphene nanoplatelet. The distance \(L\) is the ``characteristic size'' (on the order of 15-20 nm) that usually is seen as a peak in small-angle X-ray scattering experiments. \label{fig-thumbnail_image}}
\end{figure}

We hypothesize that one such condition is as follows. When analyzing the DMA results, we model the material as a simple viscoelastic Maxwell (or fractional Maxwell) element, where the whole material moves as a single viscoelastic body, with no internal vibrations. We stipulate that the co-continuous nature of the two-phase composite creates a \say{phononic bandgap} – the phonons with frequencies below \(\omega^*\) scatter at the multiple hard-soft interfaces and can be neglected when analyzing the material response to external forces. Stimuli with frequencies greater than \(\omega^*\) are likely to excite additional phonons with the wave number \(k > 0\), leading to the loss of the TTS property. This implies that the characteristic oscillation frequency, \(\omega\), must be smaller than the threshold frequency, \(\omega^*\), given by
\begin{equation}
\begin{aligned}
\omega^* = ck^* \cong \sqrt{\frac{E}{\rho}} \frac{2\pi}{L},
\end{aligned}
\end{equation}
where, \(E\) is the Young's modulus, \(\rho\) is the density, and \(L\) is the \say{characteristic size} (see Figure \ref{fig-thumbnail_image}). Let us now substitute \(\omega \cong \left(a_T(T)\tau_c\right)^{-1}\) and introduce indices \(i=1,2\) for elements 1 and 2,
\begin{equation} \label{eq-justify2}
\begin{aligned}
\left(a_T(T)\tau_{c,i}\right)^{-1} \leq \sqrt{\frac{E_{c,i}}{\rho}} \frac{2\pi}{L}.
\end{aligned}
\end{equation}
The shift factor, \(a_T(T)\) is a monotonically decreasing function of temperature. We can thus assume that the temperature at which the two sides of \eqref{eq-justify2} are equal is the upper limit for the validity of the TTS (it is of course possible that TTS could break down for other reasons at lower temperatures). If we label that temperature \(T_{max}\) and re-arrange, we can re-write \eqref{eq-justify2} in the following form,
\begin{equation}
\begin{aligned}
a_T(T_{max})\tau_{c,i} \sqrt{\frac{E_{c,i}}{\rho}} \frac{2\pi}{L} = 1.
\end{aligned}
\end{equation}
From this equation, and neglecting the small differences in the densities of the hard and soft phases, we obtain the condition described in \eqref{eq-constraint},
\begin{equation} \begin{aligned}
\tau_{c,1} \sqrt{E_{c,1}} \cong \tau_{c,2} \sqrt{E_{c,2}}
\end{aligned}
\end{equation}

Additionally, the subsequent parameter ranges are considered for all nanocomposite systems for both models. The characteristic moduli are limited to the ranges of \(0 \leq E_{c_1} \leq 10^4 \: \text{MPa}\), and \(0 \leq E_{c_2} \leq 10^3 \: \text{MPa}\). The characteristic time-scale of the first branch falls within the range of \(10^{-3} \: \text{s} \leq \tau_{c_1} \leq 10^2 \: \text{s}\), and the power exponents \(\alpha_{1,2}\) and \(\beta_{1}\) are allowed to vary from 0 to 1. Importantly, since the optimization is performed in the frequency domain, the constraint on the power exponents in the form of  \(0 \leq \beta \leq \alpha \leq 1\), as defined in the time-domain formulation of the fractional Maxwell model, are not applicable to \eqref{eq-EpSum} and \eqref{eq-EppSum}.

The objective of this optimization is to concurrently fit the model to both storage and loss moduli. To achieve this goal, a cost function is defined, comprising two terms, each representing one of the moduli as
\begin{equation}
\begin{aligned}
\min_{\theta} \left(w_1 f_1(\theta) + w_2 f_2(\theta)\right),
\end{aligned}
\end{equation}
where \(w_1\) and \(w_2\) are weights, each set to \(1/2\), and \(\theta\) represents the vector of fitting (model) parameters. The cost function for each modulus is computed as the squares of logarithmic difference between the experimental data and model predictions summed over all data points (\(N_d\))
\begin{equation}
\begin{aligned}
f_1(\theta) &= \sum_{i=1}^{N_d} \left(\log(\frac{E_{\text{exp}}^{'}}{E_{\text{model}}^{'}})\right)^2, 
f_2(\theta) &= \sum_{i=1}^{N_d} \left(\log(\frac{E_{\text{exp}}^{''}}{E_{\text{model}}^{''}})\right)^2.
\end{aligned}
\end{equation}
Additionally, the relative error, used to assess the fitting quality, is defined as
\begin{equation}
\text{Error} = \frac{\omega_1 f_1(\theta) + \omega_2 f_2(\theta)}{\omega_1 \sum_{i=1}^{N_d} (\log(E_{\text{exp}}'))^2 + \omega_2 \sum_{i=1}^{N_d} (\log(E_{\text{exp}}''))^2}.
\end{equation}
Both the models and PSO code were developed in MATLAB R2022b and executed in ICER MSU HPCC system with 1 node, 8 CPUs, and 64GB RAM.

\subsection{Sensitivity Analyses (SA)} \label{sec.SA}
The constitutive models discussed thus far have not explicitly accounted for the variability in the outputs attributable to the uncertainties in the input model parameters. These uncertainties stem from the errors inherent in the experimental datasets and the employed models, causing the deterministic model parameters to become random variables and models to render stochastic. Sensitivity analysis can provide valuable insights into the effect of input parameters on the model outputs.

Sensitivity analysis is defined as \say{the study of how uncertainty in the output of a model can be apportioned to different sources of uncertainty in the model input}, \cite{saltelli_sensitivity_2002}. One of the primary goals of SA is \say{factor prioritization}, aiming to quantify the contribution of the model inputs to the uncertainty in the model outputs and identify the most influential factors. Parameters exerting minimal influence on the model output can then be assigned reasonable \say{deterministic} values, which marks another goal of SA, known as factor fixing, \cite{saltelli_global_2008}. To achieve these objectives, performing both local and global sensitivity analyses are indispensable.

\subsubsection{Local Sensitivity Analysis (LSA)} \label{sec.LSA}
A partial derivative serves as a local sensitivity measure, which can be expressed for a model with a generic form of \(y = f(\mathbf{q})\) as \(\frac{\partial y}{\partial q_i}\vert_{\mathbf{q}^0} = \frac{\partial f(\mathbf{q}^0)}{\partial q_i} \), where \(y\) denotes the scalar model output, \(\mathbf{q}\) represent the vector of model parameters, and \(q_i\) is the i-th model parameter with a realized value of \(q^0_i\). This local sensitivity measure is referred to as the local sensitivity index, denoted by \(S_i\). To ensure the comparability of these indices for parameters with different orders of magnitude and to eliminate the presence of physical units in their definition, they can be normalized as follows \cite{zayernouri2011coherent, norton_introduction_2015}:
\begin{equation}
\begin{aligned}
\bar{S}_i = \frac{q^0_i}{f(\mathbf{q}^0)} \cdot \frac{\partial f(\mathbf{q}^0)}{\partial q_i}.
\end{aligned}
\end{equation}
This normalized form of a local sensitivity index is also known as the \textit{elasticity of \(y\) with respect to \(q_i\)}, denoted by \(E_i(q_i)\), \cite{borgonovo_sensitivity_2017}. Alternative normalization methods for differentiation-based sensitivity measures can be found in \cite{saltelli_global_2008, borgonovo_new_2001}.

In our current analysis, given the explicit formulation of our models (refer to \eqref{eq-EpSum} and \eqref{eq-EppSum}), computing the LS indices for both the storage and loss moduli across decades of shifted frequency \((x=a_T\omega)\) is straightforward utilizing the MATLAB symbolic toolbox. These indices are defined as follows:

\begin{subequations}
\begin{equation}  \label{eq-SEp}
\bar{S}_{E^{\prime}, q_i} = \frac{q_i^0}{E^{\prime}(\mathbf{q}^0; x)} \cdot \frac{\partial E^{\prime}(\mathbf{q}^0; x)}{\partial q_i},
\end{equation}

\begin{equation}  \label{eq-SEpp}
\bar{S}_{E^{\prime\prime}, q_i} = \frac{q_i^0}{E^{\prime\prime}(\mathbf{q}^0; x)} \cdot \frac{\partial E^{\prime\prime}(\mathbf{q}^0; x)}{\partial q_i}.
\end{equation}
\end{subequations}
Up to now, a single set of model parameters has been responsible for two scalar model outputs (i.e., storage and loss moduli). To simplify the factor prioritization process, it is advantageous to focus on a single scalar model output that incorporates the combined effect of model parameters on both the storage and loss moduli simultaneously. This approach will facilitate the identification of the least and most influential model parameters. This objective can be achieved by calculating the magnitude of the normalized local sensitivity indices of the complex modulus. This index is defined as
\begin{linenomath}
\begin{equation} \label{eq-SEstar}
\begin{aligned}
\left|\bar{S}_{E^*, q_i}\right| &= \frac{q_i^0}{\left|E^*(\mathbf{q}^0; x)\right|} \cdot \left|\frac{\partial E^*(\mathbf{q}^0; x)}{\partial q_i}\right| \\
&= \frac{q_i^0}{\left|E^*(\mathbf{q}^0; x)\right|} \cdot \sqrt{\left(\frac{\partial E^{\prime}(\mathbf{q}^0; x)}{\partial q_i}\right)^2 + \left(\frac{\partial E^{\prime\prime}(\mathbf{q}^0; x)}{\partial q_i}\right)^2}.
\end{aligned}
\end{equation}
\end{linenomath}

Moreover, the range of variability for each model parameter is constructed based on the following assumptions. First, the baseline (mean) value (\(\mu_i\)) for each model parameter is set to its optimized value obtained from the PSO. Second, the corresponding standard deviation (\(\sigma_i\)) is assumed to be 5\% of the mean value. Lastly, parameters are distributed uniformly over a range with \(\mu_i=\frac{b_i + a_i}{2}\) and \(\sigma_i=\frac{b_i-a_i}{\sqrt{12}}\), where \(a_i\) and \(b_i\) represent the lower and upper limits, respectively.

With the explicit form of all LS indices at hand, they can be efficiently evaluated through Monte Carlo (MC) simulation \cite{keramati2024monte}. In each realization of the MC simulation, all model parameters are uniformly sampled from their respective range of variability \(\left(q_i \sim \mathcal{U}(a_i, b_i)\right)\), and each LS index is computed at all shifted frequency points. Subsequently, the mean value and standard deviation of LS indices are determined across the entire shifted frequency domain. Since LS indices are defined by means of partial derivatives, they provide information localized around the baseline or \say{nominal} value of only one model parameter at a time, leaving the rest of the parameter variability space unexplored, \cite{saltelli_how_2010}.

\subsubsection{Global Sensitivity Analysis (GSA)} \label{sec.GSA}
Global sensitivity analysis aims to comprehensively explore the variability of model outputs by systematically altering all uncertain input model parameters within their entire range of variability, \cite{saltelli_sensitivity_2002}. GSA methods fall into statistical-based, derivative-based, and variance-based categories, \cite{douglas-smith_certain_2020}, with the latter being widely recognized as the most standard approach due to its lack of constraints on the model forms and ease of interpretation, \cite{puy_comprehensive_2022}. Among the variance-based approaches, the Sobol' method stands out as the most popular implementation of GSA, \cite{razavi_what_2015}, where Sobol' proposed indices capable of measuring the first, total, and in general any order effect of model parameters on the model output, \cite{sobol_sensitivity_1993, sobol_global_2001}. This method is also utilized in the present work, and the definition of the first- and total-order indices comes next.

To maintain consistency with the notation used in the LSA section, we consider a similar model with a generic form \(y = f(\mathbf{q})\), where the scalar model output, \(y\), may represent the storage, loss, or magnitude of the complex modulus. The model parameters are treated as random variables with the realized values \(q_1^0, ..., q_i^0, ... q_n^0\) uniformly distributed over their respective ranges \((a_i, b_i)\) as determined previously. A first-order sensitivity index is defined as
\begin{equation}
\begin{aligned}
S_i = \frac{\mathbb{V}_{q_i} \left( \mathbb{E}_{\mathbf{q}_{\sim i}}\left(y|q_i\right) \right)}{\ \mathbb{V}\left(y\right)},
\end{aligned}
\end{equation}
where \(\mathbb{V}\left(y\right)\) represents the total (unconditional) variance of the output, \(\mathbb{E}_{\mathbf{q}_{\sim i}}\left(\cdot\right)\) denotes the conditional expectation of its argument taken over all factors except \(q_i\), and \(\mathbb{V}_{q_i}\left(\cdot\right)\) is the variance of its argument taken over \(q_i\). Notably, the notation \(\mathbf{q}_{\sim i}\) is equivalent to a vector of yellow\(\{q_1, ..., q_{i-1}, q_{i+1}, ..., q_{k}\}\) random variables. Further, a more complete form of notation for the conditional expectation is \(\mathbb{E}_{\mathbf{q}_{\sim i}}\left(y|q_i=q_i^0\right)\), where the random i-th model parameter is realized and sampled as \(q_i^0\) from its respective range of variability. According to a well-known theorem \cite{mood_introduction_1950},
\begin{equation} \label{eq-identitiy}
\begin{aligned}
\mathbb{V}_{q_i} \left(\mathbb{E}_{\mathbf{q}_{\sim i}} (y|q_i) \right) + \mathbb{E}_{q_i} \left(\mathbb{V}_{\mathbf{q}_{\sim i}} (y|q_i) \right) = \mathbb{V}(y),
\end{aligned}
\end{equation}
\(S_i\) is a normalized index that ranges between zero and one, solely measuring the main effect of factor \(q_i\) on the variability of the model output, \cite{saltelli_global_2008}. This index quantifies how much the variance of the model output would decrease if the model parameter \(q_i\) were fixed, \cite{saltelli_variance_2010}.

To capture the \say{total effect} of factor \(q_i\), which includes its first-order effect alongside its interactions with other model parameters \cite{saltelli_global_2008}, another sensitivity measure known as the total-order sensitivity index is commonly employed. It is defined as \cite{saltelli_variance_2010,sobol_global_2001}:
\begin{equation}
\begin{aligned}
S_{Ti} = \frac{ \mathbb{E}_{\mathbf{q}_{\sim i}} \left( \mathbb{V}_{q_{i}} (y|\mathbf{q}_{\sim i}) \right) }{\mathbb{V}\left(y\right)} = 1 - \frac{\mathbb{V}_{\mathbf{q}_{\sim i}} \left( \mathbb{E}_{q_{i}} (y|\mathbf{q}_{\sim i}) \right)}{\mathbb{V}(y)},
\end{aligned}
\end{equation}
where \(\mathbb{E}_{q_i}\left(\cdot\right)\) represents the conditional expectation of its argument taken over \(q_i\), and \(\mathbb{V}_{\mathbf{q}_{\sim i}}\left(\cdot\right)\) is the conditional variance of its argument taken over all factors except \(q_i\). \(S_{Ti}\) quantifies the reduction in the variance of the model output when all model parameters except the input factor \(q_i\) are held constant \cite{saltelli_variance_2010}.

Due to the high computational cost of obtaining sensitivity indices of all possible higher order effects for each input factor, the common practice is to compute only the first- and total-order sensitivity indices associated with each model parameter \cite{saltelli_global_2008, homma_importance_1996}. A more detailed derivation of these two indices through the Analysis Of Variances (ANOVA) \cite{yang2012adaptive, zhang2014enabling} is provided in Appendix \ref{appendix-A}.

In order to further reduce the computational cost of calculating the first- and total-order indices, the conditional expectations and variances are not computed directly; instead, various estimators have been developed for the precise yet computationally effective evaluation of these indices. For simultaneous computation of Sobol's indices, the methodology proposed by Salteli et al. \cite{saltelli_variance_2010, herman_salib_2017} is implemented (See Appendix \ref{appendix-B}).

\section{Experimental and Modeling Results} \label{sec.Exp&Mod}
\subsection{Modeling}
We focus exclusively on the results corresponding the FMM-FMG model. For a comprehensive review of the FMG-FMG model, readers are referred to our earlier work \cite{tzelepis_polyureagraphene_2023} where the results are provided and discussed in detail.

\begin{figure}[H]
    \centering
    \makebox[\textwidth][c]{
        \begin{minipage}{1\textwidth}
            \begin{subfigure}{0.5\textwidth}
                \centering
                \includegraphics[width=\linewidth, trim=0.2cm 0.4cm 4cm 1.5cm, clip]{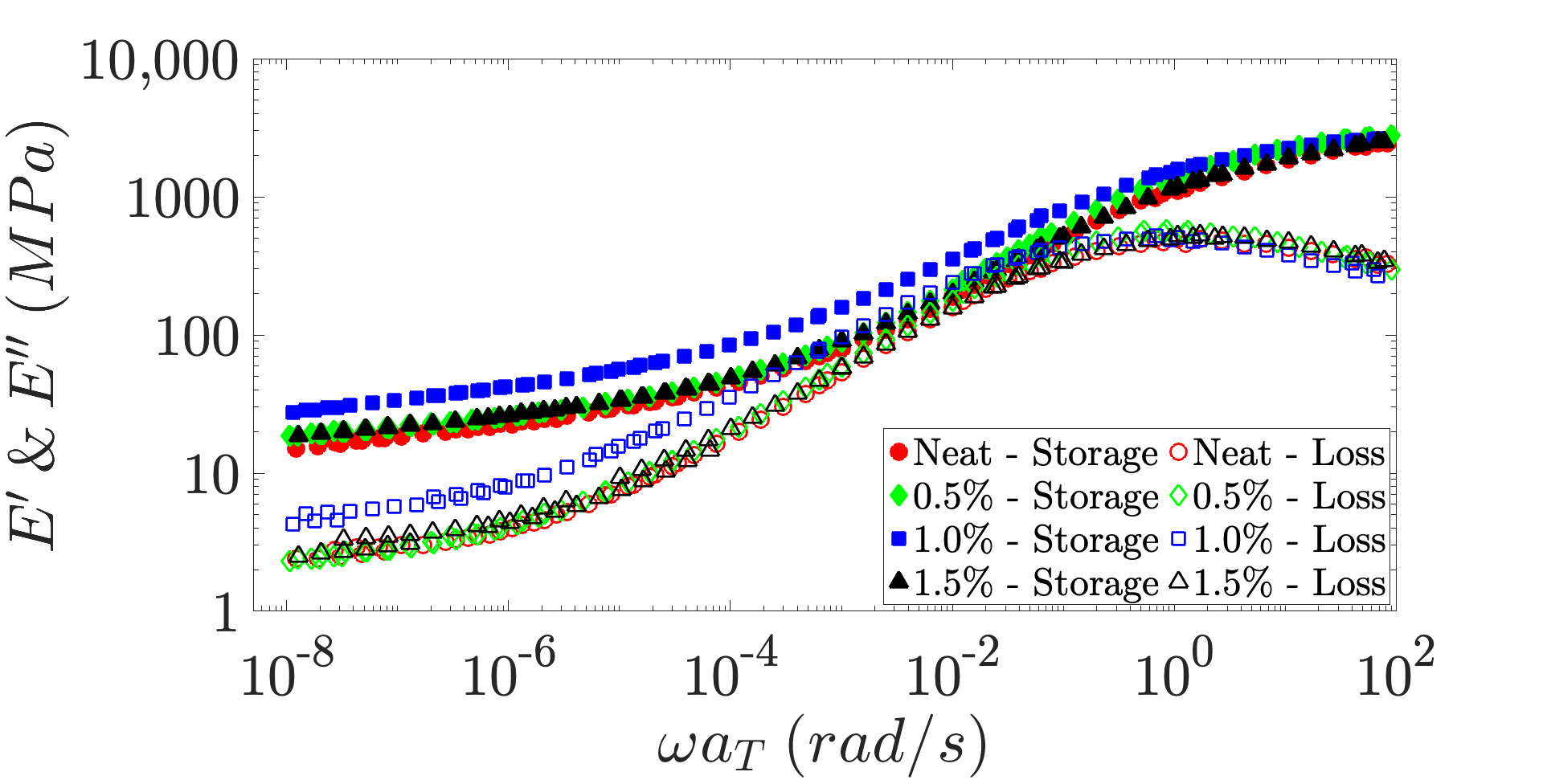}
                \caption{\label{fig-mastercurves-20HS}}
            \end{subfigure}
            \begin{subfigure}{0.5\textwidth}
                \centering
                \includegraphics[width=\linewidth, trim=0.2cm 0.4cm 4cm 1.5cm, clip]{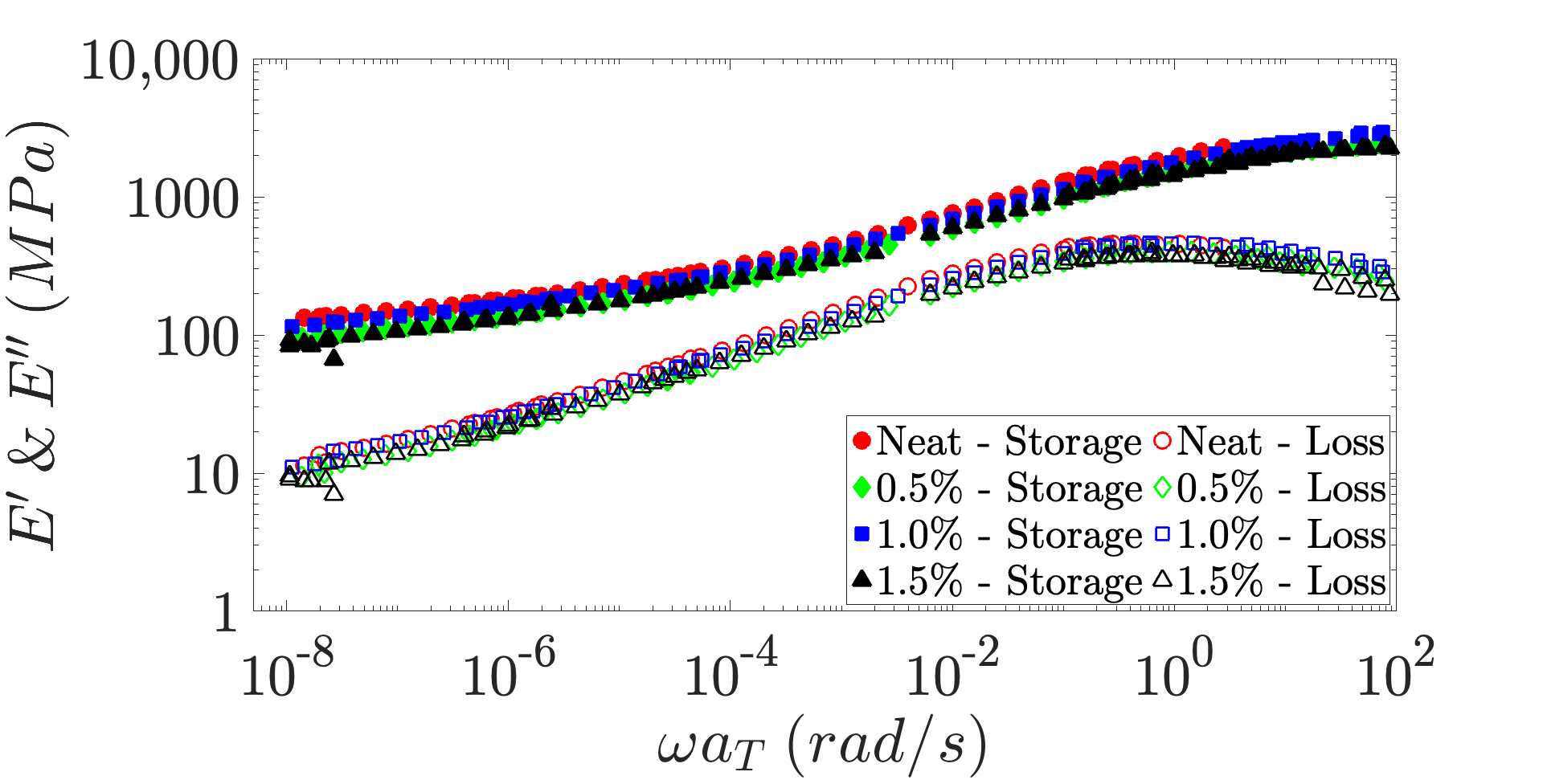}
                \caption{\label{fig-mastercurves-30HS}}
            \end{subfigure}
        \end{minipage}
    }
    \begin{subfigure}{0.5\textwidth}
        \centering
        \includegraphics[width=\linewidth, trim=0.2cm 0.4cm 4cm 1.5cm, clip]{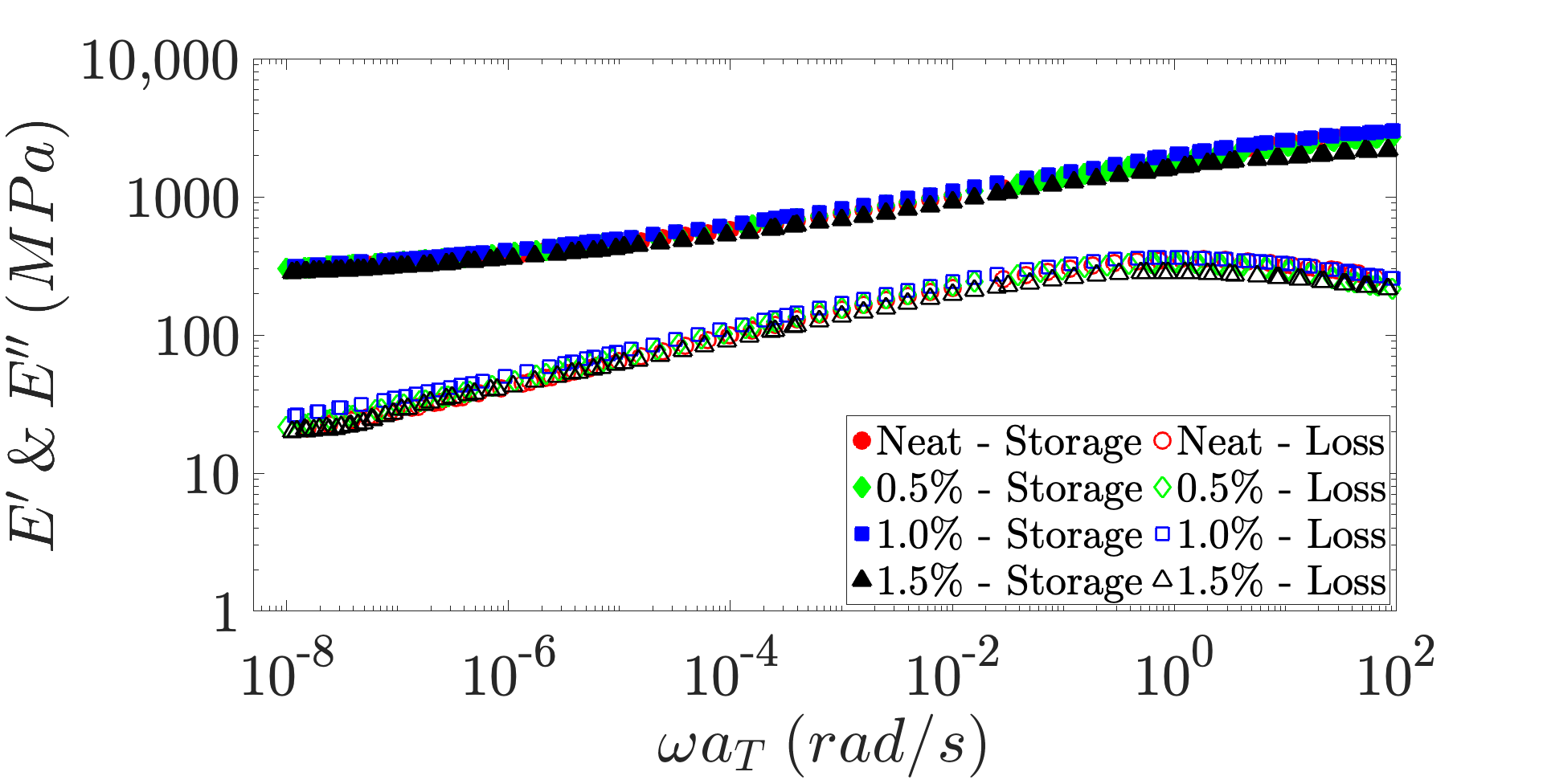}
        \caption{\label{fig-mastercurves-40HS}}
    \end{subfigure}
    \caption{Master curves for the tensile storage and loss moduli. (a) 20\% HSWF matrix with 0, 0.5, 1.0, and 1.5 wt.\% xGnP. (b) Same as (a) for the 30\% HSWF matrix. (c) Same as (a) for the 40\% HSWF matrix. \label{fig-Exp-mastercurves}}
\end{figure}

The storage and loss master curves are plotted for all nanocomposite systems: (a) IPDI-2k-20HS matrix, (b) IPDI-2k-30HS matrix; and (c) IPDI-2k-40HS matrix, in Figure \ref{fig-Exp-mastercurves}. In each system with identical hard segment volume fraction, all master curves exhibit a very close consistency, with a possible exception observed in the IPDI-2k-20HS, 1\% xGnP system (indicated by blue symbols in Figure \ref{fig-mastercurves-20HS}). The possible reasons for this unique behavior are discussed in further detail in \cite{tzelepis_polyureagraphene_2023}.

\begin{figure}[H]
\includegraphics[width=\linewidth, trim=3.25cm 1.1cm 5.5cm 2.1cm, clip]{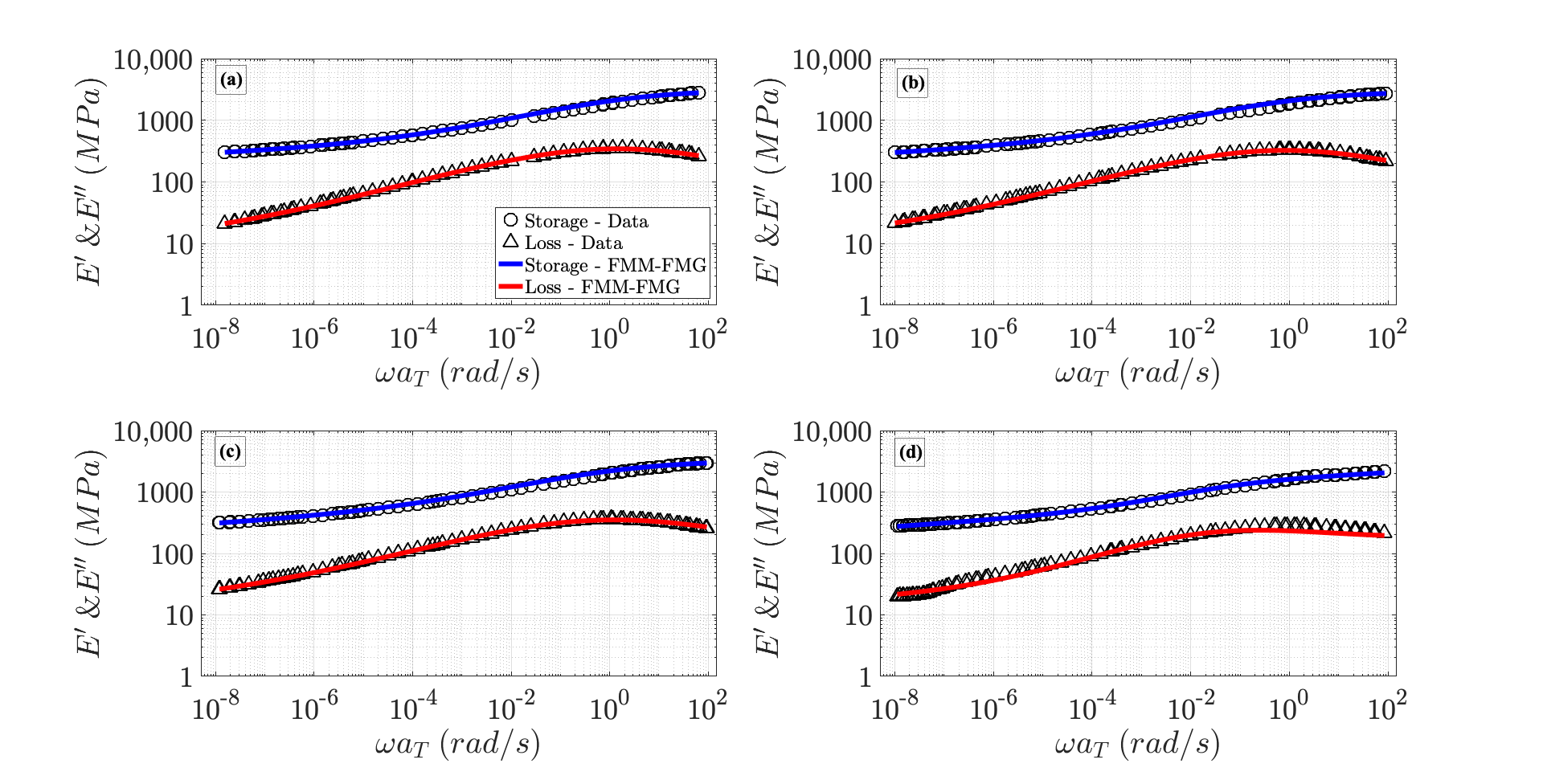}
\centering
\caption{Experimental and FMM-FMG master curves fitted for 40\% HSWF nanocomposite systems with (a) no added nanofillers; (b) 0.5 wt.\% xGnP; (c) 1.0 wt.\% xGnP; (d) 1.5 wt.\% xGnP. Circles and triangles represent the storage and loss modulus, respectively; Blue and red solid lines are the storage and loss moduli model fits, respectively. \label{fig-FMM-FMG-40HS-fit}}
\end{figure}

Moving forward, we present the results of the FMM-FMG model fitting to the master curves next. Table \ref{tab-opt} provides the mean values and standard deviation for all seven FMM-FMG model parameters. The optimization runs demonstrated excellent convergence and reproducibility, as evidenced by the low standard deviation values, with a maximum order of magnitude of \(\mathcal{O}(10^{-4})\), across all systems considered. Figure \ref{fig-FMM-FMG-40HS-fit} illustrates the FMM-FMG model fits to the experimental shifted data for 40HSWF nanocomposites, while the results for 20HSWF and 30HSWF nanocomposites are provided in Figures S1 and S2, respectively. All the fitted curves were generated using the mean values for the model parameters (Table \ref{tab-opt}), given that the standard deviation of each model parameter was negligible.

For all systems, the relative error between the model and data was less than 1.98\%, covering a broad range of frequencies (between \(10^{-8}\) and \(10^2\) \(\text{ rad s}^{-1}\)). In contrast to our previous work, in which the FMG-FMG model was used, no deviation between the present FMM-FMG model and the experimental data was observed above the glass transition point in the loss modulus for the 20HSWF sample (see Figure S3). This improvement can be attributed to the additional power exponent included in the FMM branch (\(\beta_1\)), which better captures the viscoelastic behavior of the soft phase in the samples.

\begin{figure}[H]
\includegraphics[width=\linewidth, trim=5.25cm 1.25cm 5.5cm 2.2cm, clip]{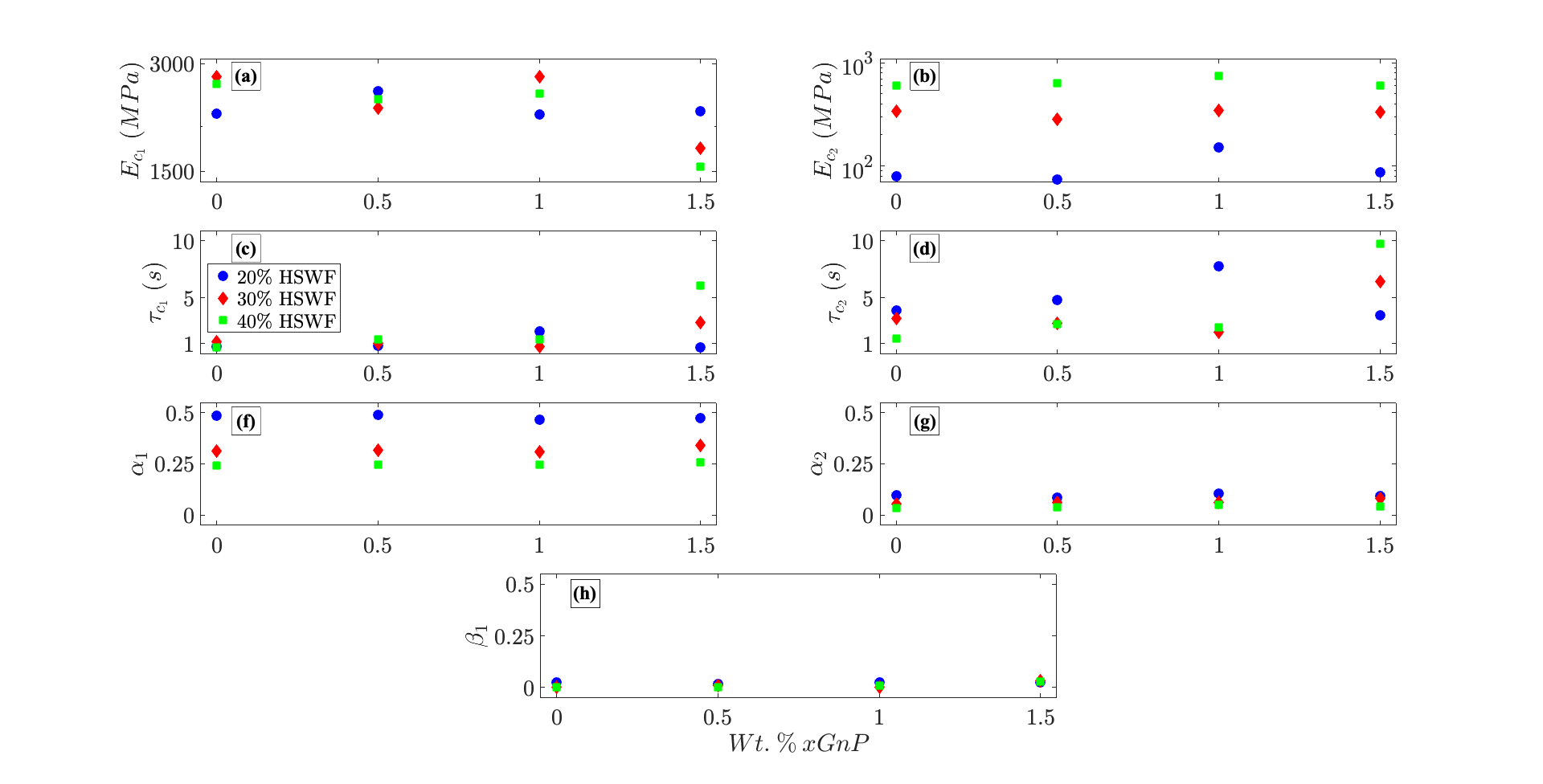}
\centering
\caption{The effect of the nanofiller loading on the (a) characteristic modulus of the first branch (\(E_{c_1}\)) (b) characteristic modulus of the second branch (\(E_{c_2}\)), (c) characteristic time of the first branch (\(\tau_{c_1}\)), (d) characteristic time of the second branch (\(\tau_{c_2}\)), (e) first power law exponent of the first branch (\(\alpha_{1}\)), (f) first power law exponent of the second branch (\(\alpha_{2}\)), and (g) second power law exponent of the second branch (\(\beta_{1}\)). \label{fig-Exp-params}}
\end{figure}

Figure \ref{fig-Exp-params} illustrates the influence of the nanofiller content on the variation of the optimized model parameters across all families of nanocomposites. Figure \ref{fig-Exp-params}a depicts the minor effects of nanoparticles weight fraction on the characteristic modulus of the FMM branch (\(E_{c_1}\)). Figure \ref{fig-Exp-params}b shows an increase in the characteristic modulus of the FMG branch (\(E_{c_2}\)) as the xGnP content increases up to 1 wt.\%, followed by a decrease when the xGnP further increases to 1.5 wt.\% for all nanocomposite systems. This trend can be linked to the fact that the impact of nanofillers typically reaches a peak before diminishing. This reduction is likely due to the aggregation of platelets, which decreases their aspect ratio and reduces their overall effect.

Furthermore, time-scales of both branches stay relatively unchanged by the addition of the nanoparticles as shown in Figure \ref{fig-Exp-params}c and \ref{fig-Exp-params}d. In Figure \ref{fig-Exp-params}e, \(\alpha_1\) shows no statistically meaningful dependence on the xGnP weight fractions but exhibits a strong dependence on hard segment content, decreasing as HSWF increases, which is in line with the fact that the material becomes \say{more elastic}. Additionally, Figure \ref{fig-Exp-params}f reveals that the power exponent of hard phase, \(\alpha_2\), remains close to zero in all nanocomposite systems, never exceeding a maximum value of 0.11. This indicates that the behavior of this branch is very similar to an almost pure elastic spring. Finally, Figure \ref{fig-Exp-params}g shows the variation of the \(\beta_1\), the second power exponent of the soft phase. This parameter is either zero (as in the 30HSWF 0.0xGnP, 40HSWF 0.0xGnP and 0.5xGnP systems) or very close to zero in the remaining nanocomposite systems with a maximum value of 0.033 (Table \ref{tab-opt}). This also suggests that the second spring-pot element associated with the soft phase acts similarly to a pure elastic spring, although with a slight inclusion of viscoelastic behavior.


\begin{table}[H]
    \centering
    \caption{Optimized FMM-FMG model parameters. Rows represent nanocomposite systems with 20\%, 30\%, and 40\% HSWF, while columns correspond to the model parameters.}
    \label{tab-opt}
        \begin{adjustbox}{max width=1.1\textwidth,center=\textwidth}
        \begin{tabular}{ccccccccc}
            \toprule
            \multicolumn{2}{c}{Samples} & \(E_{c_1}\) & \(\tau_{c_1}\) & \(\alpha_{1}\) & \(\beta_{1}\) & \(E_{c_2}\) & \(\tau_{c_2}\) & \(\alpha_{2}\)\\
            \midrule
            \multirow{8}{*}{\centering 20HS} & \multirow{2}{*}{0.0xGnP} & 2180 & 0.75 & 0.48 & 0.023 & 79 & 3.92 & 0.10 \\
                                                            & & {$\pm1.0\times10^{-5}$} & {$\pm1.6\times10^{-8}$} & {$\pm1.3\times10^{-9}$} & {$\pm3.4\times10^{-10}$} & {$\pm4.6\times10^{-7}$} & {$\pm6.8\times10^{-8}$} & {$\pm4.9\times10^{-10}$} \vspace{8pt}\\
                                  & \multirow{2}{*}{0.5xGnP} & 2513 & 0.82 & 0.49 & 0.018 & 74 & 4.79 & 0.08 \\
                                                            & & {$\pm5.9\times10^{-6}$} & {$\pm1.1\times10^{-8}$} & {$\pm1.1\times10^{-9}$} & {$\pm1.6\times10^{-10}$} & {$\pm4.6\times10^{-7}$} & {$\pm4.6\times10^{-8}$} & {$\pm5.5\times10^{-10}$} \vspace{8pt}\\
                                  & \multirow{2}{*}{1.0xGnP} & 2166 & 2.04 & 0.46 & 0.025 & 150 & 7.74 & 0.11 \\
                                                            & & {$\pm8.8\times10^{-6}$} & {$\pm3.8\times10^{-8}$} & {$\pm1.2\times10^{-9}$} & {$\pm2.8\times10^{-10}$} & {$\pm7.9\times10^{-7}$} & {$\pm1.2\times10^{-7}$} & {$\pm4.6\times10^{-10}$} \vspace{8pt}\\ 
                                  & \multirow{2}{*}{1.5xGnP} & 2211 & 0.69 & 0.47 & 0.025 & 87 & 3.49 & 0.09 \\
                                                            & & {$\pm1.1\times10^{-5}$} & {$\pm1.6\times10^{-8}$} & {$\pm1.4\times10^{-9}$} & {$\pm4.3\times10^{-10}$} & {$\pm5.5\times10^{-7}$} & {$\pm6.5\times10^{-8}$} & {$\pm5.7\times10^{-10}$} \vspace{8pt}\\
            \midrule
            \multirow{8}{*}{\centering 30HS} & \multirow{2}{*}{0.0xGnP} & 2758 & 1.14 & 0.31 & 0.0 & 342 & 3.92 & 0.05 \\
                                                            & & {$\pm2.2\times10^{-5}$} & {$\pm5.0\times10^{-8}$} & {$\pm1.3\times10^{-9}$} & {0} & {$\pm2.0\times10^{-6}$} & {$\pm1.2\times10^{-7}$} & {$\pm4.4\times10^{-10}$} \vspace{8pt}\\
                                  & \multirow{2}{*}{0.5xGnP} & 2251 & 0.98 & 0.32 & 0.009 & 285 & 2.76 & 0.06 \\
                                                            & & {$\pm2.2\times10^{-5}$} & {$\pm5.4\times10^{-8}$} & {$\pm1.6\times10^{-9}$} & {$\pm1.2\times10^{-9}$} & {$\pm2.2\times10^{-6}$} & {$\pm1.3\times10^{-7}$} & {$\pm5.9\times10^{-10}$} \vspace{8pt}\\
                                  & \multirow{2}{*}{1.0xGnP} & 2758 & 0.71 & 0.31 & 0.002 & 343 & 2.02 & 0.06 \\
                                                            & & {$\pm5.5\times10^{-5}$} & {$\pm8.1\times10^{-8}$} & {$\pm3.3\times10^{-9}$} & {$\pm2.3\times10^{-9}$} & {$\pm5.4\times10^{-7}$} & {$\pm2.0\times10^{-7}$} & {$\pm1.1\times10^{-9}$} \vspace{8pt}\\
                                  & \multirow{2}{*}{1.5xGnP} & 1741 & 2.83 & 0.34 & 0.033 & 334 & 6.46 & 0.08 \\
                                                            & & {$\pm1.0\times10^{-4}$} & {$\pm7.6\times10^{-7}$} & {$\pm6.9\times10^{-9}$} & {$\pm7.1\times10^{-9}$} & {$\pm5.9\times10^{-6}$} & {$\pm1.5\times10^{-6}$} & {$\pm2.1\times10^{-9}$} \vspace{8pt}\\
            \midrule
            \multirow{8}{*}{\centering 40HS} & \multirow{2}{*}{0.0xGnP} & 2636 & 0.69 & 0.24 & 0.0 & 604	& 1.44 & 0.03 \\
                                                            & & {$\pm2.6\times10^{-6}$} & {$\pm2.7\times10^{-9}$} & {$\pm3.4\times10^{-10}$} & {0} & {$\pm1.8\times10^{-6}$} & {$\pm6.2\times10^{-9}$} & {$\pm1.7\times10^{-10}$} \vspace{8pt}\\
                                  & \multirow{2}{*}{0.5xGnP} &2389	&1.39 & 0.25 & 0.0 & 638 & 2.68 & 0.04 \\
                                                            & & {$\pm3.8\times10^{-5}$} & {$\pm1.6\times10^{-7}$} & {$\pm3.7\times10^{-9}$} & {$\pm1.0\times10^{-9}$} & {$\pm9.5\times10^{-6}$} & {$\pm2.7\times10^{-7}$} & {$\pm9.1\times10^{-10}$} \vspace{8pt}\\
                                  & \multirow{2}{*}{1.0xGnP} & 2475 & 1.35 & 0.25 & 0.01 & 748 & 2.45 & 0.05\\
                                                            & & {$\pm5.9\times10^{-5}$} & {$\pm2.2\times10^{-7}$} & {$\pm4.8\times10^{-9}$} & {$\pm1.5\times10^{-9}$} & {$\pm1.4\times10^{-5}$} & {$\pm3.4\times10^{-7}$} & {$\pm1.2\times10^{-9}$} \vspace{8pt}\\
                                  & \multirow{2}{*}{1.5xGnP} & 1545 & 6.09 & 0.26 & 0.029 & 602 & 9.75 & 0.04 \\
                                                            & & {$\pm4.4\times10^{-5}$} & {$\pm8.3\times10^{-7}$} & {$\pm5.5\times10^{-9}$} & {$\pm1.4\times10^{-9}$} & {$\pm1.5\times10^{-5}$} & {$\pm1.1\times10^{-6}$} & {$\pm1.7\times10^{-9}$} \vspace{8pt}\\
            \bottomrule
        \end{tabular}
        \end{adjustbox}
\end{table}

\subsection{Local Sensitivity Analysis}
\subsubsection{Factor Prioritization and Factor Fixing based on the LSA for the FMG-FMG Model Parameters}
\begin{figure}[H]
    \centering
    \begin{subfigure}{0.45\textwidth}
        \includegraphics[width=\linewidth, trim=3.5cm 0.5cm 5.5cm 1.5cm, clip]{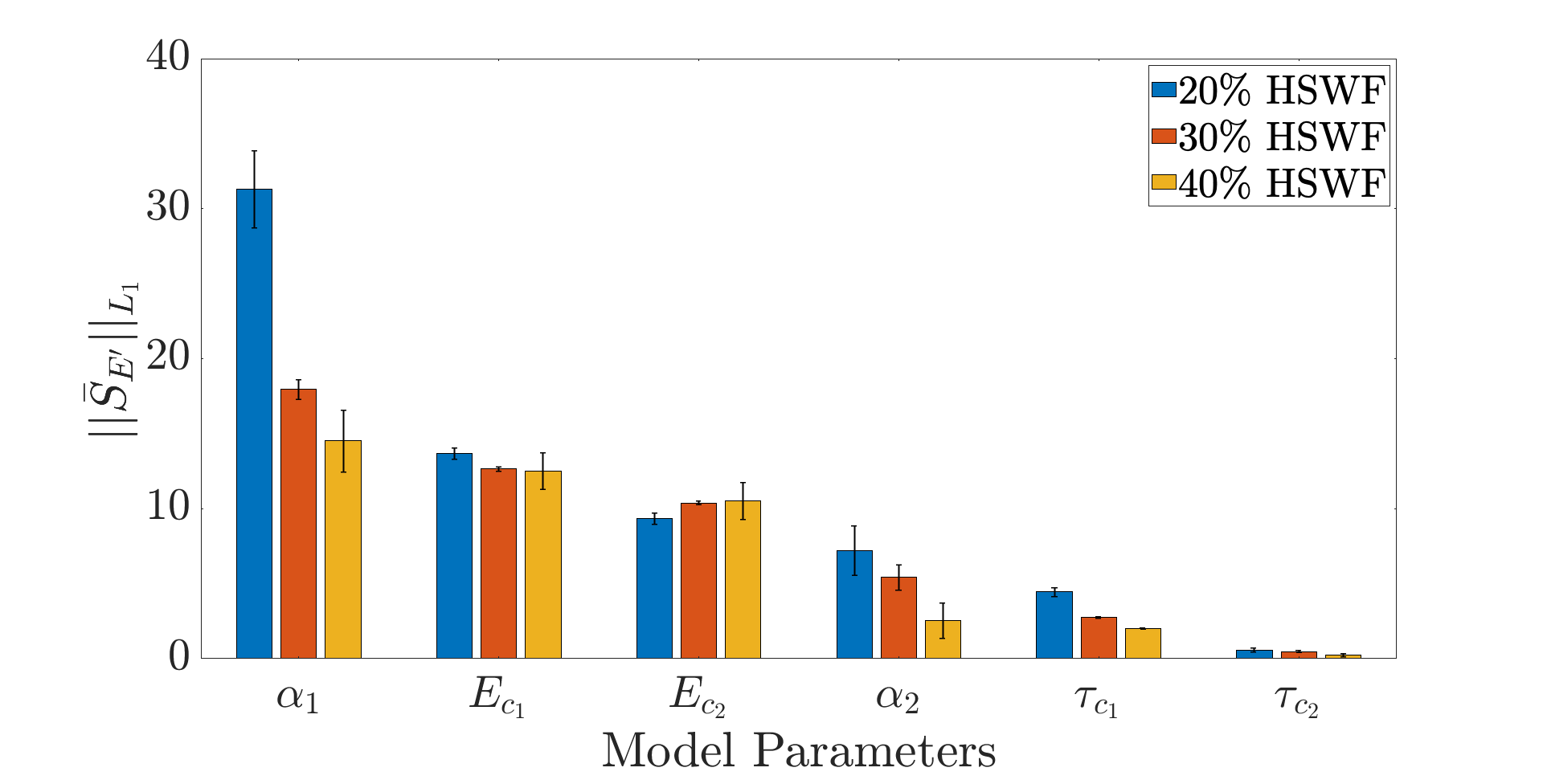}
        \caption{\label{fig-barplots-SEp-L1-FMG}}
    \end{subfigure}
    \begin{subfigure}{0.45\textwidth}
        \includegraphics[width=\linewidth, trim=3.5cm 0.5cm 5.5cm 1.5cm, clip]{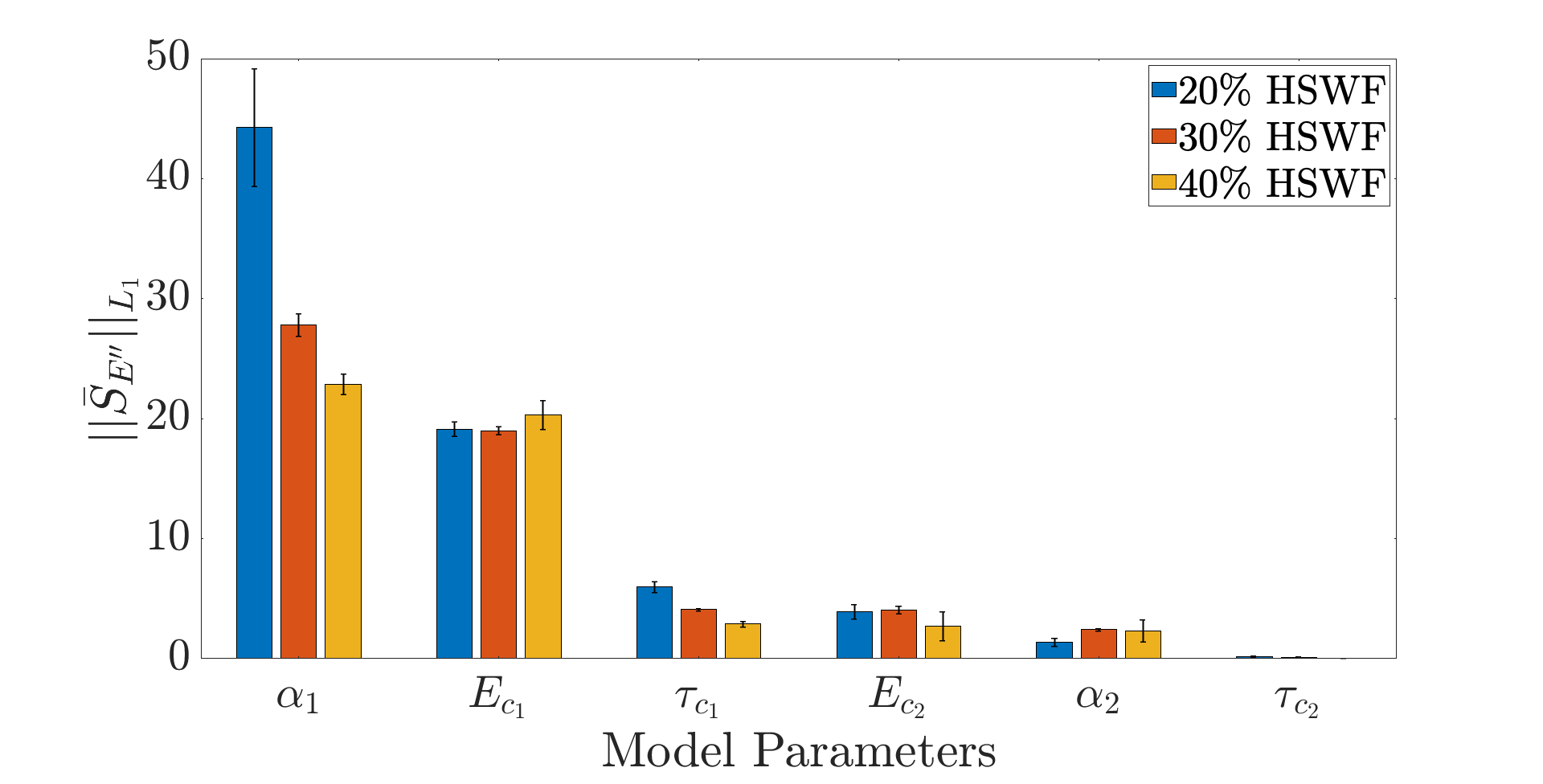}
        \caption{\label{fig-barplots-SEpp-L1-FMG}}
    \end{subfigure}
    \\
    \begin{subfigure}{0.45\textwidth}
        \includegraphics[width=\linewidth, trim=3.5cm 0.5cm 5.5cm 1.5cm, clip]{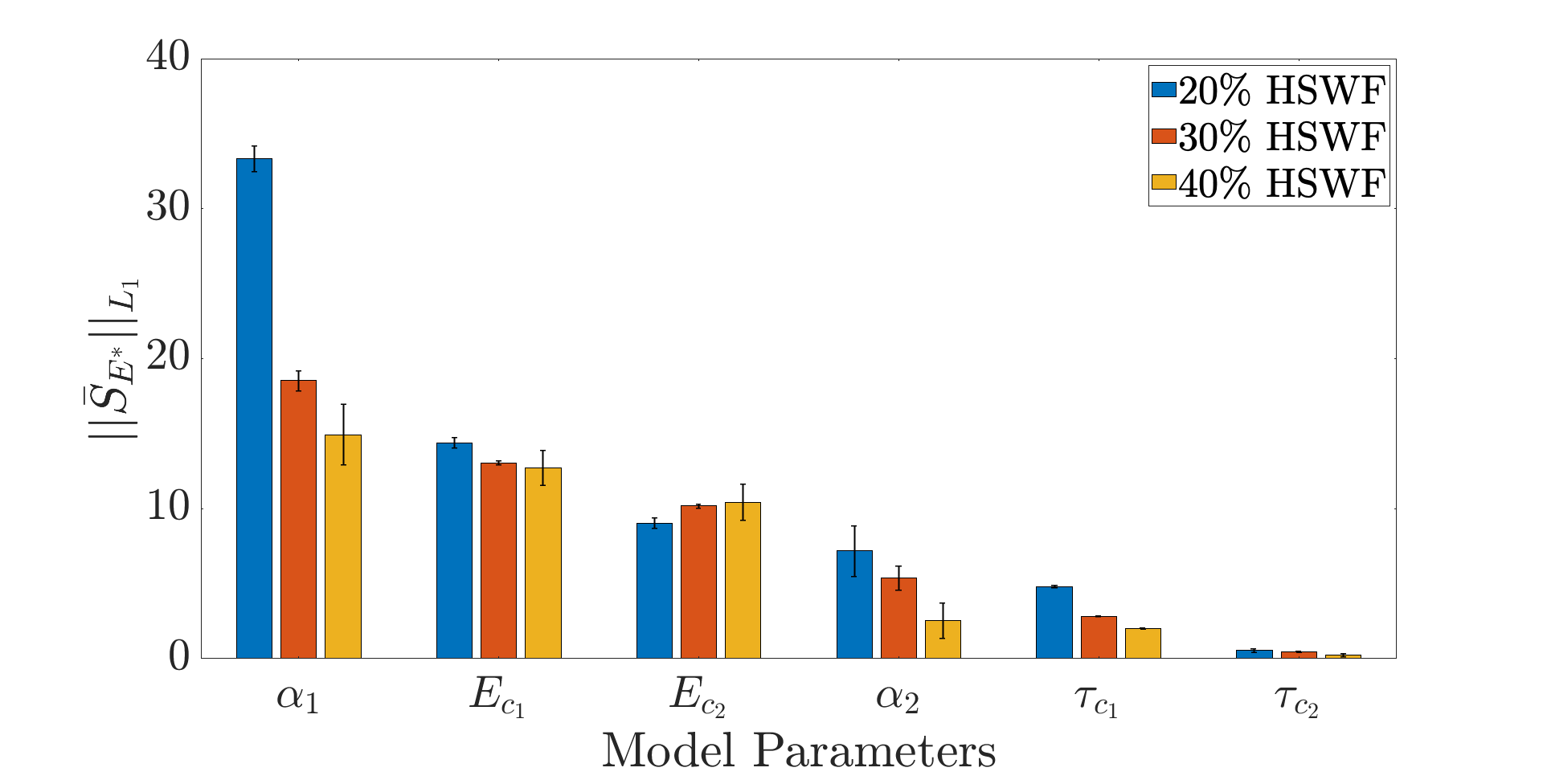}
        \caption{\label{fig-barplots-SEstar-L1-FMG}}
    \end{subfigure}
    \caption{Averaged L1-norm of the normalized local sensitivity indices associated with the (a) storage, (b) loss, (c) complex moduli with respect to each FMG-FMG model parameter.\label{fig-barplots-L1-FMG}}
\end{figure}

According to the definition of normalized LS indices, \eqref{eq-SEp} and \eqref{eq-SEpp}, they can exhibit various trends and distributions across at least ten decades of frequency (e.g., see Figures \ref{fig-LSA-Ep-FMG-20HS}, \ref{fig-LSA-Epp-FMG-20HS}, and \ref{fig-LSA-Estar-FMG-20HS}). Therefore, adoption of a reliable quantitative measure becomes crucial to take into account these variations and to facilitates the identification of the least and most influential model parameters. In this study, we employ two conventional norms, specifically L1- and L2-norms with respect to a logarithmic measure, to evaluate the significance of model parameters in terms of their impact on the local variability of the model outputs, namely storage, loss, and complex moduli. For instance, these norms can be defined for the LS indices of the storage modulus as follows:
\begin{equation}  \label{eq-L1}
\begin{aligned}
\lVert \bar{S}_{E^{\prime}, q_i} \rVert_{L_1} = \int \left| \bar{S}_{E^{\prime}, q_i} \right| \,d\log(a_T\omega), \lVert \bar{S}_{E^{\prime}, q_i} \rVert_{L_2} = \sqrt{\int \left| \bar{S}_{E^{\prime}, q_i} \right|^2 \,d\log(a_T\omega)}.
\end{aligned}
\end{equation}
A similar definition of the norms can be applied to the LS indices of the loss or complex modulus. These integrals are evaluated numerically using the trapezoidal rule available in MATLAB, as the LS indices are uniformly distributed in the logarithm space. Tables \ref{tab-L1-storage-FMG}, \ref{tab-L1-loss-FMG}, and \ref{tab-L1-complex-FMG} present the L1-norm of the LS indices for the storage, loss, complex moduli associated with each FMG-FMG model parameter, respectively. Additionally, Figure \ref{fig-barplots-L1-FMG} displays the mean value of these L1-norms for each nanocomposite system. This procedure of data visualization helps us to gain a clearer understanding of the relative influence of the model parameters.

As illustrated in Figure \ref{fig-barplots-SEp-L1-FMG}, \(\tau_{c_2}\) and \(\tau_{c_1}\) are identified as the least influential, while \(\alpha_1\) and \(E_{c_1}\) are recognized as the most influential FMG-FMG model parameters affecting the local variability of the storage modulus. However, considering the impact of parameters on the loss modulus, although \(\alpha_1\) and \(E_{c_1}\) emerge as the most influential contributors to the local variability of the loss modulus again (as shown Figure \ref{fig-barplots-SEpp-L1-FMG}), it is challenging to pinpoint two least significant model parameters. Apart from \(\tau_{c_2}\), which is also the least important model parameter for the loss modulus, the remaining model parameters, \(\tau_{c_1}\), \(E_{c_2}\), and \(\alpha_2\), demonstrate an equivalent level of importance regarding the local variability of the loss modulus.

As stated in Section \ref{sec.LSA}, we can approach our problem such that the current set of model parameters corresponds to a single scalar output. As a result, we have also calculated the magnitude of the LS indices of the complex modulus. Figure \ref{fig-barplots-SEstar-L1-FMG} vividly demonstrates that \(\alpha_1\) and \(E_{c_1}\) are the most influential model parameters, while \(\tau_{c_2}\) is identified as the least influential model parameter. The remaining parameters has the same level of influence on the output, a trend that is consistent with the observations made in the LS indices of the loss modulus. Notably, the L2-norms of these LS indices are also tabulated and visualized in a similar manner, and factor prioritization based on the L2-norms confirms the conclusions drawn from the L1-norm analysis. This consistency across different norms reinforces the reliability of our sensitivity analysis results.

\begin{table}[H]
\caption{L1-norm of the normalized LS indices of the storage modulus with respect to all FMG-FMG model parameters.} \label{tab-L1-storage-FMG}
    \begin{adjustbox}{max width=\textwidth,center=\textwidth}
    \begin{tabular}{cccccccc}
        \toprule
        \textbf{HSWF wt.\%}	& \textbf{xGnP wt.\%}
        & \textbf{\(\lVert \bar{S}_{E^{\prime}, E_{c_1}} \rVert_{L_1}\)}	
        & \textbf{\(\lVert \bar{S}_{E^{\prime}, \tau_{c_1}} \rVert_{L_1}\)}
        & \textbf{\(\lVert \bar{S}_{E^{\prime}, \alpha_{1}} \rVert_{L_1}\)}
        & \textbf{\(\lVert \bar{S}_{E^{\prime}, E_{c_2}} \rVert_{L_1}\)}
        & \textbf{\(\lVert \bar{S}_{E^{\prime}, \tau_{c_2}} \rVert_{L_1}\)}
        & \textbf{\(\lVert \bar{S}_{E^{\prime}, \alpha_{2}} \rVert_{L_1}\)}\\
        \midrule
        \multirow{4}{*}{20} & 0.0 & 13.5 & 4.6 & 32.2 & 9.6 & 0.6 & 8.2\\
			  	      & 0.5 & 14.2 & 4.7 & 34.0 & 8.9 & 0.4 & 4.8\\
			            & 1.0 & 13.8 & 4.0 & 27.9 & 9.3 & 0.6 & 8.0\\
                            & 1.5 & 13.3 & 4.4 & 31.0 & 9.7 & 0.6 & 7.8\\
        \midrule
        \multirow{4}{*}{30} & 0.0 & 12.7 & 2.7 & 17.2 & 10.3 & 0.4 & 4.4\\
			  	      & 0.5 & 12.8 & 2.8 & 18.7 & 10.2 & 0.5 & 5.3\\
			            & 1.0 & 12.6 & 2.8 & 18.2 & 10.4 & 0.5 & 5.4\\
                            & 1.5 & 12.5 & 2.7 & 17.6 & 10.5 & 0.6 & 6.5\\
        \midrule
        \multirow{4}{*}{40} & 0.0 & 11.8 & 2.0 & 13.5 & 11.2 & 0.3 & 2.8\\
			  	      & 0.5 & 11.8 & 1.9 & 12.9 & 11.2 & 0.3 & 2.9\\
			            & 1.0 & 12.1 & 2.0 & 14.2 & 11.0 & 0.3 & 3.6\\
                            & 1.5 & 14.3 & 2.1 & 17.5 & 8.7  & 0.1 & 0.9\\
        \bottomrule
    \end{tabular}
    \end{adjustbox}
\end{table}

\begin{table}[H]
\caption{L1-norm of the normalized LS indices of the loss modulus with respect to all FMG-FMG model parameters.} \label{tab-L1-loss-FMG}
    \begin{adjustbox}{max width=\textwidth,center=\textwidth}
    \begin{tabular}{cccccccc}
        \toprule
        \textbf{HSWF wt.\%}	& \textbf{xGnP wt.\%}
        & \textbf{\(\lVert \bar{S}_{E^{\prime\prime}, E_{c_1}} \rVert_{L_1}\)}	
        & \textbf{\(\lVert \bar{S}_{E^{\prime\prime}, \tau_{c_1}} \rVert_{L_1}\)}
        & \textbf{\(\lVert \bar{S}_{E^{\prime\prime}, \alpha_{1}} \rVert_{L_1}\)}
        & \textbf{\(\lVert \bar{S}_{E^{\prime\prime}, E_{c_2}} \rVert_{L_1}\)}
        & \textbf{\(\lVert \bar{S}_{E^{\prime\prime}, \tau_{c_2}} \rVert_{L_1}\)}
        & \textbf{\(\lVert \bar{S}_{E^{\prime\prime}, \alpha_{2}} \rVert_{L_1}\)}\\
        \midrule
        \multirow{4}{*}{20} & 0.0 & 18.7 & 6.0 & 45.0 & 4.3 & 0.22 & 1.1\\
			  	      & 0.5 & 20.0 & 6.5 & 49.9 & 3.0 & 0.08 & 1.8\\
			            & 1.0 & 19.0 & 5.4 & 37.9 & 4.1 & 0.22 & 1.1\\
			            & 1.0 & 19.0 & 5.4 & 37.9 & 4.1 & 0.22 & 1.1\\
                            & 1.5 & 18.7 & 6.0 & 44.2 & 4.3 & 0.20 & 1.3\\
        \midrule
        \multirow{4}{*}{30} & 0.0 & 19.3 & 4.2 & 28.0 & 3.7 & 0.08 & 2.6\\
			  	      & 0.5 & 19.1 & 4.1 & 28.6 & 3.9 & 0.11 & 2.4\\
			            & 1.0 & 19.0 & 4.1 & 28.2 & 4.0 & 0.11 & 2.4\\
                            & 1.5 & 18.5 & 4.0 & 26.5 & 4.5 & 0.16 & 2.3\\
        \midrule
        \multirow{4}{*}{40} & 0.0 & 19.9 & 3.1 & 23.8 & 3.2 & 0.03 & 2.7\\
			  	      & 0.5 & 19.7 & 3.0 & 22.3 & 3.3 & 0.04 & 2.8\\
			            & 1.0 & 19.5 & 2.8 & 21.9 & 3.5 & 0.05 & 2.8\\
                            & 1.5 & 22.1 & 2.6 & 23.4 & 0.9 & 0.00 & 0.9\\
        \bottomrule
    \end{tabular}
    \end{adjustbox}
\end{table}

\begin{table}[H]
\caption{L1-norm of the magnitude of the normalized LS indices of the complex modulus for the FMG-FMG model with respect to all model parameters.} \label{tab-L1-complex-FMG}
    \begin{adjustbox}{max width=\textwidth,center=\textwidth}
    \begin{tabular}{cccccccc}
        \toprule
        \textbf{HSWF wt.\%} & \textbf{xGnP wt.\%}
        & \textbf{\(\lVert \bar{S}_{E^{*}, E_{c_1}} \rVert_{L_1}\)}
        & \textbf{\(\lVert \bar{S}_{E^{*}, \tau_{c_1}} \rVert_{L_1}\)}
        & \textbf{\(\lVert \bar{S}_{E^{*}, \alpha_{1}} \rVert_{L_1}\)}
        & \textbf{\(\lVert \bar{S}_{E^{*}, E_{c_2}} \rVert_{L_1}\)}
        & \textbf{\(\lVert \bar{S}_{E^{*}, \tau_{c_2}} \rVert_{L_1}\)}
        & \textbf{\(\lVert \bar{S}_{E^{*}, \alpha_{2}} \rVert_{L_1}\)} \\
        \midrule
        \multirow{4}{*}{20}   & 0.0 & 14.2 & 4.8 & 32.9 & 9.2 & 0.6 & 8.0 \\
                              & 0.5 & 14.9 & 4.9 & 34.6 & 8.5 & 0.4 & 4.7 \\
                              & 1.0 & 14.2 & 4.8 & 32.9 & 9.2 & 0.6 & 8.0 \\
                              & 1.5 & 14.2 & 4.8 & 32.9 & 9.2 & 0.6 & 8.0 \\
        \midrule
        \multirow{4}{*}{30}   & 0.0 & 13.1 & 2.8 & 17.8 & 10.2 & 0.4 & 4.4 \\
                              & 0.5 & 13.2 & 2.9 & 19.3 & 10.0 & 0.5 & 5.3 \\
                              & 1.0 & 13.0 & 2.8 & 18.8 & 10.2 & 0.5 & 5.4 \\
                              & 1.5 & 12.9 & 2.8 & 18.2 & 10.4 & 0.6 & 6.4 \\
        \midrule
        \multirow{4}{*}{40}   & 0.0 & 12.0 & 2.0 & 13.9 & 11.1 & 0.2 & 2.8 \\
                              & 0.5 & 12.1 & 2.0 & 13.3 & 11.1 & 0.3 & 2.9 \\
                              & 1.0 & 12.3 & 2.0 & 14.6 & 10.9 & 0.3 & 3.6 \\
                              & 1.5 & 14.5 & 2.1 & 17.8 & 8.6  & 0.1 & 0.9 \\
        \bottomrule
    \end{tabular}
    \end{adjustbox}
\end{table}

\subsubsection{Factor Prioritization and Factor Fixing based on the LSA for the FMM-FMG Model Parameters}
\begin{figure}[H]
    \centering
    \begin{subfigure}{0.49\textwidth}
        \includegraphics[width=\linewidth, trim=3cm 0.5cm 5.5cm 1.5cm, clip]{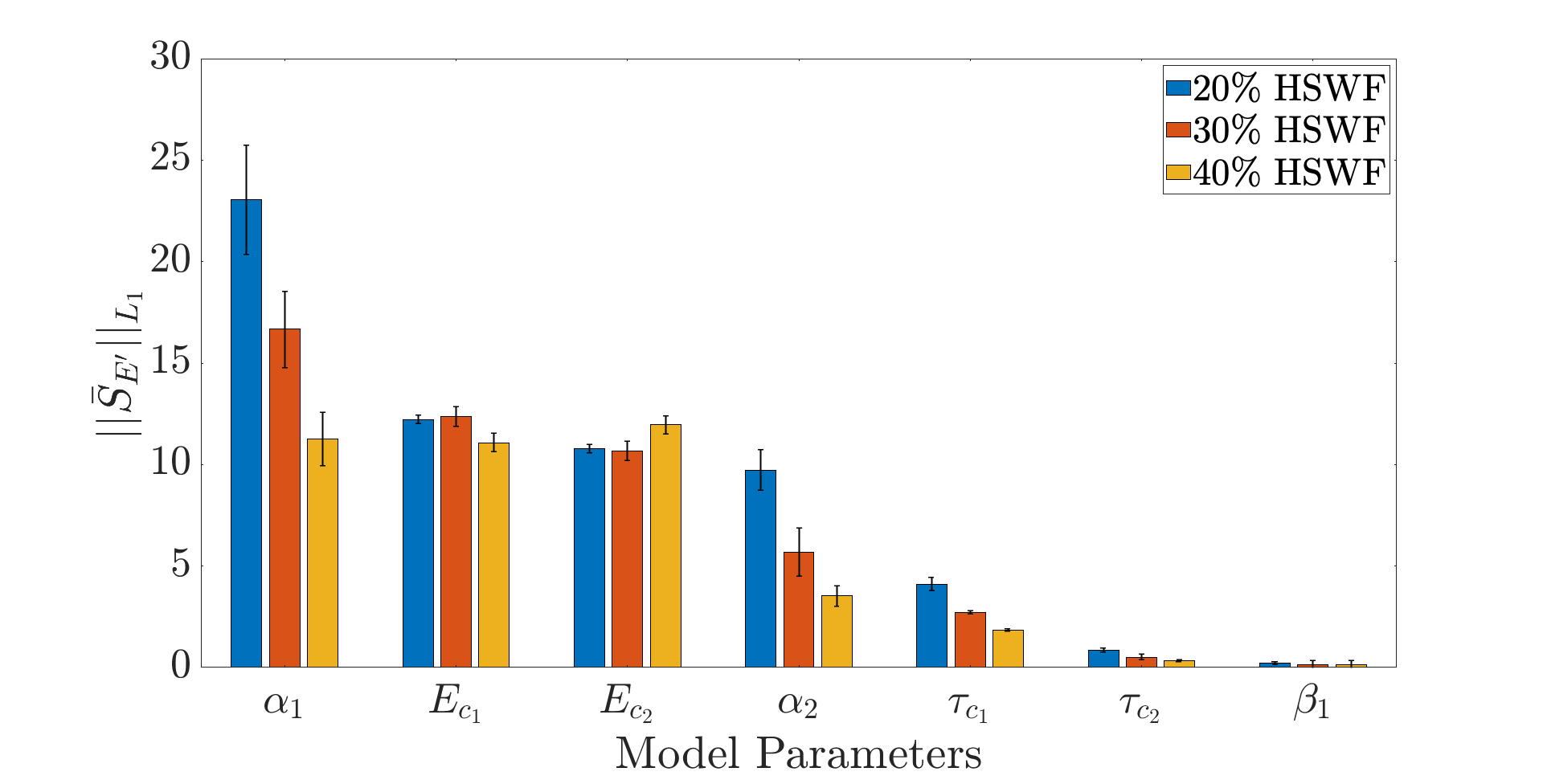}
        \caption{\label{fig-barplots-SEp-L1-FMM}}
    \end{subfigure}
    \begin{subfigure}{0.49\textwidth}
        \includegraphics[width=\linewidth, trim=3cm 0.5cm 5.5cm 1.5cm, clip]{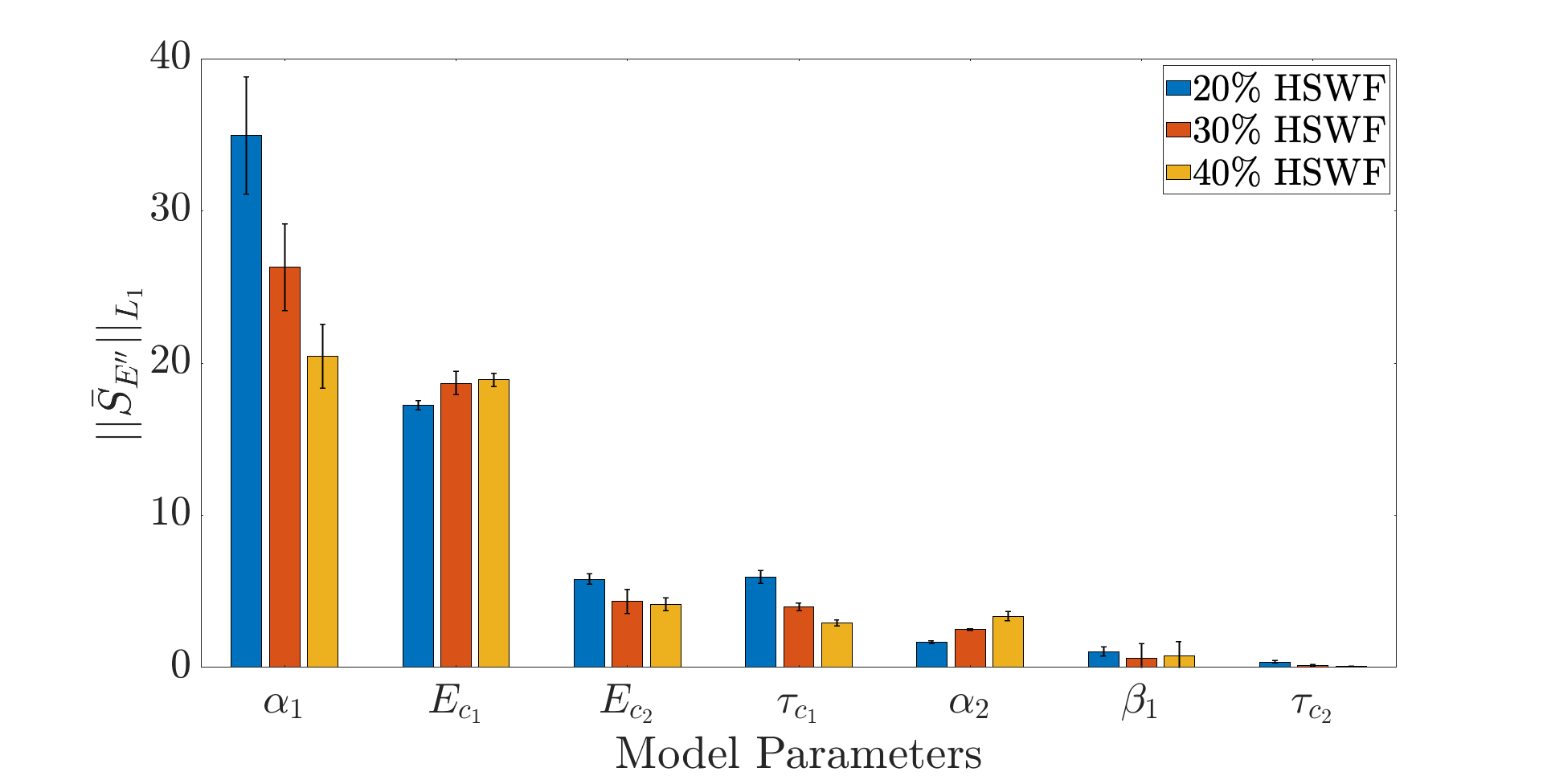}
        \caption{\label{fig-barplots-SEpp-L1-FMM}}
    \end{subfigure}
    \begin{subfigure}{0.49\textwidth}
        \includegraphics[width=\linewidth, trim=3cm 0.5cm 5.5cm 1.5cm, clip]{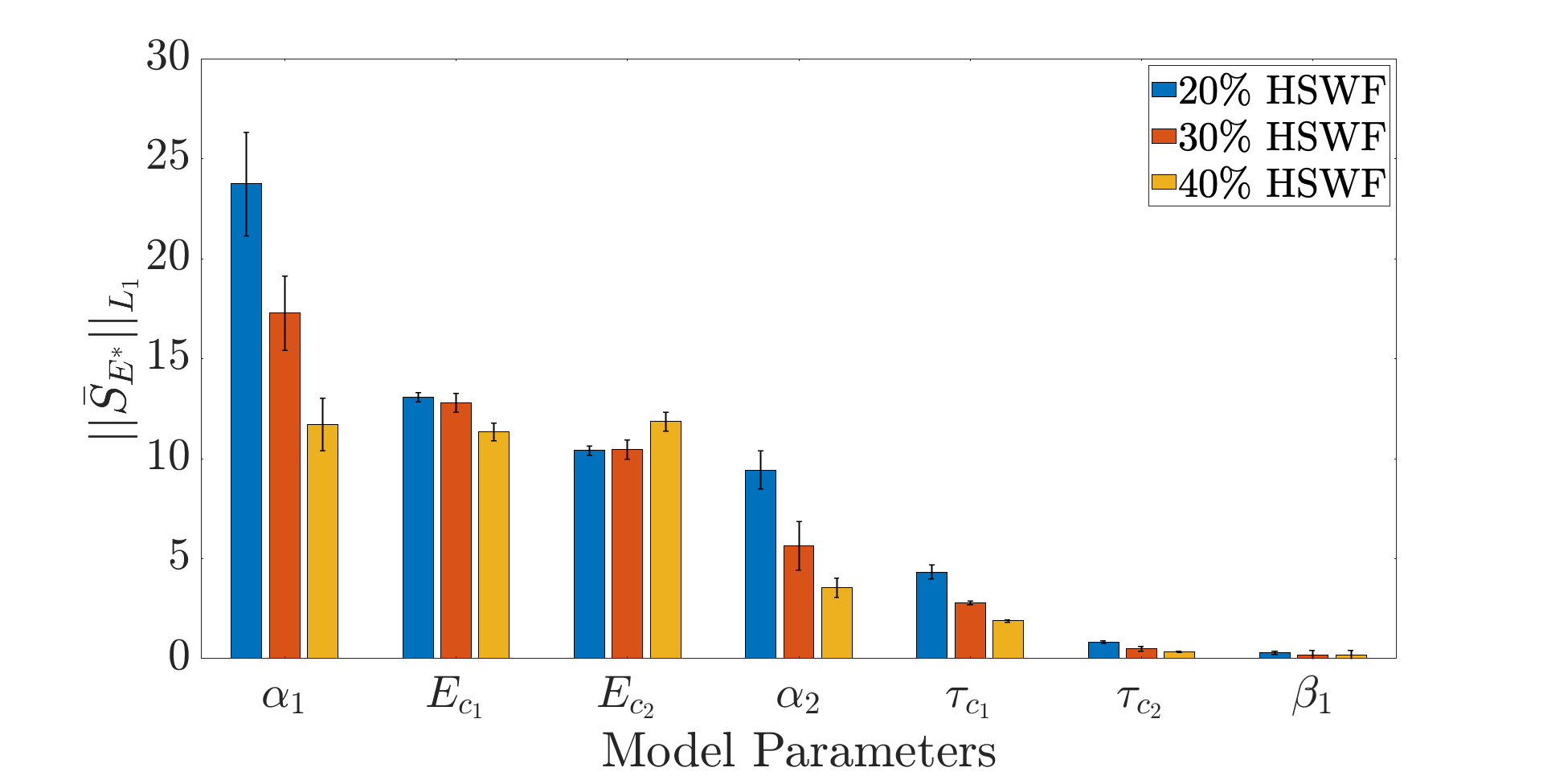}
        \caption{\label{fig-barplots-SEstar-L1-FMM}}
    \end{subfigure}
    \caption{Averaged L1-norm of the normalized local sensitivity indices associated with (a) storage, (b) loss, and (c) complex moduli with respect to each FMM-FMG model parameter.\label{fig-barplots-L1-FMM}}
\end{figure}

Employing a similar procedure, we aim to identify the least and most influential FMM-FMG model parameters that contribute to the local variability of the storage, loss, and complex moduli. As outlined in Section \ref{sec.MM}, the FMM-FMG model includes an extra parameter, \(\beta_1\), representing the power exponent of the second spring-pot element in the FMM branch, in addition to the six model parameters present in the FMG-FMG model. Tables \ref{tab-L1-storage-FMM}, \ref{tab-L1-loss-FMM}, and \ref{tab-L1-complex-FMM} list the L1-norm of the LS indices for the storage, loss, and complex moduli associated with all seven FMM-FMG model parameters, respectively. Additionally, Figure \ref{fig-barplots-L1-FMM} depicts the mean value of these L1-norms for all nanocomposite systems, providing a comprehensive visual representation of the parameters influence.

Considering the storage modulus (Figure \ref{fig-barplots-SEp-L1-FMM}), the most influential model parameter is identified as \(\alpha_1\), while \(\beta_1\) and \(\tau_{c_2}\) are deemed the least influential. \(E_{c_1}\) and \(E_{c_2}\) exhibit the equal level of influence on the storage modulus, and neither \(\alpha_2\) nor \(\tau_{c_1}\) has an L1-norm as low as \(\beta_1\) or \(\tau_{c_2}\) to be considered among the least influential model parameters. Moreover, from the perspective of the loss modulus (Figure \ref{fig-barplots-SEpp-L1-FMM}), \(\alpha_1\) and \(E_{c_1}\) can be recognized as the most influential, while \(\tau_{c_2}\) and \(\beta_1\) are the least influential model parameters. The remaining model parameters have equivalent effects on this modulus.

To complete the parameter prioritization of the FMM-FMG model, we shift our focus to the L1-norm of the LS indices for the complex modulus, since it reflects the combined effect of the model parameters on both the storage and loss moduli. As observed from Figure \ref{fig-barplots-SEstar-L1-FMM}, \(\alpha_1\) and \(E_{c_1}\) are the most influential and \(\tau_{c_2}\) and \(\beta_1\) are the least influential model parameters. The remaining model parameters (\(E_{c_2}\), \(\alpha_2\), and \(\tau_{c_2}\)) share an equivalent impact on the variation of the complex modulus. Remarkably, factor prioritization based on the L2-norms of these LS indices yields identical results, reinforcing the consistency and reliability of our findings.

To recap our efforts in identifying the least and most influential model parameters for both the FMG-FMG and FMM-FMG models, they are presented once again in Table \ref{tab-prio-LSA}. Given that the hard phase (FMG2) is almost perfectly elastic (see Table \ref{tab-opt}) with minimal dissipation, \(\tau_{c_2}\) is indeed the least influential parameter of both the FMG-FMG and FMM-FMG models in the local sense. Furthermore, since the second spring-pot element of the FMM branch also behaves nearly elastically, \(\beta_1\) naturally becomes the other least influential parameter of the FMM-FMG model from the local sensitivity perspective.

Thus far, our analysis has only considered the local effect of the model parameters on the variability of the model outputs. However, a global sensitivity analysis would further complete our tasks of factor prioritization and factor fixing through the calculation of the Sobol' indices, providing a more comprehensive understanding of parameters impacts.

\begin{table}[H]
\caption{L1-norm of the normalized LS indices of the storage modulus for the FMM-FMG model with respect to all model parameters.} \label{tab-L1-storage-FMM}
    \begin{adjustbox}{max width=\textwidth,center=\textwidth}
    \begin{tabular}{ccccccccc}
        \toprule
        \textbf{HSWF wt.\%}	& \textbf{xGnP wt.\%}
        & \textbf{\(\lVert \bar{S}_{E^{\prime}, E_{c_1}} \rVert_{L_1}\)}	
        & \textbf{\(\lVert \bar{S}_{E^{\prime}, \tau_{c_1}} \rVert_{L_1}\)}
        & \textbf{\(\lVert \bar{S}_{E^{\prime}, \alpha_{1}} \rVert_{L_1}\)}
        & \textbf{\(\lVert \bar{S}_{E^{\prime}, \beta_{1}} \rVert_{L_1}\)}
        & \textbf{\(\lVert \bar{S}_{E^{\prime}, E_{c_2}} \rVert_{L_1}\)}
        & \textbf{\(\lVert \bar{S}_{E^{\prime}, \tau_{c_2}} \rVert_{L_1}\)}
        & \textbf{\(\lVert \bar{S}_{E^{\prime}, \alpha_{2}} \rVert_{L_1}\)}\\
        \midrule
        \multirow{4}{*}{20} & 0.0 & 12.2 & 4.2 & 24.0 & 0.2 & 10.9 & 0.9 & 10.2\\
			  	      & 0.5 & 12.5 & 4.4 & 25.1 & 0.2 & 10.5 & 0.7 & 8.2\\
			            & 1.0 & 12.1 & 3.6 & 19.1 & 0.3 & 11.0 & 1.0 & 10.3\\
                            & 1.5 & 12.2 & 4.2 & 24.0 & 0.2 & 10.9 & 0.9 & 10.2\\
        \midrule
        \multirow{4}{*}{30} & 0.0 & 12.7 & 2.7 & 17.2 & 0.00 & 10.3 & 0.4 & 4.4\\
			  	      & 0.5 & 12.6 & 2.8 & 17.6 & 0.08 & 10.4 & 0.5 & 5.5\\
			            & 1.0 & 12.5 & 2.8 & 18.0 & 0.01 & 10.5 & 0.5 & 5.5\\
                            & 1.5 & 11.6 & 2.6 & 13.9 & 0.43 & 11.4 & 0.7 & 7.3\\
        \midrule
        \multirow{4}{*}{40} & 0.0 & 10.6 & 1.9 & 11.6 & 0.0 & 12.5 & 0.4 & 3.8\\
			  	      & 0.5 & 11.6 & 1.9 & 12.3 & 0.0 & 11.4 & 0.3 & 3.0\\
			            & 1.0 & 11.2 & 1.9 & 11.8 & 0.1 & 11.8 & 0.4 & 4.0\\
                            & 1.5 & 10.9 & 1.8 & 9.4  & 0.4 & 12.1 & 0.3 & 3.2\\
        \bottomrule
    \end{tabular}
    \end{adjustbox}
\end{table}

\begin{table}[H]
\caption{L1-norm of the normalized LS indices of the loss modulus for the FMM-FMG model with respect to all model parameters.} \label{tab-L1-loss-FMM}
    \begin{adjustbox}{max width=\textwidth,center=\textwidth}
    \begin{tabular}{ccccccccc}
        \toprule
        \textbf{HSWF wt.\%}	& \textbf{xGnP wt.\%}
        & \textbf{\(\lVert \bar{S}_{E^{\prime\prime}, E_{c_1}} \rVert_{L_1}\)}	
        & \textbf{\(\lVert \bar{S}_{E^{\prime\prime}, \tau_{c_1}} \rVert_{L_1}\)}
        & \textbf{\(\lVert \bar{S}_{E^{\prime\prime}, \alpha_{1}} \rVert_{L_1}\)}
        & \textbf{\(\lVert \bar{S}_{E^{\prime\prime}, \beta_{1}} \rVert_{L_1}\)}
        & \textbf{\(\lVert \bar{S}_{E^{\prime\prime}, E_{c_2}} \rVert_{L_1}\)}
        & \textbf{\(\lVert \bar{S}_{E^{\prime\prime}, \tau_{c_2}} \rVert_{L_1}\)}
        & \textbf{\(\lVert \bar{S}_{E^{\prime\prime}, \alpha_{2}} \rVert_{L_1}\)}\\
        \midrule
        \multirow{4}{*}{20} & 0.0 & 17.1 & 6.0 & 35.8 & 1.0 & 5.9 & 0.4 & 1.6\\
			  	      & 0.5 & 17.7 & 6.4 & 38.6 & 0.8 & 5.3 & 0.3 & 1.7\\
			            & 1.0 & 17.0 & 5.4 & 29.5 & 1.5 & 6.1 & 0.4 & 1.8\\
                            & 1.5 & 17.1 & 6.0 & 35.8 & 1.0 & 5.9 & 0.4 & 1.6\\
        \midrule
        \multirow{4}{*}{30} & 0.0 & 19.3 & 4.2 & 28.0 & 0.0 & 3.7 & 0.1 & 2.6\\
			  	      & 0.5 & 18.9 & 4.0 & 27.3 & 0.5 & 4.1 & 0.1 & 2.5\\
			            & 1.0 & 18.9 & 4.0 & 28.0 & 0.1 & 4.1 & 0.1 & 2.5\\
                            & 1.5 & 17.5 & 3.6 & 22.0 & 2.0 & 5.5 & 0.2 & 2.5\\
        \midrule
        \multirow{4}{*}{40} & 0.0 & 18.5 & 3.1 & 22.2 & 0.2 & 4.5 & 0.1 & 3.6\\
			  	      & 0.5 & 19.5 & 3.1 & 21.8 & 0.2 & 3.5 & 0.0 & 3.0\\
			            & 1.0 & 18.7 & 2.8 & 20.2 & 0.6 & 4.3 & 0.1 & 3.4\\
                            & 1.5 & 18.8 & 2.7 & 17.5 & 2.2 & 4.2 & 0.1 & 3.5\\
        \bottomrule
    \end{tabular}
    \end{adjustbox}
\end{table}

\begin{table}[H]
\caption{L1-norm of the magnitude of the normalized LS indices of the complex modulus for the FMM-FMG model with respect to all model parameters.} \label{tab-L1-complex-FMM}
    \begin{adjustbox}{max width=\textwidth,center=\textwidth}
    \begin{tabular}{cccccccccc}
        \toprule
        \textbf{HSWF wt.\%} & \textbf{xGnP wt.\%} 
        & \textbf{\(\lVert \bar{S}_{E^{\prime}, E_{c_1}} \rVert_{L_1}\)}
        & \textbf{\(\lVert \bar{S}_{E^{\prime}, \tau_{c_1}} \rVert_{L_1}\)}
        & \textbf{\(\lVert \bar{S}_{E^{\prime}, \alpha_{1}} \rVert_{L_1}\)}
        & \textbf{\(\lVert \bar{S}_{E^{\prime}, \beta_{1}} \rVert_{L_1}\)}
        & \textbf{\(\lVert \bar{S}_{E^{\prime}, E_{c_2}} \rVert_{L_1}\)}
        & \textbf{\(\lVert \bar{S}_{E^{\prime}, \tau_{c_2}} \rVert_{L_1}\)}
        & \textbf{\(\lVert \bar{S}_{E^{\prime}, \alpha_{2}} \rVert_{L_1}\)} \\
        \midrule
        \multirow{4}{*}{20HS} & 0.0 & 13.1 & 4.5 & 24.8 & 0.3 & 10.4 & 0.9 & 10.0 \\
                            & 0.5 & 13.4 & 4.6 & 25.9 & 0.2 & 10.1 & 0.7 & 8.0 \\
                            & 1.0 & 12.8 & 3.8 & 20.0 & 0.4 & 10.6 & 0.9 & 10.1 \\
                            & 1.5 & 13.0 & 4.4 & 24.3 & 0.3 & 10.5 & 0.8 & 9.6 \\
        \midrule
        \multirow{4}{*}{30HS} & 0.0 & 13.1 & 2.8 & 17.8 & 0.00 & 10.2 & 0.4 & 4.4 \\
                            & 0.5 & 13.0 & 2.9 & 18.2 & 0.12 & 10.2 & 0.5 & 5.5 \\
                            & 1.0 & 13.0 & 2.8 & 18.6 & 0.02 & 10.3 & 0.5 & 5.4 \\
                            & 1.5 & 12.1 & 2.7 & 14.5 & 0.56 & 11.2 & 0.7 & 7.3 \\
        \midrule
        \multirow{4}{*}{40HS} & 0.0 & 10.8 & 1.9 & 12.0 & 0.0 & 12.4 & 0.4 & 3.9 \\
                            & 0.5 & 11.9 & 1.9 & 12.7 & 0.0 & 11.3 & 0.3 & 3.1 \\
                            & 1.0 & 11.5 & 1.9 & 12.3 & 0.1 & 11.7 & 0.4 & 4.1 \\
                            & 1.5 & 11.2 & 1.8 & 9.8 & 0.5 & 12.0 & 0.3 & 3.2 \\
        \bottomrule
    \end{tabular}
    \end{adjustbox}
\end{table}

\begin{table}[H]
\caption{The least and most influential model parameters for both FMG-FMG and FMM-FMG models based on the LSA.} \label{tab-prio-LSA}
    \centering
    \begin{adjustbox}{max width=\textwidth,center=\textwidth}
        \begin{tabular}{ccc}
            \toprule
            \textbf{Model} & \textbf{Least Influential Parameters} & \textbf{Most Influential Parameters} \\
            \midrule
            FMG-FMG & \makecell{\(\tau_{c_2}\) \\ -} & \makecell{\(\alpha_1\) \\ \(E_{c_1}\)} \\
            FMM-FMG & \makecell{\(\tau_{c_2}\) \\ \(\beta_1\)} & \makecell{\(\alpha_1\) \\ \(E_{c_1}\)} \\
            \bottomrule
        \end{tabular}
    \end{adjustbox}
\end{table}

\subsection{Variation of LS Indices}
\begin{figure}[H]
    \centering
    \begin{subfigure}{0.74\textwidth}
        \centering
        \includegraphics[width=\textwidth, trim=4cm 0.25cm 4.5cm 1.5cm, clip]{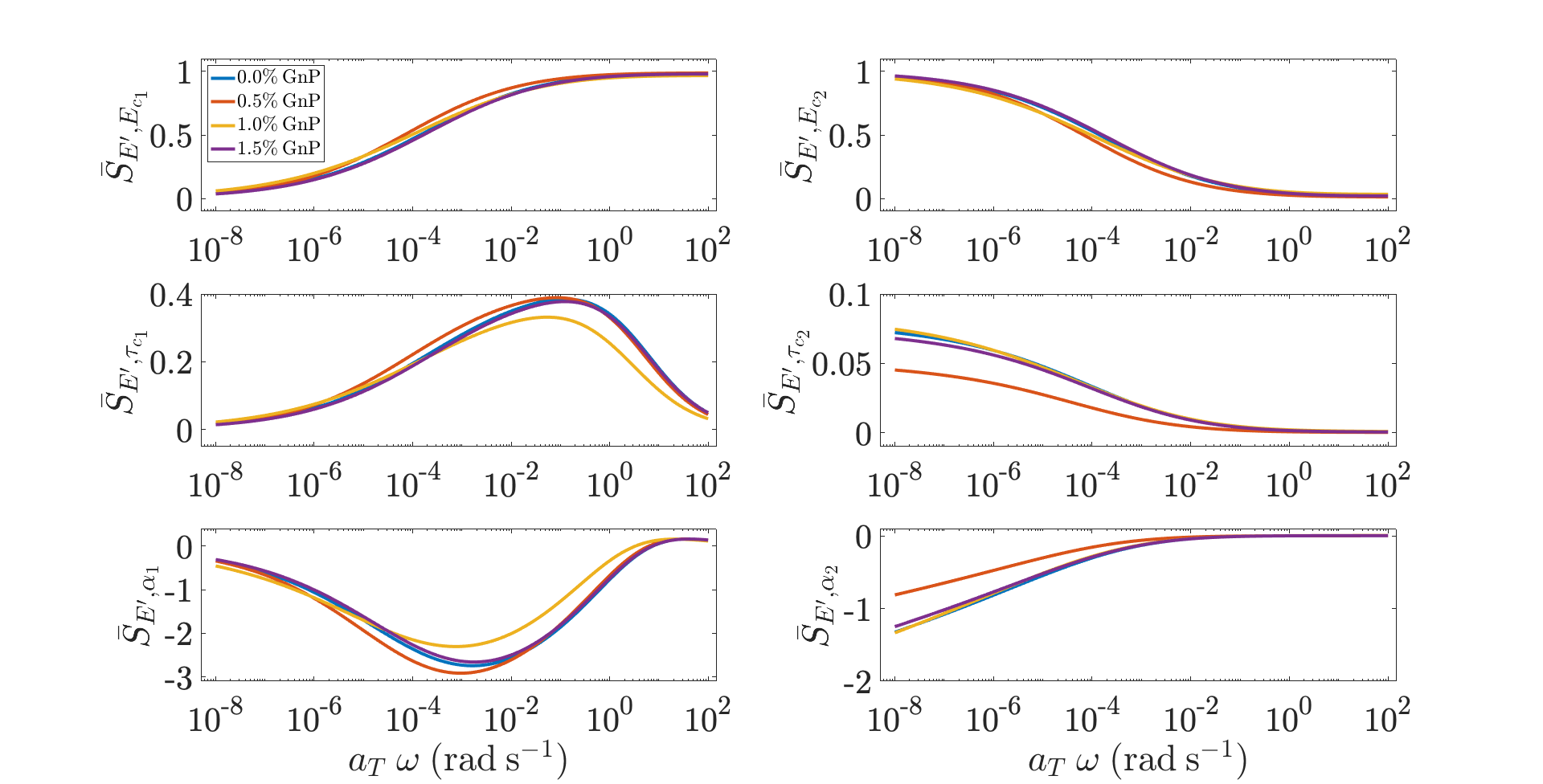}
        \caption{\label{fig-LSA-Ep-FMG-20HS}}
    \end{subfigure}
    \begin{subfigure}{0.74\textwidth}
        \centering
        \includegraphics[width=\textwidth, trim=4cm 0.25cm 4.5cm 1.5cm, clip]{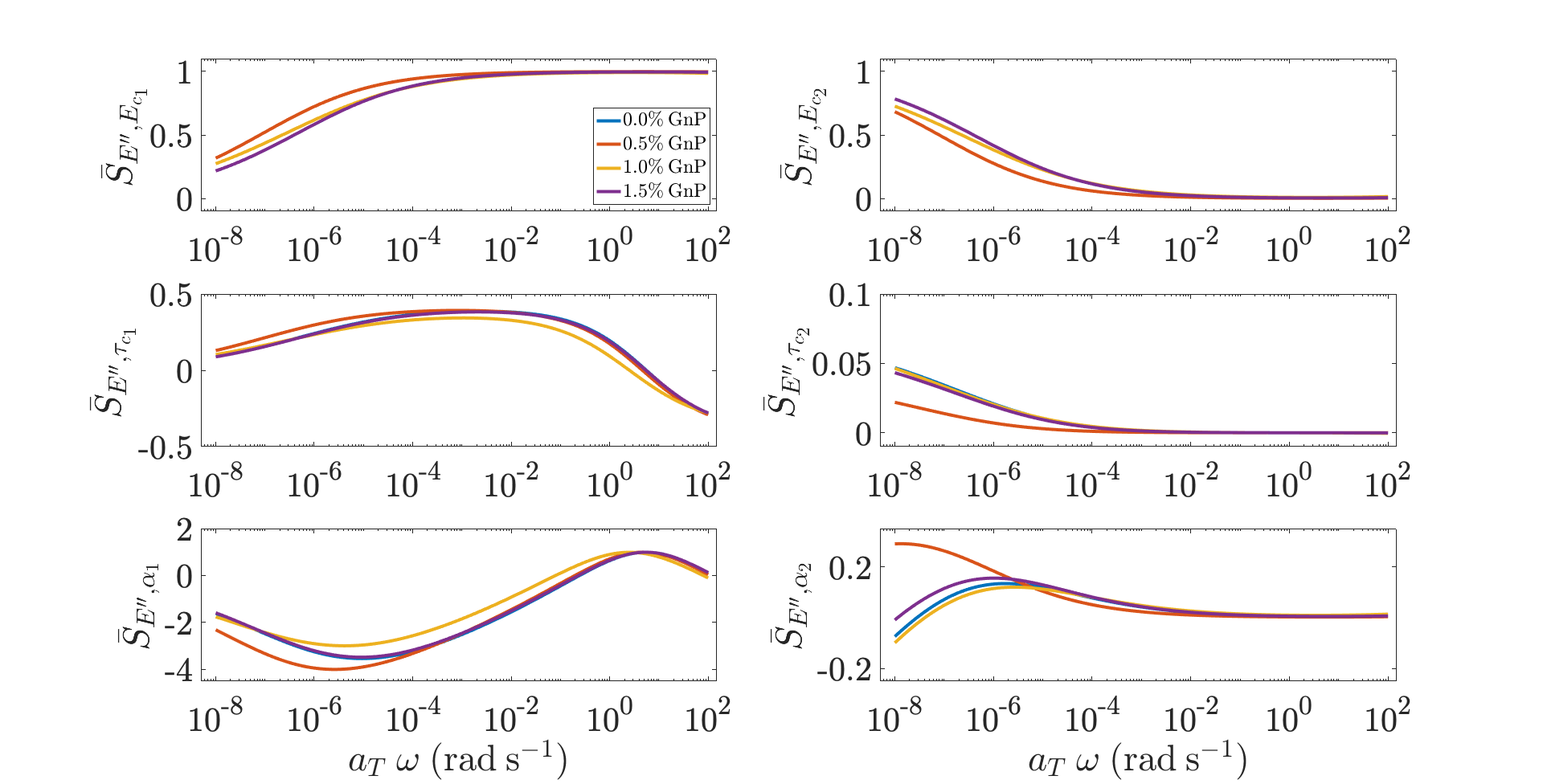}
        \caption{\label{fig-LSA-Epp-FMG-20HS}}
    \end{subfigure}
    \begin{subfigure}{0.74\textwidth}
        \centering
        \includegraphics[width=\textwidth, trim=4cm 0.25cm 4.5cm 1.5cm, clip]{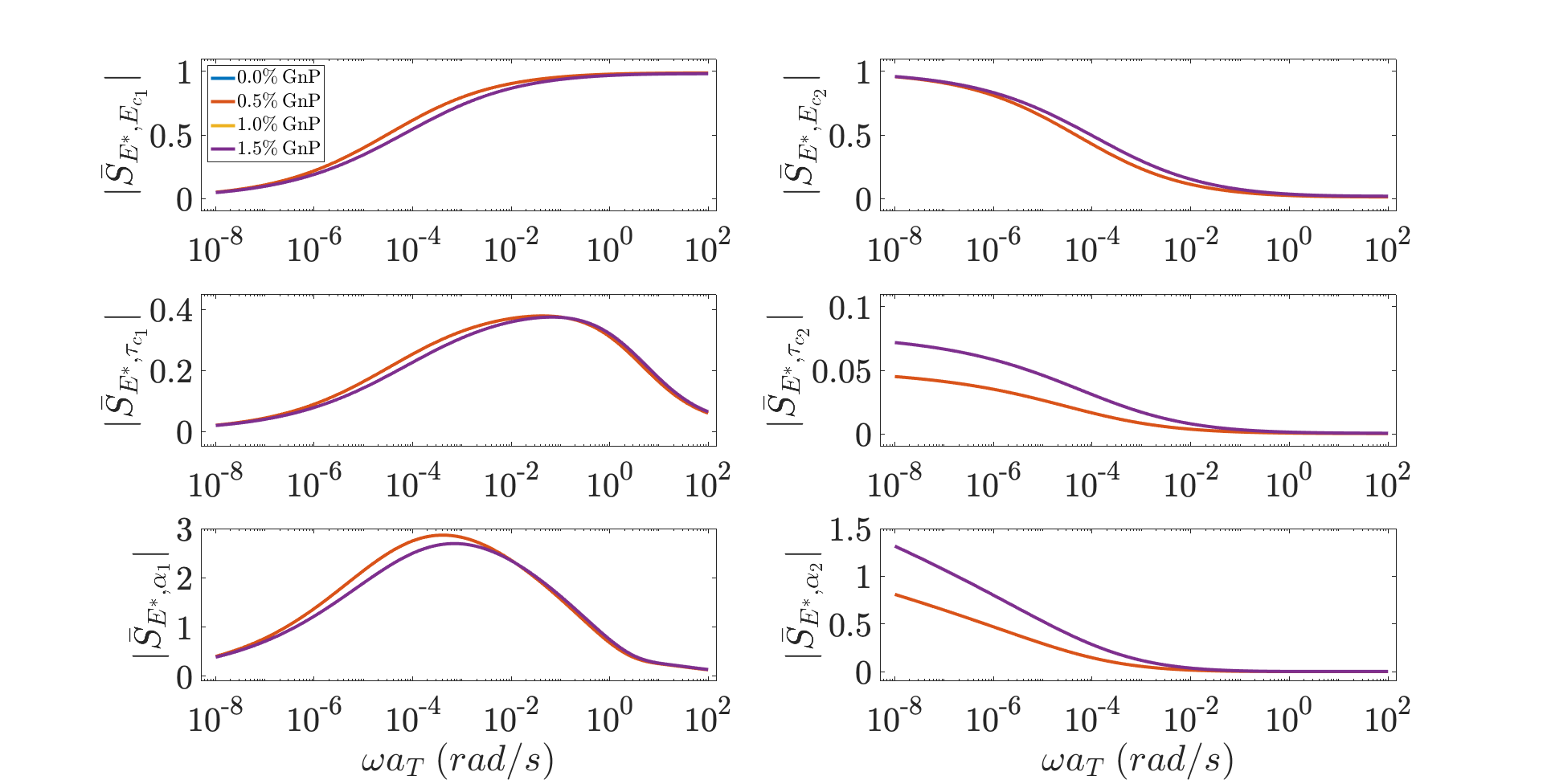}
        \caption{\label{fig-LSA-Estar-FMG-20HS}}
    \end{subfigure}
    \caption{Variation of the averaged normalized LS indices of the (a) storage, (b) loss, and (c) complex moduli of the FMG-FMG model corresponding to 20\% HSWF nanocomposite systems. The left and right columns display the normalized LS indices with respect to the model parameters in the first and second branch, respectively.} \label{fig-LSA-FMG-20HS}
\end{figure}

Focusing initially on the FMG-FMG model, Figures \ref{fig-LSA-Ep-FMG-20HS}, \ref{fig-LSA-Epp-FMG-20HS}, and \ref{fig-LSA-Estar-FMG-20HS} show the variation of the averaged normalized local sensitivity (LS) indices for the storage, loss, and complex modulus with respect to all model parameters for the 20\% HSWF nanocomposite systems. Analyzing these indices reveals intriguing behaviors and trends.

Beginning with the analysis of the normalized LS indices corresponding to the storage modulus (Figure \ref{fig-LSA-Ep-FMG-20HS}), both \(\bar{S}_{E^{\prime}, E_{c_1}}\) and \(\bar{S}_{E^{\prime}, E_{c_2}}\) exhibit complementary behaviors. \(\bar{S}_{E^{\prime}, E_{c_1}}\) displays a monotone increasing trend, starting from a nearly zero value around a shifted frequency of \(a_t\omega \approx 10^{-8} \text{ rad s}^{-1}\), and plateaus at value of one around \(a_t\omega \approx 10^0 \text{ rad s}^{-1}\), maintaining this value for the subsequent frequency decades. Conversely, \(\bar{S}_{E^{\prime}, E_{c_2}}\) begins at a value of one in low-frequency decades, 
monotonically decreases to zero in high-frequency decades, and reaches a similar plateau trend thereafter. Additionally, the behavior of \(\bar{S}_{E^{\prime}, \tau_{c_1}}\) and \(\bar{S}_{E^{\prime}, \alpha_{1}}\) appears correlated, with both indices exhibiting a bell-shaped distribution. Peak values are observed in the mid-range shifted frequencies (\(a_t\omega \approx 10^{-5}\) to \(\approx 10^{1}\) rad s\(^{-1}\)), while both indices approach zero at very low and very high shifted frequency decades. The remaining two local sensitivity indices, \(\bar{S}_{E^{\prime}, \tau_{c_2}}\) and \(\bar{S}_{E^{\prime}, \alpha_{2}}\), demonstrate a monotone decrease in their magnitude as the shifted frequency increases, eventually reaching a plateau at a value of zero.

Continuing with the normalized LS indices associate with the loss modulus (Figure \ref{fig-LSA-Epp-FMG-20HS}), a consistent trend emerges similar to the LS indices of the storage modulus, particularly with respect to the model parameters of the second branch (FMG2). All \(\bar{S}_{E^{\prime\prime},FMG2}\) indices exhibit a decreasing trend, ultimately reaching a near-zero plateau at high frequency decades. Moreover, akin to the behavior observed for \(\bar{S}_{E^{\prime}, E_{c_1}}\), \(\bar{S}_{E^{\prime\prime}, E_{c_1}}\) also demonstrates a plateau at high frequency ranges with a value of one. However, unlike the correlated behavior of \(\bar{S}_{E^{\prime}, \tau_{c_1}}\) and \(\bar{S}_{E^{\prime}, \alpha_{1}}\), both \(\bar{S}_{E^{\prime\prime}, \tau_{c_1}}\) and \(\bar{S}_{E^{\prime\prime}, \alpha_{1}}\) display more drastic variations across frequency decades, with relatively large values observed at both low and high frequency ranges.

Finally, the magnitude of the normalized LS indices of the complex modulus are shown in Figure \ref{fig-LSA-Estar-FMG-20HS}. As their definition suggests \eqref{eq-SEstar}, their behaviors and trends are a synthesis of the LS indices of both the storage and loss moduli. This results in a monotone decrease in the LS indices associated with the second branch parameters, \(\bar{S}_{E^*, FMG2}\), a monotone increase in the \(\bar{S}_{E^*, E_{c_1}}\), and significant variations of the \(\bar{S}_{E^*, \alpha_1}\) and \(\bar{S}_{E^*, \tau_{c,1}}\). The remaining normalized LS indices corresponding to the 30\% and 40\% HSWF nanocomposite systems are provided in the Figures S4 and S5.

\begin{figure}[H]
    \centering
    \begin{subfigure}{0.74\textwidth}
        \includegraphics[width=\textwidth, trim=3cm 0.1cm 5.5cm 2cm, clip]{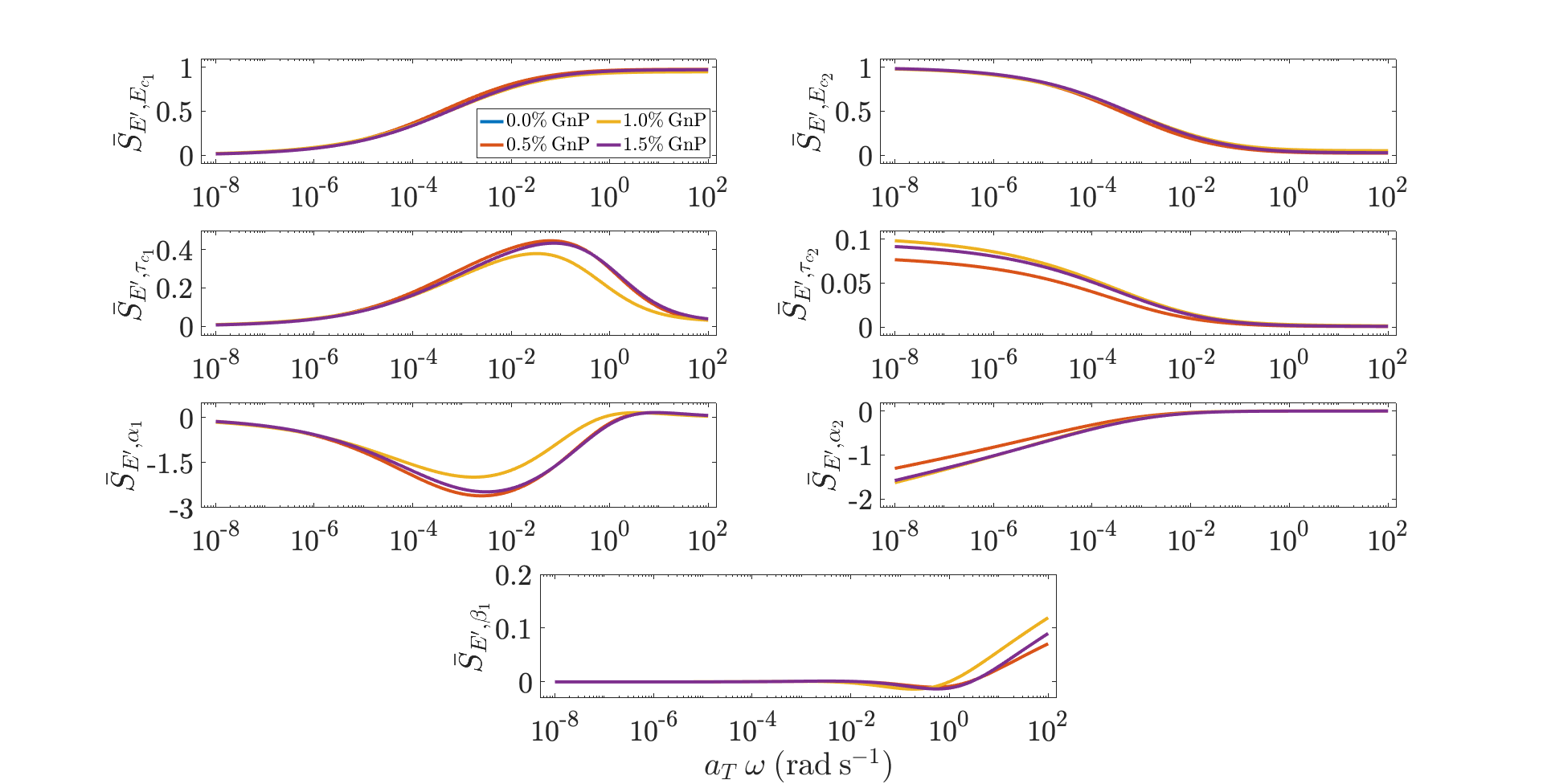}
        \caption{\label{fig-LSA-Ep-FMM-20HS}}
    \end{subfigure}
    \centering
    \begin{subfigure}{0.74\textwidth}
        \includegraphics[width=\textwidth, trim=3cm 0.1cm 5.5cm 2cm, clip]{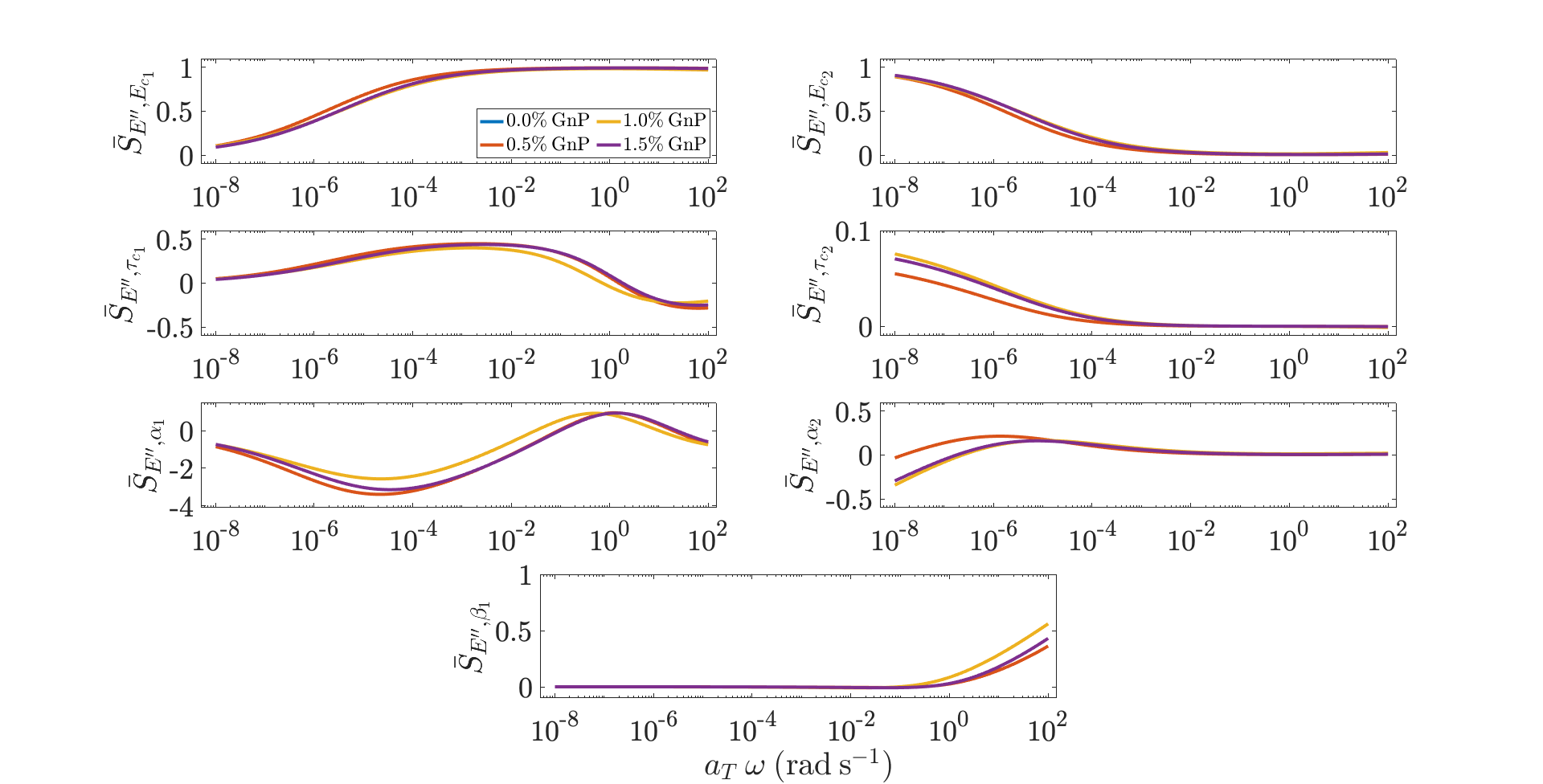}
        \caption{\label{fig-LSA-Epp-FMM-20HS}}
    \end{subfigure}
    \begin{subfigure}{0.74\textwidth}
        \includegraphics[width=\textwidth, trim=3cm 0.1cm 5.5cm 2cm, clip]{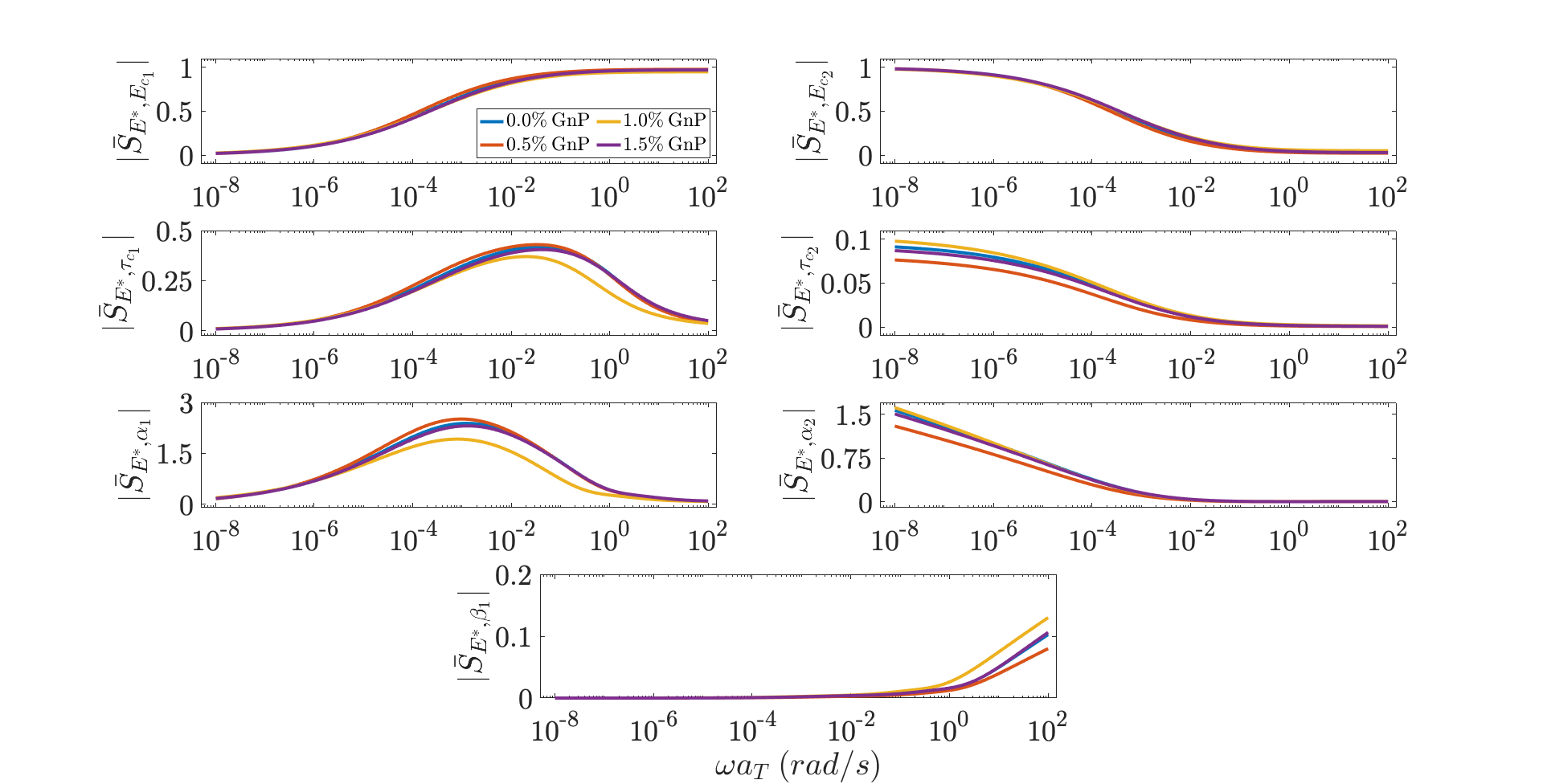}
        \caption{\label{fig-LSA-Estar-FMM-20HS}}
    \end{subfigure}
    \caption{Variation of the averaged normalized LS indices of the (a) storage, (b) loss moduli, and (c) complex moduli of the FMM-FMG model corresponding to the 20\% HSWF nanocomposite systems. The left and right columns display the normalized LS indices with respect to the model parameters in the first and second branch, respectively, with the last row dedicated to \(\beta_1\).} \label{fig-LSA-20HS-FMM}
\end{figure}

Switching to the FMM-FMG model, Figures \ref{fig-LSA-Ep-FMM-20HS}, \ref{fig-LSA-Epp-FMM-20HS}, and \ref{fig-LSA-Estar-FMM-20HS} illustrate the variation of the averaged normalized local sensitivity indices for the storage, loss, and complex moduli with respect to all model parameters for the 20\% HSWF nanocomposite systems. A comparison of these results with the LS indices of the FMG-FMG model uncovers numerous similarities in their trends and behaviors. In terms of the LS indices of the storage modulus (Figure \ref{fig-LSA-Ep-FMM-20HS}), the complementary behavior of \(\bar{S}_{E^{\prime}, E_{c_1}}\) and \(\bar{S}_{E^{\prime}, E_{c_2}}\) is also observed in this model. Additionally, the correlated behavior of \(\bar{S}_{E^{\prime}, \tau_{c_1}}\) and \(\bar{S}_{E^{\prime}, \alpha_{1}}\), with a Gaussian-like distribution and peak values at mid-range frequency magnitudes, is replicated in this model. Furthermore, the decreasing trend in the magnitude of the LS indices respective to the second branch model parameters is repeated in this model as well. Lastly, the first branch of the FMM-FMG model includes an additional model parameter, \(\beta_1\), which exhibits a near-zero LS index across multiple decades of frequency up to \(a_t\omega \approx 10^0 \) rad s\(^{-1}\), beyond which it begins to rise quite drastically.

Transitioning to the LS indices of the loss modulus for the FMM-FMG model (Figure \ref{fig-LSA-Epp-FMM-20HS}), we observe trends that mirror those described for the FMG-FMG model. Specifically, the increasing trend of \(\bar{S}_{E^{\prime\prime}, E_{c_1}}\), the drastic variations of \(\bar{S}_{E^{\prime\prime}, \tau_{c_1}}\) and \(\bar{S}_{E^{\prime\prime}, \alpha_{1}}\), and the decreasing trend in the magnitude of the LS indices corresponding to the parameters of the second branch, \(\bar{S}_{E^{\prime\prime},FMG2}\), are consistent with previous observations. Similar to the behavior of \(\bar{S}_{E^{\prime}, \beta_{1}}\), \(\bar{S}_{E^{\prime\prime}, \beta_{1}}\) also maintains a value close to zero up to \(a_t\omega \approx 10^0 \) rad s\(^{-1}\), but it increases thereafter at a more pronounced rate. Lastly, the variation of the magnitude of the LS indices of the complex modulus is presented in Figure \ref{fig-LSA-Estar-FMM-20HS}, which reflects the combined behavior of the LS indices of the storage and loss modulus, as it definition suggests in \eqref{eq-SEstar}. The remaining normalized LS indices for the 30\% and 40\% HSWF nanocomposite systems are detailed in Figures S6 and S7.

\subsection{Global Sensitivity Analysis}
Next, we present the procedure for identifying the influential and non-influential model parameters through the evaluation of Sobol' indices for all three moduli associated with each model parameter. As introduced in Section \ref{sec.GSA}, two common Sobol' indices that are used in GSA are the first- and total-order sensitivity indices. The first-order index measures the main effect contribution of each input factor to the variance of the output, while the total-order index additionally accounts for all higher-order effects due to interactions. However, in the current problem, although our models are non-linear, they are additive. Therefore, the total-order sensitivity index is equal to the first-order sensitivity index (refer to Figure S8), and the contribution of the interaction effects on the variability of the model output is minimal. Consequently, \(S_i=0\) is a necessary and sufficient condition for \(q_i\) to be considered as a non-influential factor. In other words, if \(S_i \approx 0\), then \(q_i\) can be fixed at any value within its range of variability without significantly influencing the variability of the output \cite{saltelli_global_2008}.

\subsubsection{Factor Prioritization and Factor Fixing based on the GSA for the FMG-FMG Model Parameters}
\begin{figure}[H]
    \centering
    \begin{subfigure}{\textwidth}
        \centering
        \includegraphics[width=0.74\linewidth, trim=4cm 0.25cm 4.5cm 1.5cm, clip]{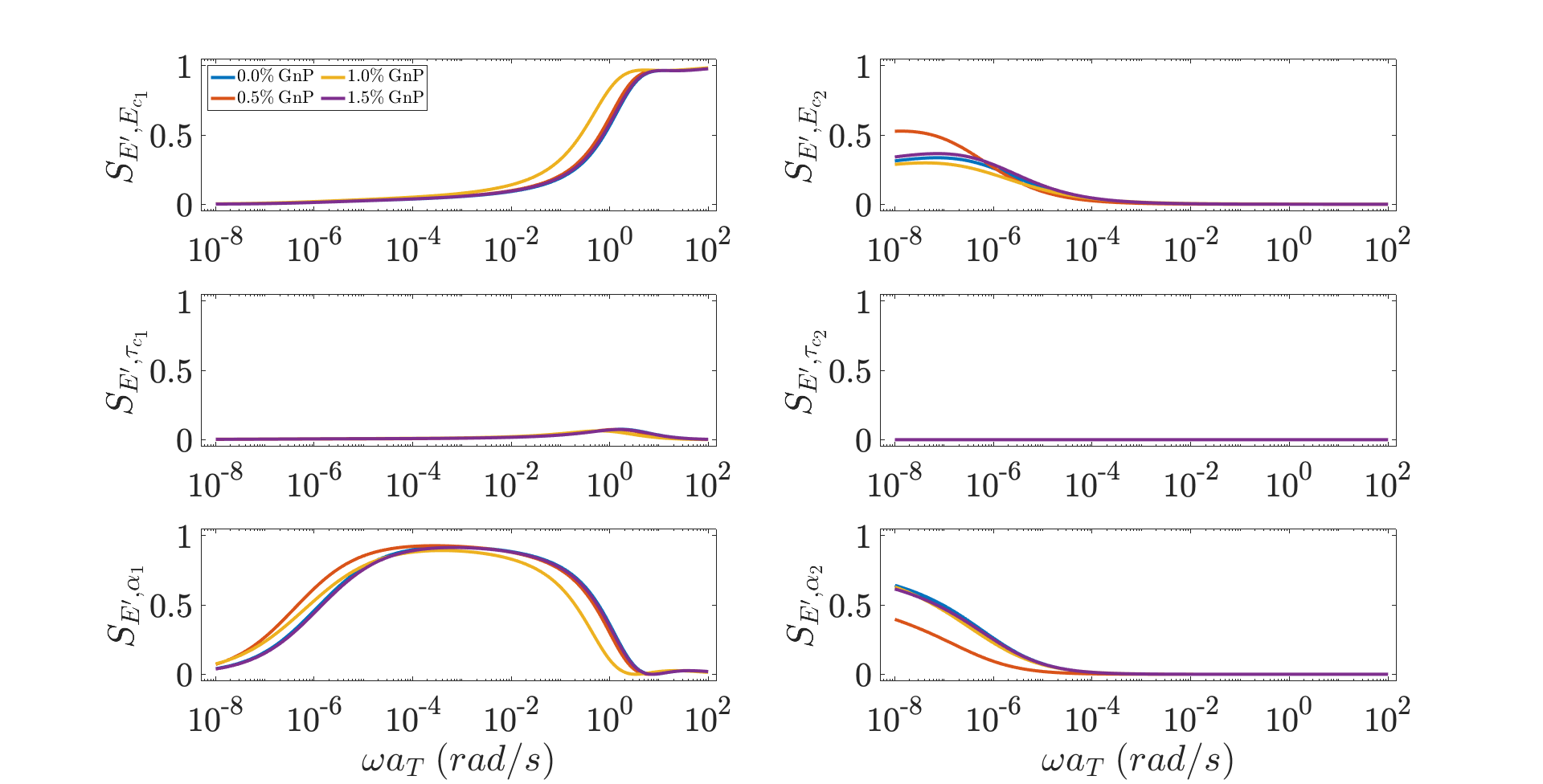}
        \caption{\label{fig-GSA-SEp-FMG}}
    \end{subfigure}
    \begin{subfigure}{\textwidth}
        \centering
        \includegraphics[width=0.74\linewidth, trim=4cm 0.25cm 4.5cm 1.5cm, clip]{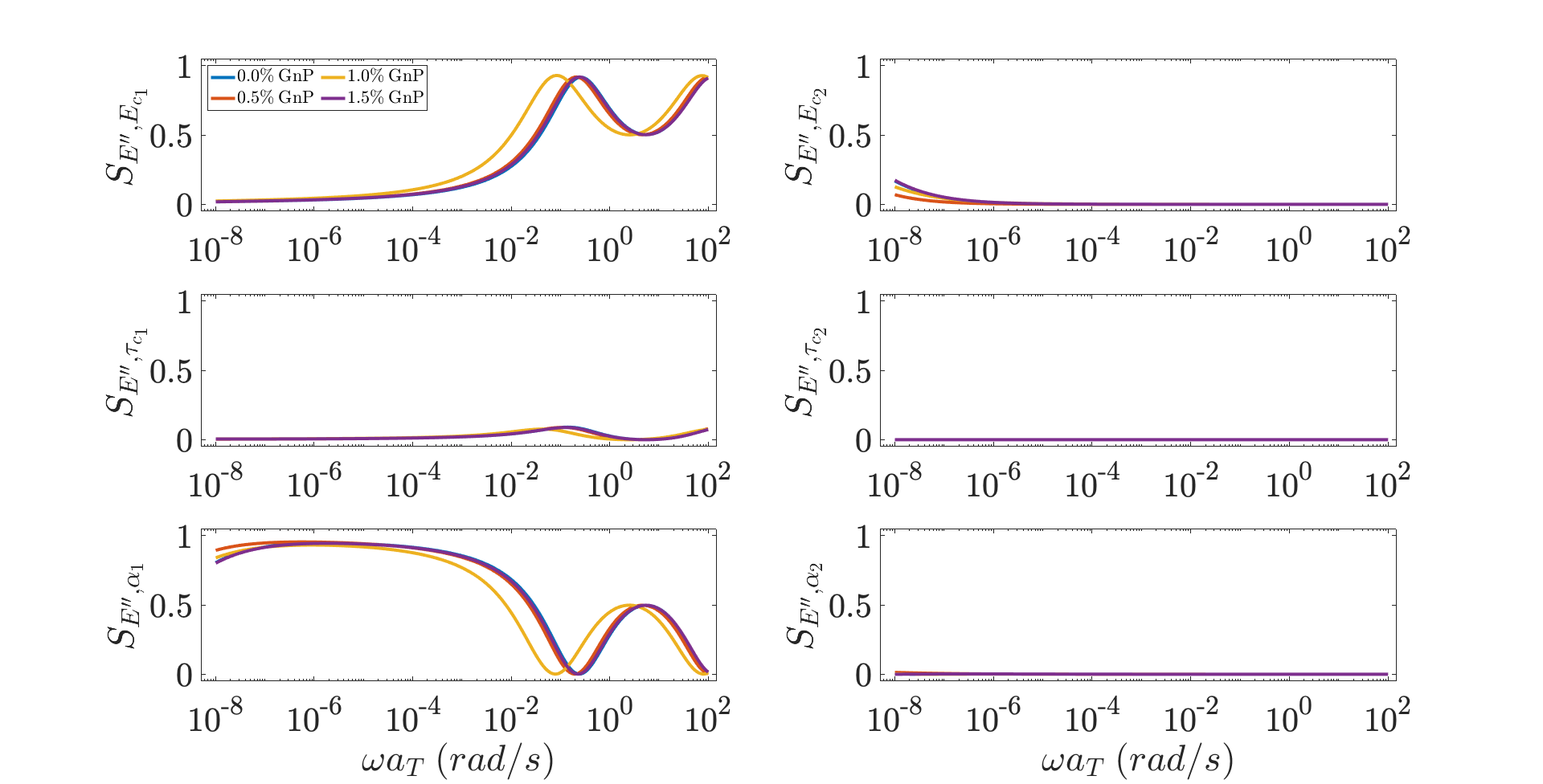}
        \caption{\label{fig-GSA-SEpp-FMG}}
    \end{subfigure}
    \begin{subfigure}{0.74\textwidth}
        \includegraphics[width=\linewidth, trim=4cm 0.25cm 4.5cm 1.5cm, clip]{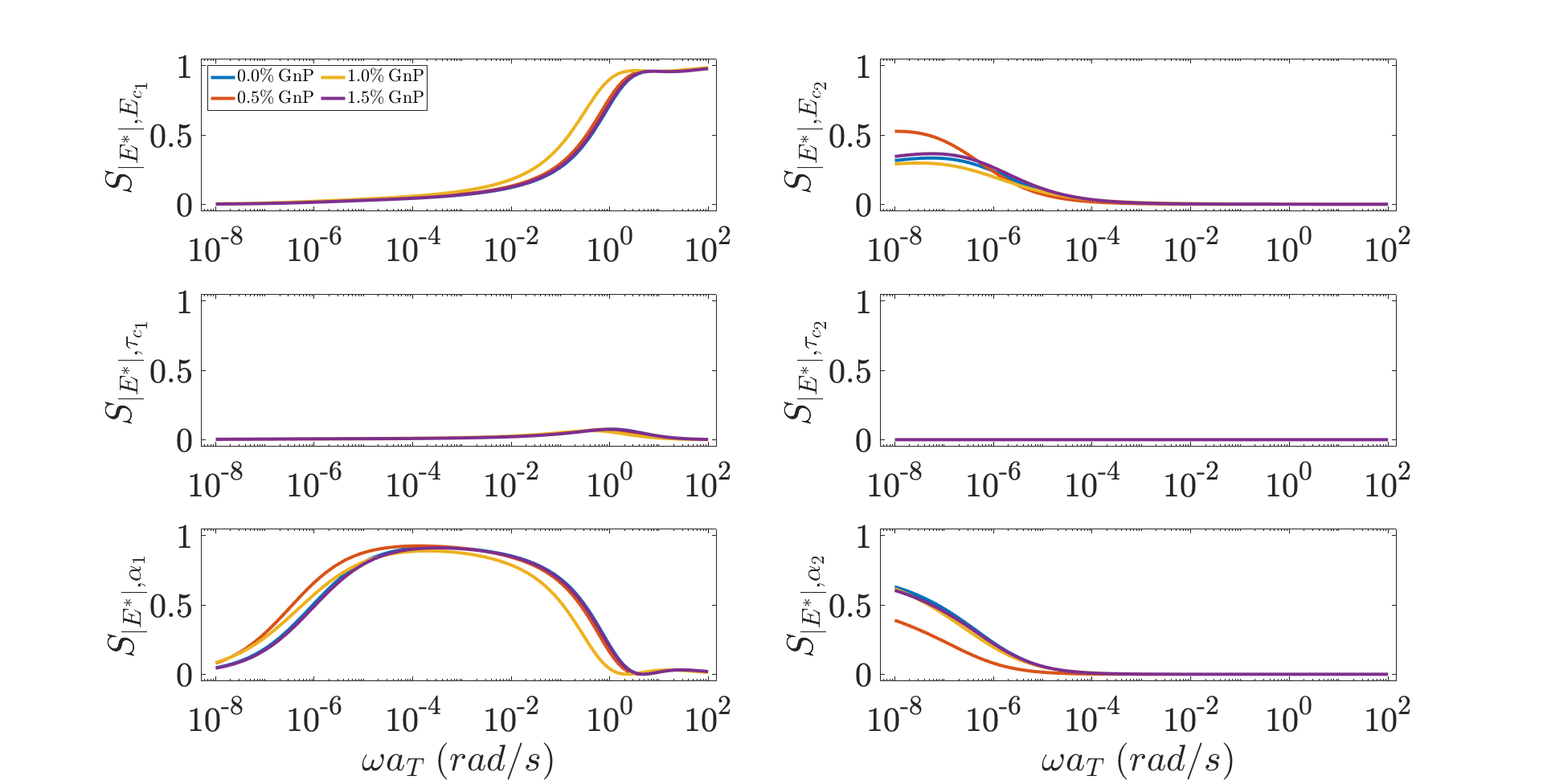}
        \caption{\label{fig-GSA-SEstar-FMG}}
    \end{subfigure}
    \caption{Variation of the first-order sensitivity indices of the (a) storage, (b) loss moduli, and (c) magnitude of the complex modulus of the FMG-FMG model corresponding to the 20\% HSWF nanocomposite systems. The left and right columns display the first-order sensitivity indices with respect to the model parameters in the first and second branch, respectively.} \label{fig-GSA-FMG}
\end{figure}

The variation of the first-order sensitivity indices for the storage, loss, and magnitude of the complex moduli for the 20\% HSWF nanocomposite systems is shown in Figure \ref{fig-GSA-FMG}, with the remaining data illustrated in Figures S9 and S10. It should be noted that while the variation of the first-order sensitivity indices is important, their maximum value is a decisive factor for the prioritization and fixing tasks as \(S_i \approx 0\) is a necessary and sufficient condition for a model parameter \(q_i\) to be considered non-influential.

\begin{figure}[H]
    \centering
    \begin{subfigure}{0.45\textwidth}
        \includegraphics[width=\linewidth, trim=3cm 0.5cm 5.5cm 1.5cm, clip]{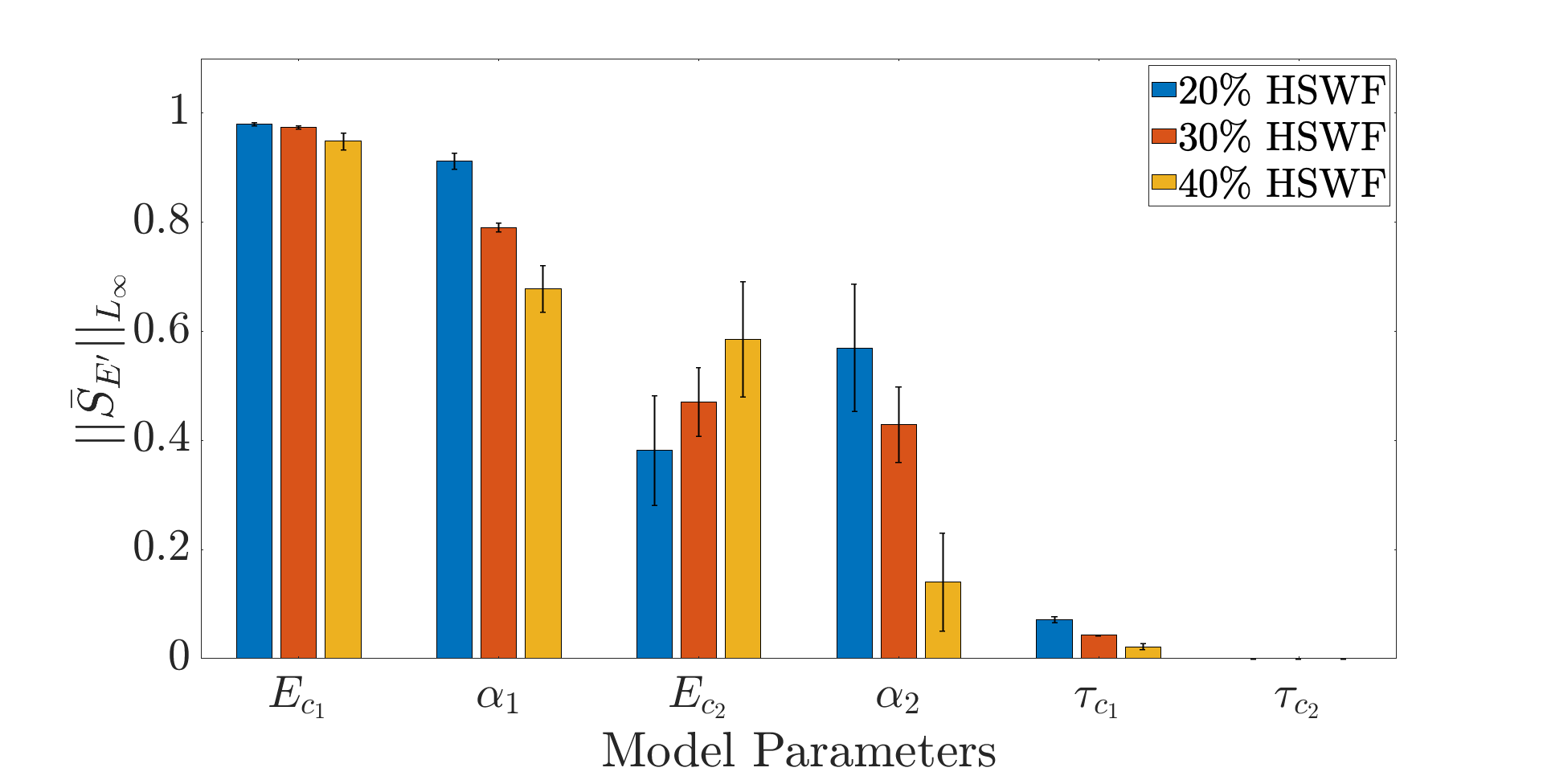}
        \caption{\label{fig-barplots-SEp-Linf-FMG}}
    \end{subfigure}
    \begin{subfigure}{0.45\textwidth}
        \includegraphics[width=\linewidth, trim=3cm 0.5cm 5.5cm 1.5cm, clip]{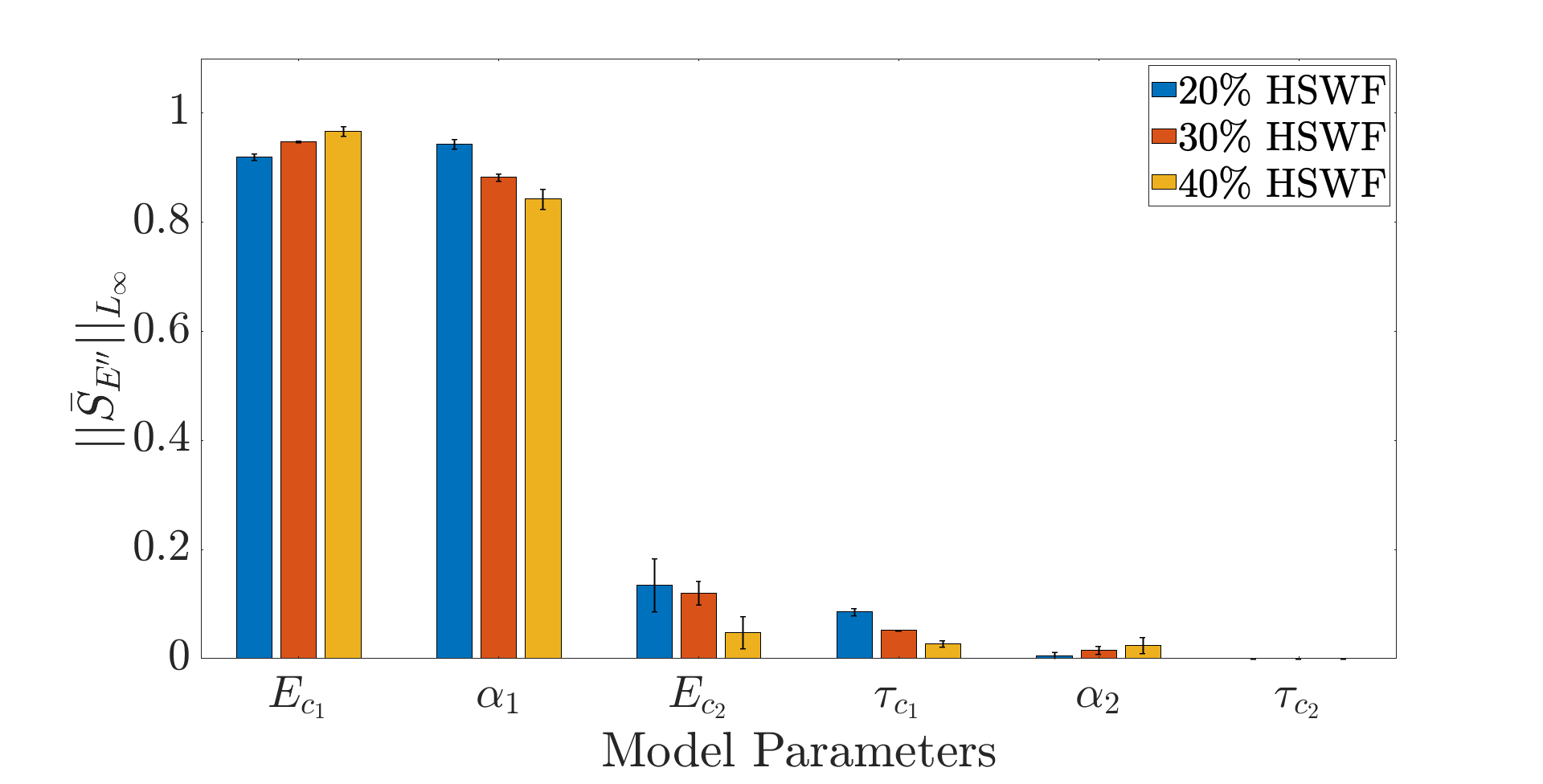}
        \caption{\label{fig-barplots-SEpp-Linf-FMG}}
    \end{subfigure}
    \\
    \begin{subfigure}{0.45\textwidth}
        \includegraphics[width=\linewidth, trim=3cm 0.5cm 5.5cm 1.5cm, clip]{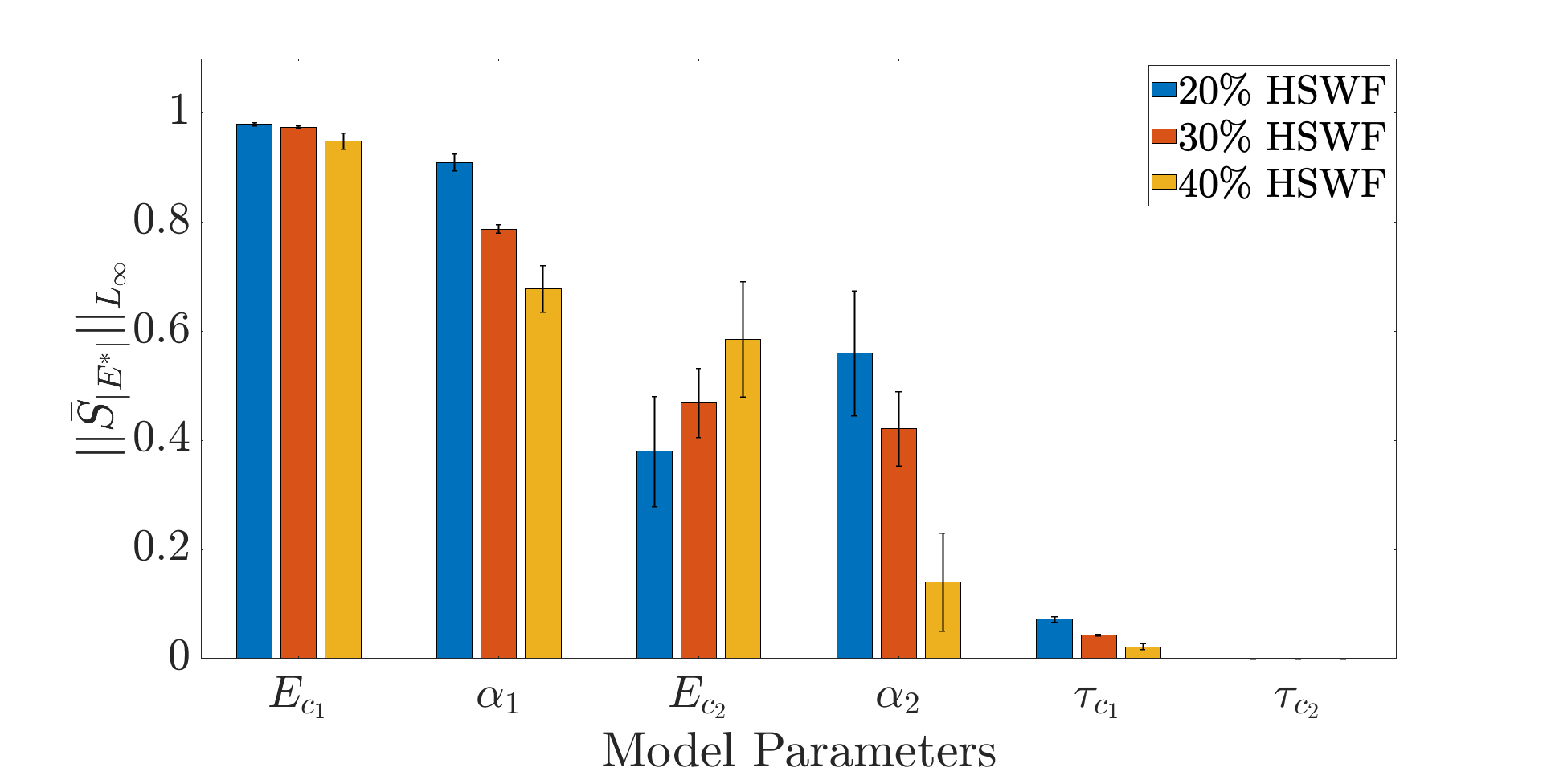}
        \caption{\label{fig-barplots-SEstar-Linf-FMG}}
    \end{subfigure}
    \caption{Averaged \(\text{L}_{\infty}\)-norm of the first-order sensitivity indices associated with the (a) storage, (b) loss, (c) magnitude of the complex moduli with respect to each FMG-FMG model parameter.} \label{fig-barplots-Linf-FMG}
\end{figure}

Given that the maximum values of the first-order sensitivity indices for all moduli associated with \(\alpha_1\), \(E_{c_1}\), \(E_{c_2}\) are significantly greater than zero over multiple decades of frequencies, these parameters are deemed influential. Furthermore, while the first-order sensitivity index of \(\alpha_2\) for the loss modulus is near zero, its effect on the storage modulus and magnitude of the complex modulus cannot be ignored. This discrepancy was anticipated based on the observations made in the LSA section regarding this parameter. Finally, both parameters \(\tau_{c_2}\), with zero first-order sensitivity index value, and \(\tau_{c_1}\), with a maximum first-order sensitivity index no greater than 0.09, are considered non-influential model parameters for the FMG-FMG model. The negligible effect of \(\tau_{c_2}\) on the variance of the moduli is due to the constraint imposed on this parameter during optimization. The \(\text{L}_{\infty}\) norm of the first-order sensitivity indices for all three nanocomposite systems presented in Tables \ref{tab-Linf-storage-FMG}, \ref{tab-Linf-loss-FMG}, and \ref{tab-Linf-complex-FMG} along with their corresponding bar plots depicted in Figure \ref{fig-barplots-Linf-FMG} further support this conclusion.

\begin{table}[H]
\caption{\(\text{L}_{\infty}\)-norm of the first-order sensitivity indices of the storage modulus for the FMG-FMG model with respect to all model parameters.} \label{tab-Linf-storage-FMG}
    \begin{adjustbox}{max width=\textwidth,center=\textwidth}
    \begin{tabular}{cccccccccc}
        \toprule
        \textbf{HSWF wt.\%} & \textbf{xGnP wt.\%} 
        & \textbf{\(\lVert S_{E^{\prime}, E_{c_1}} \rVert_{L_{\infty}}\)} 
        & \textbf{\(\lVert S_{E^{\prime}, \tau_{c_1}} \rVert_{L_{\infty}}\)} 
        & \textbf{\(\lVert S_{E^{\prime}, \alpha_{1}} \rVert_{L_{\infty}}\)} 
        & \textbf{\(\lVert S_{E^{\prime}, E_{c_2}} \rVert_{L_{\infty}}\)} 
        & \textbf{\(\lVert S_{E^{\prime}, \tau_{c_2}} \rVert_{L_{\infty}}\)} 
        & \textbf{\(\lVert S_{E^{\prime}, \alpha_{2}} \rVert_{L_{\infty}}\)} \\
        \midrule
        \multirow{4}{*}{20} & 0.0 & 0.98 & 0.08 & 0.92 & 0.34 & 0.00 & 0.64 \\
                            & 0.5 & 0.98 & 0.07 & 0.93 & 0.53 & 0.00 & 0.40 \\
                            & 1.0 & 0.98 & 0.06 & 0.89 & 0.30 & 0.00 & 0.63 \\
                            & 1.5 & 0.98 & 0.07 & 0.91 & 0.37 & 0.00 & 0.61 \\
        \midrule
        \multirow{4}{*}{30} & 0.0 & 0.98 & 0.04 & 0.79 & 0.55 & 0.00 & 0.35 \\
                            & 0.5 & 0.97 & 0.04 & 0.80 & 0.46 & 0.00 & 0.42 \\
                            & 1.0 & 0.97 & 0.04 & 0.79 & 0.47 & 0.00 & 0.43 \\
                            & 1.5 & 0.97 & 0.04 & 0.78 & 0.40 & 0.00 & 0.52 \\
        \midrule
        \multirow{4}{*}{40} & 0.0 & 0.94 & 0.03 & 0.67 & 0.68 & 0.00 & 0.16 \\
                            & 0.5 & 0.94 & 0.03 & 0.64 & 0.66 & 0.00 & 0.17 \\
                            & 1.0 & 0.94 & 0.02 & 0.66 & 0.54 & 0.00 & 0.22 \\
                            & 1.5 & 0.97 & 0.02 & 0.74 & 0.45 & 0.00 & 0.01 \\
        \bottomrule
    \end{tabular}
    \end{adjustbox}
\end{table}

\FloatBarrier
\newgeometry{top=3cm, bottom=2cm}
\begin{table}[H]
\caption{\(\text{L}_{\infty}\)-norm of first-order sensitivity indices of the loss modulus for the FMG-FMG model with respect to all model parameters.} \label{tab-Linf-loss-FMG}
    \begin{adjustbox}{max width=\textwidth,center=\textwidth}
    \begin{tabular}{cccccccccc}
        \toprule
        \textbf{HSWF wt.\%} & \textbf{xGnP wt.\%} 
        & \textbf{\(\lVert S_{E^{\prime\prime}, E_{c_1}} \rVert_{L_{\infty}}\)} 
        & \textbf{\(\lVert S_{E^{\prime\prime}, \tau_{c_1}} \rVert_{L_{\infty}}\)} 
        & \textbf{\(\lVert S_{E^{\prime\prime}, \alpha_{1}} \rVert_{L_{\infty}}\)} 
        & \textbf{\(\lVert S_{E^{\prime\prime}, E_{c_2}} \rVert_{L_{\infty}}\)} 
        & \textbf{\(\lVert S_{E^{\prime\prime}, \tau_{c_2}} \rVert_{L_{\infty}}\)} 
        & \textbf{\(\lVert S_{E^{\prime\prime}, \alpha_{2}} \rVert_{L_{\infty}}\)} \\
        \midrule
        \multirow{4}{*}{20} & 0.0 & 0.92 & 0.09 & 0.95 & 0.17 & 0.00 & 0.00 \\
                            & 0.5 & 0.92 & 0.09 & 0.95 & 0.07 & 0.00 & 0.01 \\
                            & 1.0 & 0.93 & 0.08 & 0.93 & 0.13 & 0.00 & 0.00 \\
                            & 1.5 & 0.92 & 0.09 & 0.94 & 0.17 & 0.00 & 0.00 \\
        \midrule
        \multirow{4}{*}{30} & 0.0 & 0.95 & 0.05 & 0.89 & 0.10 & 0.00 & 0.03 \\
                            & 0.5 & 0.95 & 0.05 & 0.89 & 0.11 & 0.00 & 0.01 \\
                            & 1.0 & 0.95 & 0.05 & 0.88 & 0.12 & 0.00 & 0.01 \\
                            & 1.5 & 0.95 & 0.05 & 0.87 & 0.15 & 0.00 & 0.01 \\
        \midrule
        \multirow{4}{*}{40} & 0.0 & 0.96 & 0.03 & 0.85 & 0.06 & 0.00 & 0.03 \\
                            & 0.5 & 0.96 & 0.03 & 0.83 & 0.06 & 0.00 & 0.03 \\
                            & 1.0 & 0.96 & 0.03 & 0.83 & 0.07 & 0.00 & 0.03 \\
                            & 1.5 & 0.98 & 0.02 & 0.87 & 0.00 & 0.00 & 0.00 \\
        \bottomrule
    \end{tabular}
    \end{adjustbox}
\end{table}

\begin{table}[H]
\caption{\(\text{L}_{\infty}\)-norm of the first-order sensitivity indices of the magnitude of the complex modulus for the FMG-FMG model with respect to all model parameters.} \label{tab-Linf-complex-FMG}
    \begin{adjustbox}{max width=\textwidth,center=\textwidth}
    \begin{tabular}{cccccccccc}
        \toprule
        \textbf{HSWF wt.\%} & \textbf{xGnP wt.\%} 
        & \textbf{\(\lVert S_{E^{*}, E_{c_1}} \rVert_{L_{\infty}}\)} 
        & \textbf{\(\lVert S_{E^{*}, \tau_{c_1}} \rVert_{L_{\infty}}\)} 
        & \textbf{\(\lVert S_{E^{*}, \alpha_{1}} \rVert_{L_{\infty}}\)} 
        & \textbf{\(\lVert S_{E^{*}, E_{c_2}} \rVert_{L_{\infty}}\)} 
        & \textbf{\(\lVert S_{E^{*}, \tau_{c_2}} \rVert_{L_{\infty}}\)} 
        & \textbf{\(\lVert S_{E^{*}, \alpha_{2}} \rVert_{L_{\infty}}\)} \\
        \midrule
        \multirow{4}{*}{20} & 0.0 & 0.98 & 0.08 & 0.91 & 0.33 & 0.00 & 0.63 \\
                            & 0.5 & 0.98 & 0.07 & 0.92 & 0.53 & 0.00 & 0.39 \\
                            & 1.0 & 0.98 & 0.06 & 0.89 & 0.30 & 0.00 & 0.62 \\
                            & 1.5 & 0.98 & 0.07 & 0.91 & 0.36 & 0.00 & 0.60 \\
        \midrule
        \multirow{4}{*}{30} & 0.0 & 0.98 & 0.04 & 0.78 & 0.55 & 0.00 & 0.34 \\
                            & 0.5 & 0.97 & 0.04 & 0.80 & 0.46 & 0.00 & 0.41 \\
                            & 1.0 & 0.97 & 0.04 & 0.79 & 0.47 & 0.00 & 0.42 \\
                            & 1.5 & 0.97 & 0.04 & 0.78 & 0.40 & 0.00 & 0.51 \\
        \midrule
        \multirow{4}{*}{40} & 0.0 & 0.94 & 0.03 & 0.66 & 0.68 & 0.00 & 0.15 \\
                            & 0.5 & 0.94 & 0.03 & 0.64 & 0.65 & 0.00 & 0.17 \\
                            & 1.0 & 0.94 & 0.02 & 0.66 & 0.54 & 0.00 & 0.22 \\
                            & 1.5 & 0.97 & 0.02 & 0.74 & 0.45 & 0.00 & 0.01 \\
        \bottomrule
    \end{tabular}
    \end{adjustbox}
\end{table}

\FloatBarrier
\restoregeometry

\FloatBarrier
\newgeometry{bottom=4cm}
\subsubsection{Factor Prioritization and Factor Fixing based on the GSA for the FMM-FMG Model Parameters}
\begin{figure}[H]
    \centering
    \begin{subfigure}{\textwidth}
        \centering
        \includegraphics[width=0.73\linewidth, trim=4cm 0.1cm 4.5cm 1.5cm, clip]{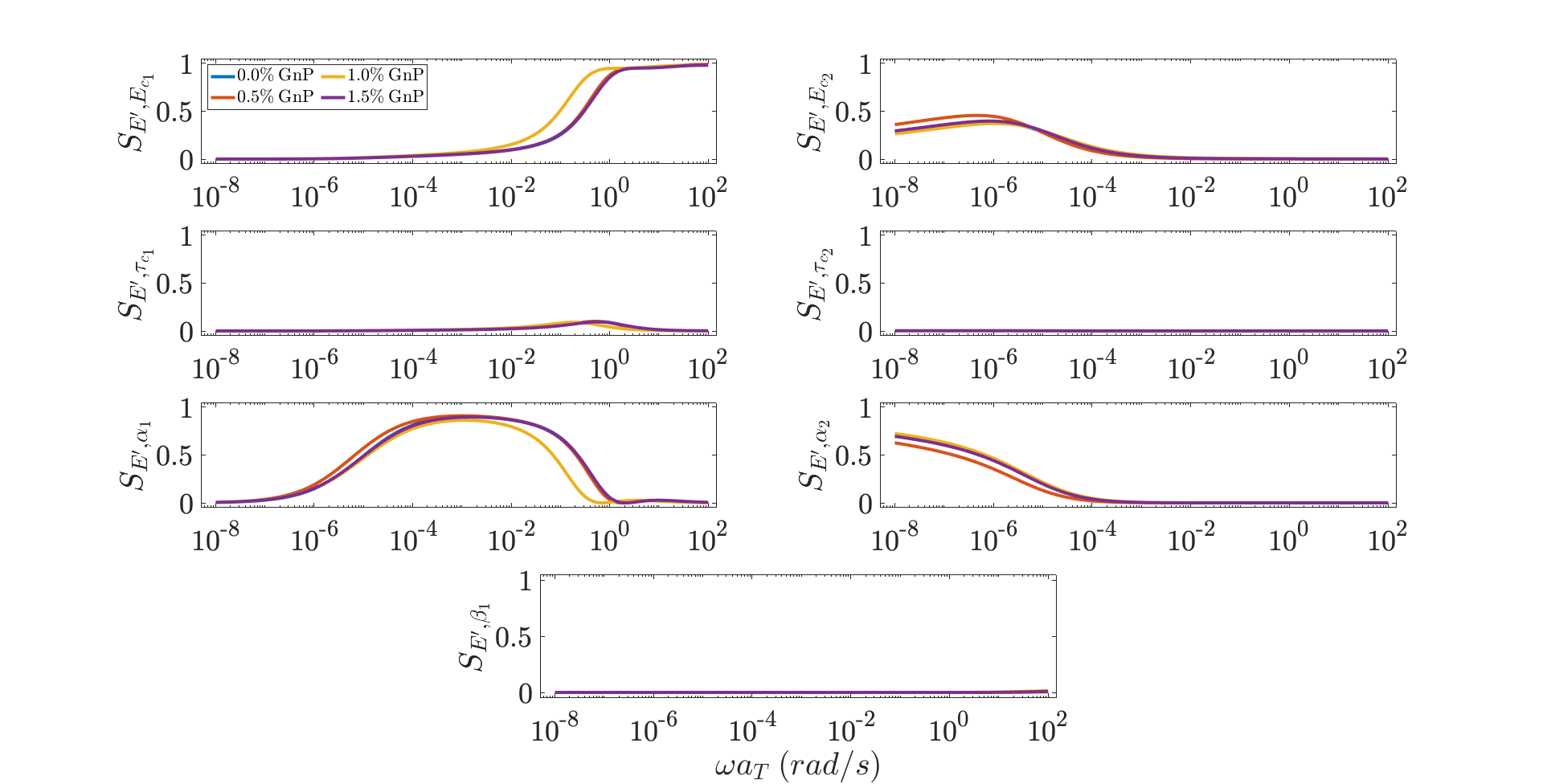}
        \caption{\label{fig-GSA-SEp-FMM}}
    \end{subfigure}
    \begin{subfigure}{\textwidth}
        \centering
        \includegraphics[width=0.73\linewidth, trim=4cm 0.1cm 4.5cm 1.5cm, clip]{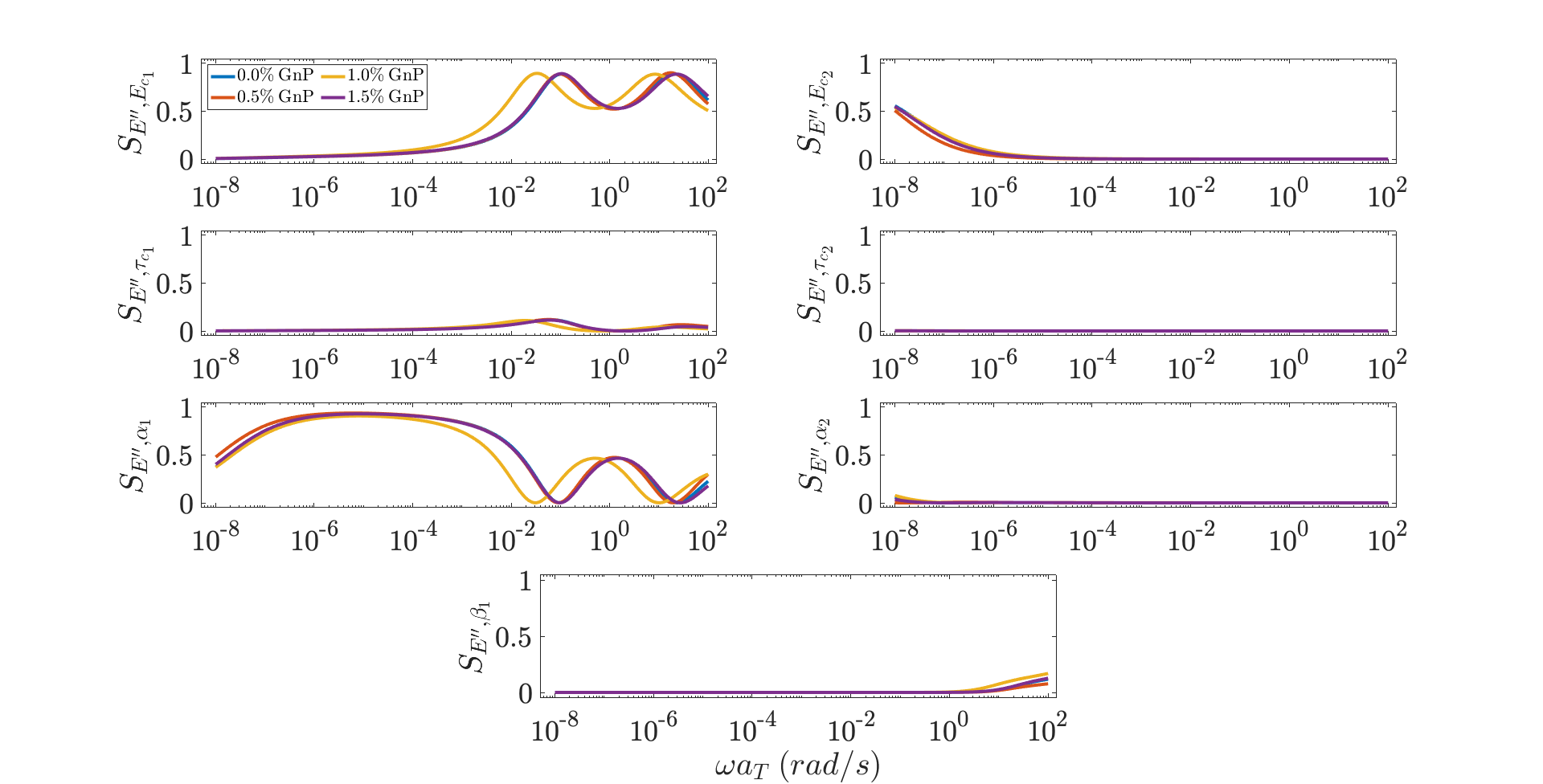}
        \caption{\label{fig-GSA-SEpp-FMM}}
    \end{subfigure}
    \begin{subfigure}{\textwidth}
        \centering
        \includegraphics[width=0.73\linewidth, trim=4cm 0.1cm 4.5cm 1.5cm, clip]{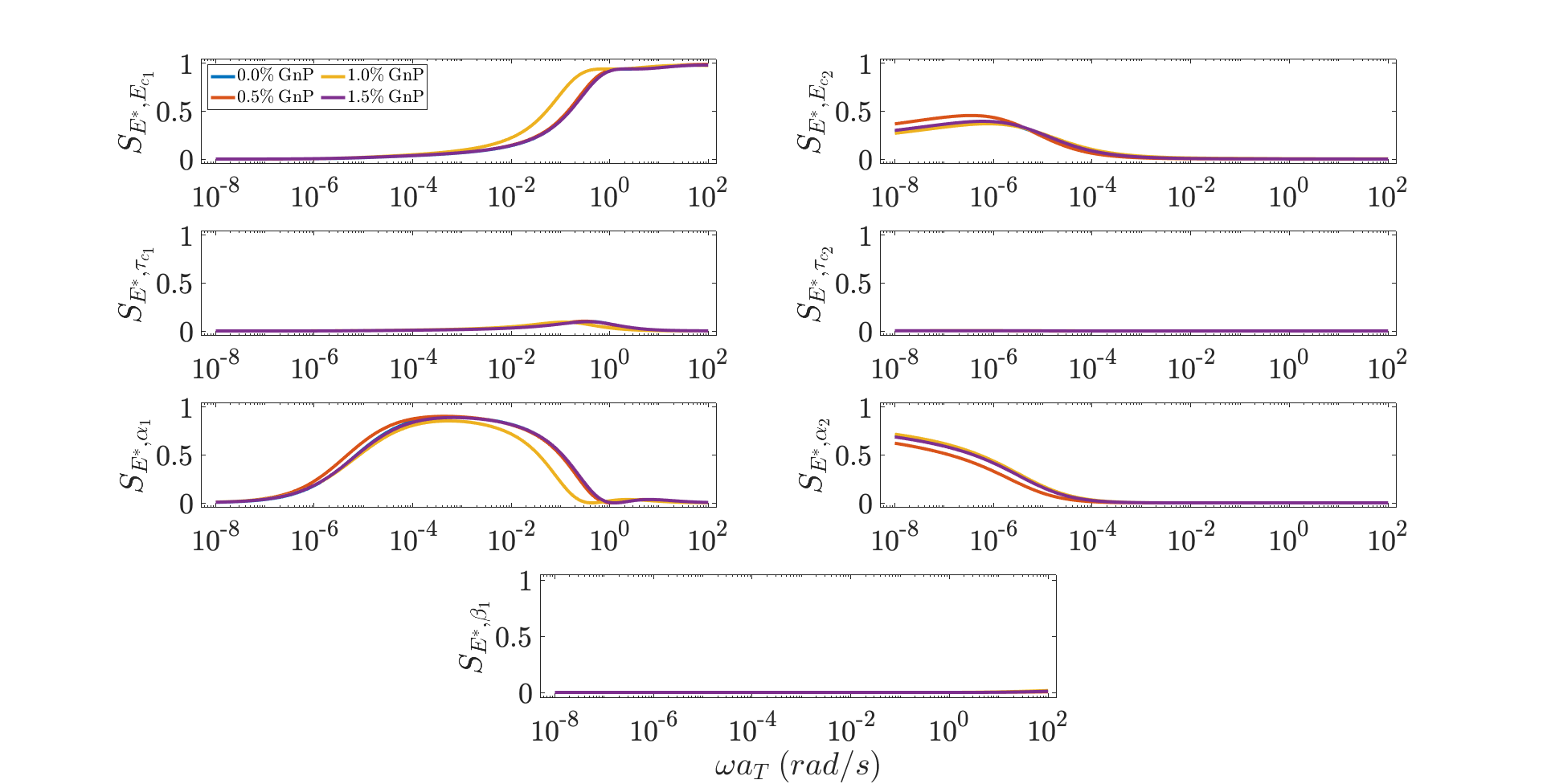}
        \caption{\label{fig-GSA-SEstar-FMM}}
    \end{subfigure}
    \caption{Variation of the first-order sensitivity indices of the (a) storage, (b) loss moduli, and (c) magnitude of the complex modulus of the FMM-FMG model corresponding to the 20\% HSWF nanocomposite systems across decades of shifted frequency. The left and right columns display the first-order sensitivity indices with respect to the model parameters in the first and second branch, respectively, with the last row dedicated to \(\beta_1\).} \label{fig-GSA-FMM}
\end{figure}

In this section, the influential and non-influential FMM-FMG model parameters are also determined in a similar fashion. The first-order sensitivity indices of all three moduli associated with the FMM-FMG model parameters for the 20\% HSWF nanocomposite systems are presented in Figure \ref{fig-GSA-FMM} (The results for the other two nanocomposite systems are provided in Figures S11 and S12). Similar to the FMG-FMG model, the parameters \(E_{c_1}\), \(\alpha_{1}\), and \(E_{c_2}\) exhibit the first-order sensitivity indices much greater than zero for all three moduli, and are thus deemed influential parameters for the FMM-FMG model. Furthermore, \(\tau_{c_2}\) consistently emerges as the least influential model parameter with respect all three moduli, leading to its classification as the least influential parameter.

\begin{figure}[H]
    \centering
    \begin{subfigure}{0.45\textwidth}
        \includegraphics[width=\linewidth, trim=3cm 0.5cm 5.5cm 1.5cm, clip]{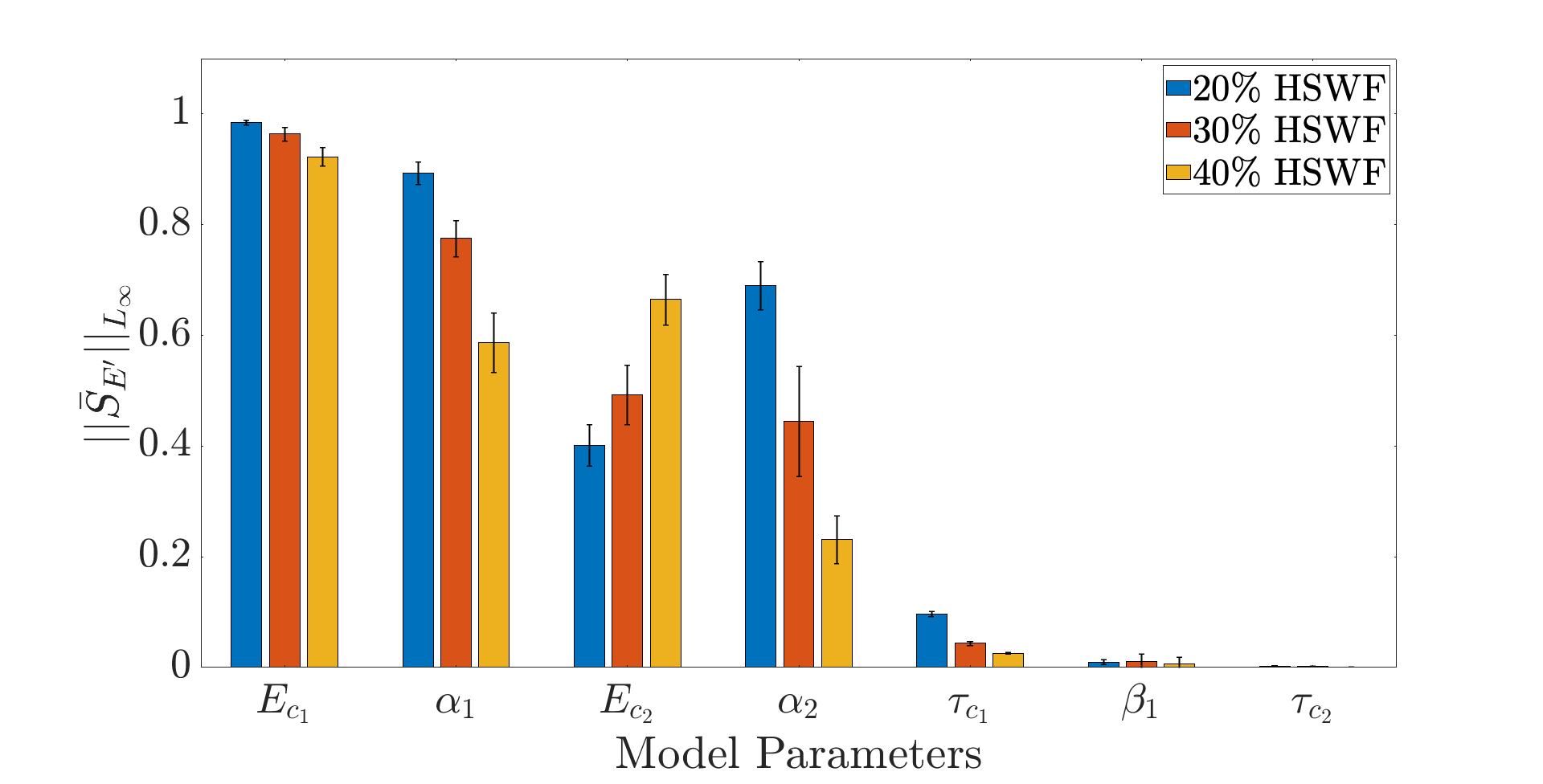}
        \caption{\label{fig-barplots-SEp-Linf-FMM}}
    \end{subfigure}
    \begin{subfigure}{0.45\textwidth}
        \includegraphics[width=\linewidth, trim=3cm 0.5cm 5.5cm 1.5cm, clip]{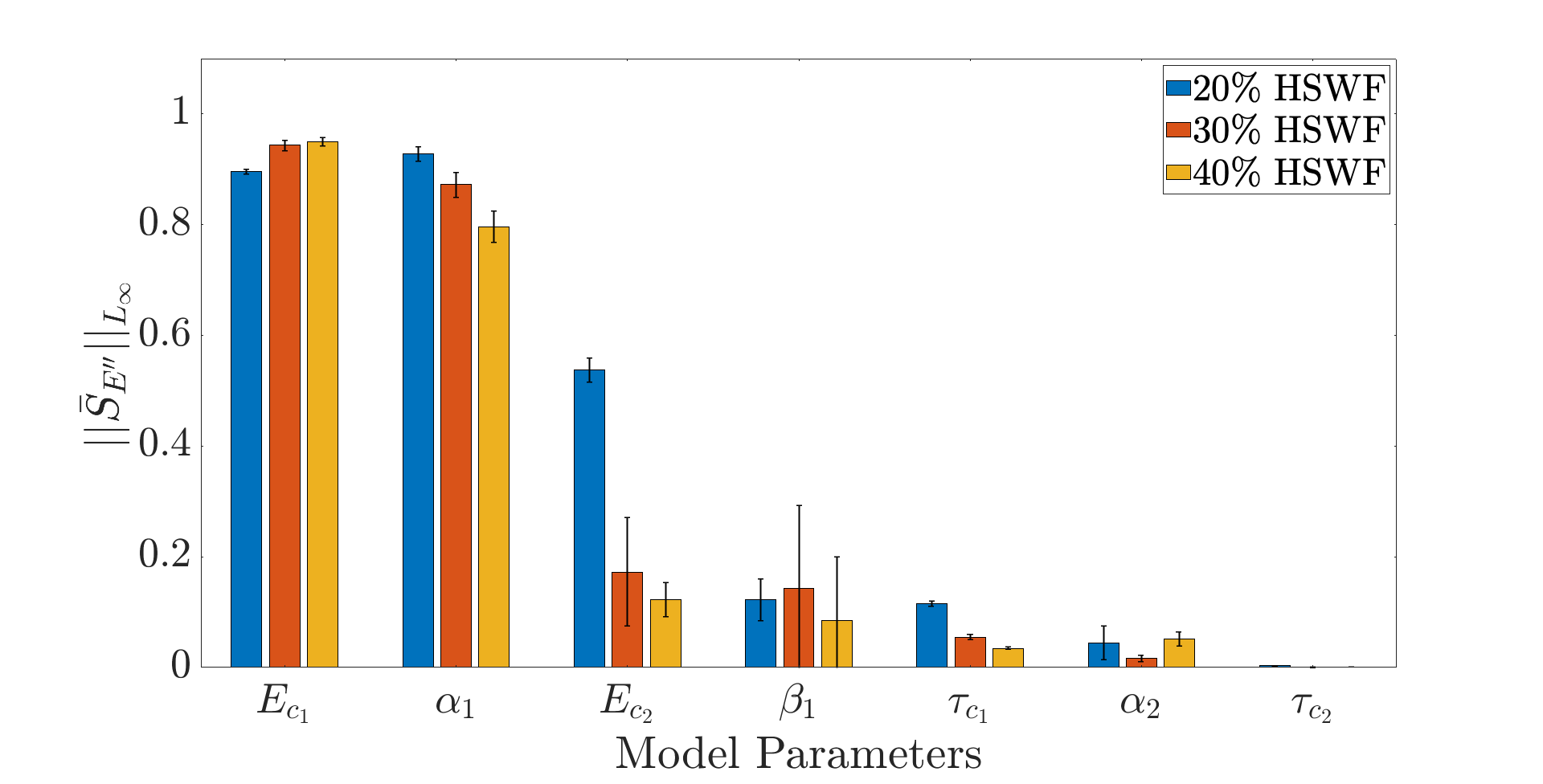}
        \caption{\label{fig-barplots-SEpp-Linf-FMM}}
    \end{subfigure}
    \\
    \begin{subfigure}{0.45\textwidth}
        \includegraphics[width=\linewidth, trim=3cm 0.5cm 5.5cm 1.5cm, clip]{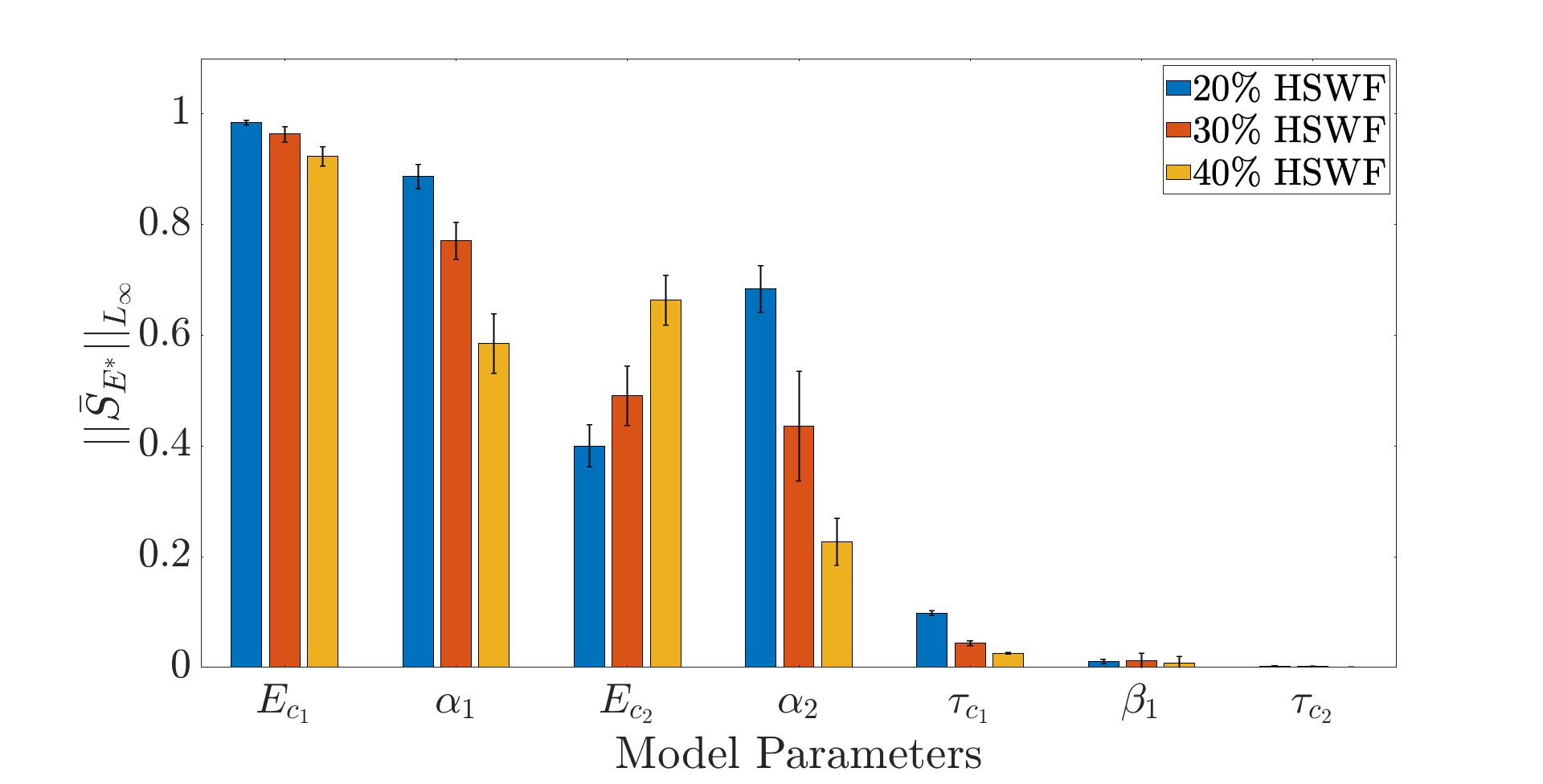}
        \caption{\label{fig-barplots-SEstar-Linf-FMM}}
    \end{subfigure}
    \caption{Averaged \(\text{L}_{\infty}\)-norm of the first-order sensitivity indices associated with the (a) storage, (b) loss, (c) magnitude of the complex moduli with respect to each FMM-FMG model parameter.} \label{fig-barplots-Linf-FMM}
\end{figure}

The remaining parameters, \(\alpha_2\), \(\tau_{c_1}\), and \(\beta_1\), exhibit inconsistent effects on the variability of the storage and loss moduli, a phenomenon also observed in the LSA section. To address this issue, we focus on the effects of these three parameters on the variance of the magnitude of the complex modulus, which reflects the combined effects of these parameters on the variance of both the storage and loss moduli. Since the magnitude of the complex modulus is a strong function of storage modulus, their first-order sensitivity indices become very similar. Consequently, \(\alpha_1\) with its high values of sensitivity indices for both storage and magnitude of the complex moduli cannot be considered as non-influential. But \(\beta_1\) and \(\tau_{c_1}\) with respective maximum values not exceeding 0.01 and 0.1, are regarded as the other two non-influential model parameters for the FMM-FMG model.  The \(\text{L}_{\infty}\) norm of the first-order sensitivity indices for all three nanocomposite systems presented in Tables \ref{tab-Linf-storage-FMM}, \ref{tab-Linf-loss-FMM}, and \ref{tab-Linf-complex-FMM} along with their corresponding bar plots depicted in Figure \ref{fig-barplots-Linf-FMM} further support this conclusion.

\begin{table}[H]
\caption{\(\text{L}_{\infty}\)-norm of the first-order sensitivity indices of the storage modulus for the FMM-FMG model with respect to all model parameters.} \label{tab-Linf-storage-FMM}
    \begin{adjustbox}{max width=\textwidth,center=\textwidth}
    \begin{tabular}{cccccccccc}
        \toprule
        \textbf{HSWF wt.\%} & \textbf{xGnP wt.\%} 
        & \textbf{\(\lVert S_{E^{\prime}, E_{c_1}} \rVert_{L_{\infty}}\)}
        & \textbf{\(\lVert S_{E^{\prime}, \tau_{c_1}} \rVert_{L_{\infty}}\)}
        & \textbf{\(\lVert S_{E^{\prime}, \alpha_{1}} \rVert_{L_{\infty}}\)}
        & \textbf{\(\lVert S_{E^{\prime}, \beta_{1}} \rVert_{L_{\infty}}\)}
        & \textbf{\(\lVert S_{E^{\prime}, E_{c_2}} \rVert_{L_{\infty}}\)}
        & \textbf{\(\lVert S_{E^{\prime}, \tau_{c_2}} \rVert_{L_{\infty}}\)}
        & \textbf{\(\lVert S_{E^{\prime}, \alpha_{2}} \rVert_{L_{\infty}}\)} \\
        \midrule
        \multirow{4}{*}{20} & 0.0 & 0.98 & 0.10 & 0.90 & 0.01 & 0.38 & 0.00 & 0.71 \\
                            & 0.5 & 0.99 & 0.10 & 0.91 & 0.01 & 0.46 & 0.00 & 0.63 \\
                            & 1.0 & 0.98 & 0.09 & 0.86 & 0.02 & 0.37 & 0.00 & 0.72 \\
                            & 1.5 & 0.98 & 0.10 & 0.90 & 0.01 & 0.40 & 0.00 & 0.70 \\
                            \midrule
        \multirow{4}{*}{30} & 0.0 & 0.96 & 0.04 & 0.79 & 0.01 & 0.57 & 0.00 & 0.34 \\
                            & 0.5 & 0.97 & 0.04 & 0.79 & 0.00 & 0.48 & 0.00 & 0.43 \\
                            & 1.0 & 0.97 & 0.04 & 0.79 & 0.00 & 0.48 & 0.00 & 0.43 \\
                            & 1.5 & 0.95 & 0.05 & 0.73 & 0.03 & 0.44 & 0.00 & 0.58 \\
                            \midrule
        \multirow{4}{*}{40} & 0.0 & 0.93 & 0.03 & 0.61 & 0.00 & 0.68 & 0.00 & 0.25 \\
                            & 0.5 & 0.94 & 0.02 & 0.63 & 0.00 & 0.68 & 0.00 & 0.18 \\
                            & 1.0 & 0.93 & 0.02 & 0.60 & 0.00 & 0.60 & 0.00 & 0.28 \\
                            & 1.5 & 0.90 & 0.02 & 0.51 & 0.02 & 0.70 & 0.00 & 0.21 \\
                            \bottomrule
    \end{tabular}
    \end{adjustbox}
\end{table}

\begin{table}[H]
\caption{\(\text{L}_{\infty}\)-norm of the first-order sensitivity indices of the loss modulus for the FMM-FMG model with respect to all model parameters.} \label{tab-Linf-loss-FMM}
    \begin{adjustbox}{max width=\textwidth,center=\textwidth}
    \begin{tabular}{cccccccccc}
        \toprule
        \textbf{HSWF wt.\%} & \textbf{xGnP wt.\%} 
        & \textbf{\(\lVert S_{E^{\prime\prime}, E_{c_1}} \rVert_{L_{\infty}}\)}
        & \textbf{\(\lVert S_{E^{\prime\prime}, \tau_{c_1}} \rVert_{L_{\infty}}\)}
        & \textbf{\(\lVert S_{E^{\prime\prime}, \alpha_{1}} \rVert_{L_{\infty}}\)}
        & \textbf{\(\lVert S_{E^{\prime\prime}, \beta_{1}} \rVert_{L_{\infty}}\)}
        & \textbf{\(\lVert S_{E^{\prime\prime}, E_{c_2}} \rVert_{L_{\infty}}\)}
        & \textbf{\(\lVert S_{E^{\prime\prime}, \tau_{c_2}} \rVert_{L_{\infty}}\)}
        & \textbf{\(\lVert S_{E^{\prime\prime}, \alpha_{2}} \rVert_{L_{\infty}}\)} \\
        \midrule
        \multirow{4}{*}{20} & 0.0 & 0.89 & 0.12 & 0.93 & 0.12 & 0.56 & 0.00 & 0.06 \\
                            & 0.5 & 0.90 & 0.12 & 0.94 & 0.08 & 0.51 & 0.00 & 0.01 \\
                            & 1.0 & 0.90 & 0.11 & 0.91 & 0.17 & 0.54 & 0.00 & 0.08 \\
                            & 1.5 & 0.89 & 0.11 & 0.93 & 0.13 & 0.55 & 0.00 & 0.04 \\
        \midrule
        \multirow{4}{*}{30} & 0.0 & 0.95 & 0.05 & 0.89 & 0.27 & 0.11 & 0.00 & 0.03 \\
                            & 0.5 & 0.95 & 0.05 & 0.88 & 0.03 & 0.14 & 0.00 & 0.01 \\
                            & 1.0 & 0.95 & 0.05 & 0.88 & 0.00 & 0.13 & 0.00 & 0.01 \\
                            & 1.5 & 0.93 & 0.06 & 0.84 & 0.27 & 0.32 & 0.00 & 0.01 \\
        \midrule
        \multirow{4}{*}{40} & 0.0 & 0.95 & 0.04 & 0.81 & 0.05 & 0.15 & 0.00 & 0.06 \\
                            & 0.5 & 0.96 & 0.03 & 0.82 & 0.00 & 0.08 & 0.00 & 0.04 \\
                            & 1.0 & 0.95 & 0.03 & 0.79 & 0.03 & 0.13 & 0.00 & 0.04 \\
                            & 1.5 & 0.94 & 0.04 & 0.76 & 0.25 & 0.14 & 0.00 & 0.07 \\
        \bottomrule
    \end{tabular}
    \end{adjustbox}
\end{table}
\FloatBarrier
\restoregeometry

\FloatBarrier
\newgeometry{top=3cm, bottom=3cm}
\FloatBarrier
\restoregeometry

\begin{table}[H]
\caption{\(\text{L}_{\infty}\)-norm of the first-order sensitivity indices of the magnitude of the complex modulus for the FMM-FMG model with respect to all model parameters.} \label{tab-Linf-complex-FMM}
    \begin{adjustbox}{max width=\textwidth,center=\textwidth}
    \begin{tabular}{cccccccccc}
        \toprule
        \textbf{HSWF wt.\%} & \textbf{xGnP wt.\%} 
        & \textbf{\(\lVert \bar{S}_{E^{*}, E_{c_1}} \rVert_{L_{\infty}}\)}
        & \textbf{\(\lVert \bar{S}_{E^{*}, \tau_{c_1}} \rVert_{L_{\infty}}\)}
        & \textbf{\(\lVert \bar{S}_{E^{*}, \alpha_{1}} \rVert_{L_{\infty}}\)}
        & \textbf{\(\lVert \bar{S}_{E^{*}, \beta_{1}} \rVert_{L_{\infty}}\)}
        & \textbf{\(\lVert \bar{S}_{E^{*}, E_{c_2}} \rVert_{L_{\infty}}\)}
        & \textbf{\(\lVert \bar{S}_{E^{*}, \tau_{c_2}} \rVert_{L_{\infty}}\)}
        & \textbf{\(\lVert \bar{S}_{E^{*}, \alpha_{2}} \rVert_{L_{\infty}}\)} \\
        \midrule
        \multirow{4}{*}{20} & 0.0 & 0.99 & 0.10 & 0.90 & 0.01 & 0.38 & 0.00 & 0.71 \\
                            & 0.5 & 0.99 & 0.10 & 0.90 & 0.01 & 0.46 & 0.00 & 0.62 \\
                            & 1.0 & 0.98 & 0.09 & 0.86 & 0.02 & 0.37 & 0.00 & 0.72 \\
                            & 1.5 & 0.98 & 0.10 & 0.89 & 0.01 & 0.39 & 0.00 & 0.69 \\
        \midrule
        \multirow{4}{*}{30} & 0.0 & 0.96 & 0.04 & 0.79 & 0.02 & 0.57 & 0.00 & 0.33 \\
                            & 0.5 & 0.97 & 0.04 & 0.79 & 0.00 & 0.48 & 0.00 & 0.43 \\
                            & 1.0 & 0.97 & 0.04 & 0.79 & 0.00 & 0.48 & 0.00 & 0.42 \\
                            & 1.5 & 0.94 & 0.05 & 0.72 & 0.03 & 0.44 & 0.00 & 0.57 \\
        \midrule
        \multirow{4}{*}{40} & 0.0 & 0.93 & 0.03 & 0.61 & 0.00 & 0.68 & 0.00 & 0.25 \\
                            & 0.5 & 0.94 & 0.02 & 0.63 & 0.00 & 0.68 & 0.00 & 0.18 \\
                            & 1.0 & 0.93 & 0.02 & 0.60 & 0.00 & 0.60 & 0.00 & 0.27 \\
                            & 1.5 & 0.90 & 0.02 & 0.51 & 0.03 & 0.70 & 0.00 & 0.21 \\
        \bottomrule
    \end{tabular}
    \end{adjustbox}
\end{table}

To summarize our findings in search of the influential and non-influential parameters, Table \ref{tab-prio-GSA} presents the most and least significant parameters in regard to their impact on the variance of the output of both models. Even though both the storage and loss moduli are the quantities of interest in our analysis, the assessment of the parameters importance is made based on their effects on the variability of the complex modulus as it includes the resultant of these separate effects.

\begin{table}[H]
\caption{The least and most influential model parameters for both FMG-FMG and FMM-FMG models based on the GSA.} \label{tab-prio-GSA}
    \centering
    \begin{adjustbox}{max width=\textwidth,center=\textwidth}
        \begin{tabular}{ccc}
            \toprule
            \textbf{Model} & \textbf{Least Influential Parameters} & \textbf{Most Influential Parameters} \\
            \midrule
            FMG-FMG & \makecell{\(\tau_{c_2}\) \\ \(\tau_{c_1}\)} & \makecell{\(E_{c_1}\) \\ \(\alpha_1\)} \\
            \midrule
            FMM-FMG & \makecell{\(\tau_{c_2}\) \\ \(\beta_1\) \\ \(\tau_{c_1}\)} & \makecell{\(E_{c_1}\) \\ \(\alpha_1\)} \\
            \bottomrule
        \end{tabular}
    \end{adjustbox}
\end{table}

A brief remark on comparing the LSA and GSA approaches is worth mentioning here. The impact of model parameters on the output, as measured from a local perspective, often aligns more closely with physical intuition. However, from a variance-based global sensitivity perspective, as suggested by the definitions of the first-order and total-order sensitivity indices, all parameters vary simultaneously, and conditional expectations and variances are used to capture the effects of these variations on the output. Consequently, it is nearly impossible to fully comprehend the effects of these simultaneous variations and their potential interactions on the model’s output. As a result, the influential and non-influential parameters identified through global sensitivity analysis might sometimes appear counter-intuitive.

\section{Conclusions} \label{sec.con}
In this study, the linear viscoelasticity of the neat polyurea (PUa) elastomers with 20\%, 30\%, and 40\% HSWF, and their corresponding nanocomposites with the addition of 0.5\%, 1.0\%, and 1.5\% weight fractions of exfoliated graphene nanoplatelets (xGnP) were modeled using a fractional constitutive model. This new model, consisting of parallel FMM and FMG branches corresponding to the soft and hard segments respectively, perfectly captures the storage and loss moduli master curves of each sample across ten decades of shifted frequency. While the FMM-FMG model performed similar to our earlier FMG-FMG model, the new proposed model revealed a more well-behaved response at the extreme limits of data, leading to a relatively faster optimization convergences, hence better integrity of the complex material throughout. Also, among the seven FMM-FMG model parameters, only the characteristic modulus of the FMG branch (\(E_{c_2}\)) exhibit a dependence on the weight fraction of xGnPs.

Furthermore, we performed a thorough derivative-based local and variance-based global sensitivity analyses seeking for the influential and non-influential parameters with respect to the storage, loss, and complex moduli for both FMG-FMG and FMM-FMG models. Introducing quantitative measures based on the L1-, L2-, and \(\text{L}_{\infty}\)-norms, we successfully accounted for the variation of the local and global sensitivity indices over multiple decades of frequency and proposed a priority list to determine the most and least influential model parameters. This decision was primarily based on the effect of model parameters on the complex modulus as it reflects parameters effects on both the storage and loss moduli. \(E_{c_1}\) and \(\alpha_1\) were found to be the most influential parameters for both models, while \(\tau_{c_2}\) and \(\tau_{c_1}\) were recognized as the least influential parameters for both models, with the addition of \(\beta_1\) as the third non-influential FMM-FMG model parameter. These non-influential model parameters can be treated deterministically, leading to a reduction in the dimensionality of the problem and potentially benefiting other aspects of model analysis, such as uncertainty quantification.

Finally, the proposed robust and well-behaved constitutive model, applicable over ten decades of frequency, can be confidently utilized to optimize the performance of PUa and PUa nanocomposites in their wide range of applications, from biomedical to defense industries.

\section*{CRediT authorship contribution statement}
\textbf{Arman Khoshnevis:} Methodology, Software, Validation, Formal analysis, Investigation, Data curation, Writing – Original draft, Writing – Review \& editing, Visualization. \textbf{Demetrios A Tzelepis:} Conceptualization, Resources, Data curation, Project administration, and Writing – Review \& editing.
\textbf{Valeriy V Ginzburg:} Conceptualization, Methodology, Investigation, Writing - review and editing.
\textbf{Mohsen Zayernouri:} Conceptualization, Methodology, Software, Validation, Formal analysis, Investigation, Resources, Data curation, Writing – Original draft, Writing – Review \& editing, Visualization, Supervision, Project administration, Funding acquisition.

\section*{Declaration of competing interest}
The authors declare that they have no known competing financial interests or personal relationships that could have appeared to influence the work reported in this paper.

\section*{Data availability}
Data will be made available on request.

\section*{Acknowledgments}
This research was funded by the ARO Young Investigator Program (YIP) award (W911NF-19-1-0444) and the NSF award (DMS-1923201).

\appendix

\section{Sobol' Indices} \label{appendix-A}
\setcounter{table}{0}
\setcounter{figure}{0}
This appendix provides a detailed derivation of Sobol' indices \cite{sobol_sensitivity_1993, sobol_global_2001, saltelli_variance_2010, smith_uncertainty_2013, yang2012adaptive, zhang2014enabling}, which are variance-based global sensitivity measures applicable to general non-linear models. Assuming, without loss of generality,  \(y=f(\mathbf{q})\) is a square integrable function and each of its model parameters, \(q_i\), is independently and uniformly distributed over the interval \([0, 1]\), that is \(q_i \sim \mathcal{U}(a_i, b_i)\). Furthermore, let \(\Omega^k = [0, 1]^k\) represent a \(k\)-dimensional unit hyper cube, over which the function \(f\) is defined. The ANOVA employs the following hierarchical expansion
\begin{equation} \label{eq-decomp}
\begin{aligned}
f(\mathbf{q}) &= f_0 + \sum_{i=1}^{k} f_i(q_i) + \sum_{1\leq i < j \leq k} f_{i,j}(q_i,q_j) + \dots \\
&+ \sum_{1\leq i_1 < \dots < i_s \leq k} f_{i_1, \dots, i_s}(q_{i_1}, \dots, q_{i_s}) + \dots + f_{1, 2, \dots, k}(q_1, \dots, q_k).
\end{aligned}
\end{equation}
To ensure the uniqueness of this decomposition, it is required that each term in the expansion has a zero mean, and each pair of terms is orthogonal. The decomposition outlined in \eqref{eq-decomp} consists of \(2^k\) terms in total, where the zeroth-, first-, and second-order terms (functions) can be derived as follows:
\begin{subequations}
\begin{align}
f_0 &= \int_{\Omega} f(\mathbf{q}) \,d\mathbf{q} = \mathbb{E}(y), \\
f_i(q_i) &= \int_{\Omega^{k-1}} f(\mathbf{q}) \, d\mathbf{q}_{\sim i} - f_0 = \mathbb{E}_{\mathbf{q}_{\sim i}}(y|q_i) - f_0, \\
f_{i,j}(q_i,q_j) &= \int_{\Omega^{k-2}} f(\mathbf{q}) \, d\mathbf{q}_{\sim i,j} - f_i(q_i) - f_j(q_j) - f_0, \nonumber \\
&= \mathbb{E}_{\mathbf{q}_{\sim i,j}}(y|q_i, q_j) - f_i(q_i) - f_j(q_j) - f_0.
\end{align}
\end{subequations}
Higher order terms are derived in a similar manner. It should be noted that the probability distribution function of each factor, \(q_i\), is embedded in the function \(f_i(q_i)\) without loss of generality.

Similar to the decomposition employed in \eqref{eq-decomp}, the total variance of \(y\) can also be decomposed as well
\begin{equation}
\begin{aligned}
\mathbb{V}\left(y\right) = \sum_{i=1}^{k} \mathbb{V}_{q_{i}}\left(y\right) + \sum_{1 \leq i < j \leq k}^{k} \mathbb{V}_{q_i, q_j}\left(y\right) + \dots + \mathbb{V}_{q_1, \dots, q_k}\left(y\right).
\end{aligned}
\end{equation}
By dividing both sides of this equation by \(\mathbb{V}\left(y\right)\), the first- and higher-order sensitivity indices are obtained
\begin{equation} \label{eq-sums}
\begin{aligned}
\sum_{i=1}^{k} S_i + \sum_{1 \leq i < j \leq k}^{k} S_{i,j} + \dots + S_{1, 2, \dots, k} = 1
\end{aligned}
\end{equation}
Hence, the first-order index, \(S_i\), which quantifies the main effect contribution of each input factor to the variance of the output, is calculated as follows:
\begin{equation}
\begin{aligned}
S_i = \frac{\mathbb{V}\left(f_i(q_i)\right)}{\mathbb{V}\left(y\right)} = \frac{\mathbb{V}_{q_i}\left(\mathbb{E}_{\mathbf{q}_{\sim i}}(y|q_i)\right)}{\mathbb{V}\left(y\right)}.
\end{aligned}
\end{equation}
Due to the computational limitations, it is customary to define the so-called total-order sensitivity index, \(S_{Ti}\), which involves the firs-order effect and all the higher-order interaction effects with the rest of the model parameters on the variance of the output. This index can be expressed as
\begin{equation}
\begin{aligned}
S_{Ti} &= S_i + \sum_{\substack{j=1 \\ j\neq i}}^{k} S_{i,j} + \sum_{\substack{1 \leq j < l \leq k \\ j\neq i, \,l \neq i}} S_{i,j,l} + \dots + S_{1, \dots, k}.
\end{aligned}
\end{equation}
This sensitivity index can be obtained by excluding all terms of any order that do no contain the factor \(q_i\) from the summations in \eqref{eq-sums}. Consequently, the total-order sensitivity index can be defined as
\begin{equation}
\begin{aligned}
S_{Ti} = 1 - \frac{\mathbb{V}_{\mathbf{q}_{\sim i}}\left( \mathbb{E}_{q_i}(y|\mathbf{q}_{\sim i}) \right)}{\mathbb{V}\left(y\right)}.
\end{aligned}
\end{equation}
Alternatively, it is possible to obtain this index using the identity introduced in \eqref{eq-identitiy}, but this time, conditioning on all factors except \(q_i\) as \(\mathbb{V}\left(y\right) = \mathbb{V}_{\mathbf{q_{\sim i}}}\left( \mathbb{E}_{q_i}(y|\mathbf{q}_{\sim i}) \right) + \mathbb{E}_{\mathbf{q_{\sim i}}} \left(\mathbb{V}_{q_i}(y|\mathbf{q}_{\sim i}) \right)\). Dividing both sides of this equation by \(\mathbb{V}(y)\) and rearranging the terms, another form of the total-order sensitivity index can be obtained as
\begin{equation}
\begin{aligned}
S_{Ti} = 1 - \frac{\mathbb{V}_{\mathbf{q}_{\sim i}}\left( \mathbb{E}_{q_i}(y|\mathbf{q_{\sim i}}) \right)}{\mathbb{V}\left(y\right)} = \frac{\mathbb{E}_{\mathbf{q_{\sim i}}} \left( \mathbb{V}_{q_i} (y|\mathbf{q}_{\sim i}) \right)}{\mathbb{V}\left(y\right)}.
\end{aligned}
\end{equation}

\section{Estimation of the first- and total-order Sobol' sensitivity indices} \label{appendix-B}
\setcounter{table}{0}
\setcounter{figure}{0}
In this paper, two estimators proposed by Salteli \cite{saltelli_variance_2010} and Jenson \cite{jansen_analysis_1999} are employed to compute the first- and total-order sensitivity indices, respectively. In this appendix, a concise overview of the algorithm implemented for these computations is presented as follows:
\begin{enumerate}
\item Generate two \((N \times k)\) matrices, \(\mathbf{A}\) and \(\mathbf{B}\), where each element \(a_i^{(j)}\) and \(b_i^{(j)}\) is independently drawn from the interval \([0, 1]\) using a low discrepancy quasi-random Sobol' sequence:

\begin{linenomath}
\begin{equation}
\begin{aligned}
\mathbf{A} &= \begin{bmatrix}
a_1^{(1)} & \cdots & a_i^{(1)} & \cdots & a_k^{(1)} \\
\vdots & \ddots & \vdots & \ddots & \vdots \\
a_1^{(N)} & \cdots & a_i^{(1)} & \cdots & a_k^{(N)} \\
\end{bmatrix},
\mathbf{B} &= \begin{bmatrix}
b_1^{(1)} & \cdots & b_i^{(1)} & \cdots & b_k^{(1)} \\
\vdots & \ddots & \vdots & \ddots & \vdots \\
b_1^{(N)} & \cdots & b_i^{(1)} & \cdots & b_k^{(N)} \\
\end{bmatrix}.
\end{aligned}
\end{equation}
\end{linenomath}
Subsequently, each \(a_i^{(j)}\) and \(b_i^{(j)}\) element is mapped to the range of variability of their respective model parameter \(q_i\). Therefore, \(a_i^{(j)}\) and \(b_i^{(j)}\) are realizations of this model parameter.

\item Construct the cross-sampled matrix \(\mathbf{A}_{\mathbf{B}}^{(i)}\) such that all columns are from \(\mathbf{A}\) except for the \(i\)-th column, which is replaced by the \(i\)-th column of \(\mathbf{B}\) as

\begin{equation}
\begin{aligned}
\mathbf{A}_{\mathbf{B}}^{(i)} = \begin{bmatrix}
a_1^{(1)} & \cdots & b_i^{(1)} & \cdots & a_k^{(1)} \\
\vdots & \ddots & \vdots & \ddots & \vdots \\
a_1^{(N)} & \cdots & b_i^{(1)} & \cdots & a_k^{(N)} \\
\end{bmatrix}.
\end{aligned}
\end{equation}

\item Evaluate the model for all input values in the sample matrices \(\mathbf{A}\), \(\mathbf{B}\), and \(\mathbf{A_{\mathbf{B}}}^{(i)}\):
\begin{equation}
\begin{aligned}
y(\mathbf{A}) = f(\mathbf{A})\,, y(\mathbf{B}) = f(\mathbf{B})\,, y(\mathbf{A}_{\mathbf{B}}^{(i)}) = f(\mathbf{A_{\mathbf{B}}}^{(i)}).
\end{aligned}
\end{equation}

\item Compute the total variance \(\mathbb{V}(y)\) by taking the variance from the row-wise concatenated model outputs \(y({\mathbf{C}}) = \begin{bmatrix}
    y(\mathbf{A}) \\
    y(\mathbf{B})
\end{bmatrix}\).

\item Utilize the following estimators to calculate the first-order and total-order sensitivity indices:
\begin{subequations}
\begin{linenomath}
\begin{equation}
S_i \approx \frac{\frac{1}{N} \sum_{j=1}^{N} f(\mathbf{B})_j (f(\mathbf{A_{\mathbf{B}}}^{(i)})_j - f(\mathbf{A})_j }{\mathbb{V}(y(\mathbf{C}))},
\end{equation}
\end{linenomath}

\begin{linenomath}
\begin{equation}
S_{Ti} \approx \frac{\frac{1}{N} \sum_{j=1}^{N} (f(\mathbf{A})_j - f(\mathbf{A_{\mathbf{B}}}^{(i)})_j)}{\mathbb{V}(y(\mathbf{C}))}.
\end{equation}
\end{linenomath}
\end{subequations}
\end{enumerate}

\bibliographystyle{unsrt}  
\bibliography{manuscript}

\begin{thebibliography}{10}

\bibitem{petrovic_polyurethane_1991}
Zoran~S Petrovi{\'c} and James Ferguson.
\newblock Polyurethane elastomers.
\newblock {\em Progress in Polymer Science}, 16(5):695--836, 1991.

\bibitem{szycher_polyurethanes_2013}
M~Szycher.
\newblock Polyurethanes.
\newblock {\em Szycher’S Handbook of Polyurethanes; CRC Press Taylor \& Francis Group: Boca Raton, FL, USA}, pages 1--12, 2013.

\bibitem{akindoyo_polyurethane_2016}
John~O Akindoyo, MdDH Beg, Suriati Ghazali, MR~Islam, Nitthiyah Jeyaratnam, and AR~Yuvaraj.
\newblock Polyurethane types, synthesis and applications--a review.
\newblock {\em Rsc Advances}, 6(115):114453--114482, 2016.

\bibitem{das_brief_2020}
Abhijit Das and Prakash Mahanwar.
\newblock A brief discussion on advances in polyurethane applications.
\newblock {\em Advanced Industrial and Engineering Polymer Research}, 3(3):93--101, 2020.

\bibitem{sonnenschein_polyurethanes_2021}
Mark~F Sonnenschein.
\newblock {\em Polyurethanes: science, technology, markets, and trends}.
\newblock John Wiley \& Sons, 2021.

\bibitem{chattopadhyay_thermal_2005}
DK~Chattopadhyay, B~Sreedhar, and KVSN Raju.
\newblock Thermal stability of chemically crosslinked moisture-cured polyurethane coatings.
\newblock {\em Journal of applied polymer science}, 95(6):1509--1518, 2005.

\bibitem{chattopadhyay_structural_2007}
Dipak~Kumar Chattopadhyay and KVSN Raju.
\newblock Structural engineering of polyurethane coatings for high performance applications.
\newblock {\em Progress in polymer science}, 32(3):352--418, 2007.

\bibitem{driffield_method_2007}
Malcolm Driffield, Emma~L Bradley, and Laurence Castle.
\newblock A method of test for residual isophorone diisocyanate trimer in new polyester-polyurethane coatings on light metal packaging using liquid chromatography with tandem mass spectrometric detection.
\newblock {\em Journal of Chromatography A}, 1141(1):61--66, 2007.

\bibitem{ding2024recent}
Lailong Ding, Yifan Wang, Jiayu Lin, Mingliang Ma, Jinhu Hu, Xishun Qiu, Chao Wu, and Chao Feng.
\newblock Recent advances in polyurea elastomers and their applications in blast protection: a review.
\newblock {\em Journal of Materials Science}, pages 1--31, 2024.

\bibitem{cognard_handbook_2005}
Philippe Cognard.
\newblock {\em Handbook of adhesives and sealants: basic concepts and high tech bonding}.
\newblock Elsevier, 2005.

\bibitem{lee_preparation_2008}
Wen-Jau Lee and Meng-Shinan Lin.
\newblock Preparation and application of polyurethane adhesives made from polyhydric alcohol liquefied taiwan acacia and china fir.
\newblock {\em Journal of applied polymer science}, 109(1):23--31, 2008.

\bibitem{jia-hu_synthesis_2015}
Guo Jia-Hu, Liu Yu-Cun, Chai Tao, Jing Su-Ming, Ma~Hui, Qin Ning, Zhou Hua, Yan Tao, and He~Wei-Ming.
\newblock Synthesis and properties of a nano-silica modified environmentally friendly polyurethane adhesive.
\newblock {\em RSC Advances}, 5(56):44990--44997, 2015.

\bibitem{meier-westhues_polyurethanes_2019}
Hans-Ulrich Meier-Westhues.
\newblock {\em Polyurethanes: coatings, adhesives and sealants}.
\newblock European Coatings, 2019.

\bibitem{bicerano_flexible_2004}
J~Bicerano, RD~Daussin, MJA Elwell, HR~van~der Wal, P~Berthevas, M~Brown, F~Casati, W~Farrissey, J~Fosnaugh, R~de~Genova, et~al.
\newblock Flexible polyurethane foams.
\newblock In {\em Polymeric Foams}, pages 164--229. CRC Press, 2004.

\bibitem{harikrishnan_modeling_2010}
G~Harikrishnan and DV~Khakhar.
\newblock Modeling the dynamics of reactive foaming and film thinning in polyurethane foams.
\newblock {\em AIChE journal}, 56(2):522--530, 2010.

\bibitem{allan_thermal_2013}
Deborah Allan, J~Daly, and JJ~Liggat.
\newblock Thermal volatilisation analysis of tdi-based flexible polyurethane foam.
\newblock {\em Polymer Degradation and Stability}, 98(2):535--541, 2013.

\bibitem{dsouza_polyurethane_2014}
Jason D'Souza, Rafael Camargo, and Ning Yan.
\newblock Polyurethane foams made from liquefied bark-based polyols.
\newblock {\em Journal of Applied Polymer Science}, 131(16), 2014.

\bibitem{furukawa_microphase-separated_2005}
Mutsuhisa Furukawa, Yoshitaka Mitsui, Tomoya Fukumaru, and Ken Kojio.
\newblock Microphase-separated structure and mechanical properties of novel polyurethane elastomers prepared with ether based diisocyanate.
\newblock {\em Polymer}, 46(24):10817--10822, 2005.

\bibitem{bagdi_quantitative_2012}
Krist{\'o}f Bagdi, Kinga Moln{\'a}r, Mih{\'a}ly K{\'a}llay, Peter Sch{\"o}n, Julius~G Vancs{\'o}, and B{\'e}la Puk{\'a}nszky.
\newblock Quantitative estimation of the strength of specific interactions in polyurethane elastomers, and their effect on structure and properties.
\newblock {\em European polymer journal}, 48(11):1854--1865, 2012.

\bibitem{koberstein_smallangle_1983}
Jeffrey~T Koberstein and Richard~S Stein.
\newblock Small-angle x-ray scattering studies of microdomain structure in segmented polyurethane elastomers.
\newblock {\em Journal of Polymer Science: Polymer Physics Edition}, 21(8):1439--1472, 1983.

\bibitem{koberstein_smallangle_1983-1}
Jeffrey~T Koberstein and Richard~S Stein.
\newblock Small-angle x-ray scattering measurements of diffuse phase-boundary thicknesses in segmented polyurethane elastomers.
\newblock {\em Journal of polymer science: polymer physics edition}, 21(10):2181--2200, 1983.

\bibitem{leung_smallangle_1985}
Louis~M Leung and Jeffrey~T Koberstein.
\newblock Small-angle scattering analysis of hard-microdomain structure and microphase mixing in polyurethane elastomers.
\newblock {\em Journal of Polymer Science: Polymer Physics Edition}, 23(9):1883--1913, 1985.

\bibitem{christenson_model_1986}
CP~Christenson, MA~Harthcock, MD~Meadows, HL~Spell, WL~Howard, MW~Creswick, RE~Guerra, and RB~Turner.
\newblock Model mdi/butanediol polyurethanes: molecular structure, morphology, physical and mechanical properties.
\newblock {\em Journal of Polymer Science Part B: Polymer Physics}, 24(7):1401--1439, 1986.

\bibitem{koberstein_compression-molded_1992}
Jeffrey~T Koberstein and Louis~M Leung.
\newblock Compression-molded polyurethane block copolymers. 2. evaluation of microphase compositions.
\newblock {\em Macromolecules}, 25(23):6205--6213, 1992.

\bibitem{koberstein_compression-molded_1992-1}
JT~Koberstein, AF~Galambos, and LM~Leung.
\newblock Compression-molded polyurethane block copolymers. 1. microdomain morphology and thermomechanical properties.
\newblock {\em Macromolecules}, 25(23):6195--6204, 1992.

\bibitem{ginzburg_theoretical_2007}
Valeriy~V Ginzburg, Jozef Bicerano, Christopher~P Christenson, Alan~K Schrock, and Alexander~Z Patashinski.
\newblock Theoretical modeling of the relationship between young's modulus and formulation variables for segmented polyurethanes.
\newblock {\em Journal of Polymer Science Part B: Polymer Physics}, 45(16):2123--2135, 2007.

\bibitem{benoit_scattering_1988}
Hadziioannoua Benoit and G~Hadziioannou.
\newblock Scattering theory and properties of block copolymers with various architectures in the homogeneous bulk state.
\newblock {\em Macromolecules}, 21(5):1449--1464, 1988.

\bibitem{heintz_spectroscopic_2005}
Amy~M Heintz, Daniel~J Duffy, Chad~M Nelson, Ying Hua, Shaw~L Hsu, Wu~Suen, and Charles~W Paul.
\newblock A spectroscopic analysis of the phase evolution in polyurethane foams.
\newblock {\em Macromolecules}, 38(22):9192--9199, 2005.

\bibitem{christenson_relationship_2005}
Elizabeth~M Christenson, James~M Anderson, Anne Hiltner, and Eric Baer.
\newblock Relationship between nanoscale deformation processes and elastic behavior of polyurethane elastomers.
\newblock {\em Polymer}, 46(25):11744--11754, 2005.

\bibitem{garrett_microdomain_2001}
James~T Garrett, Christopher~A Siedlecki, and James Runt.
\newblock Microdomain morphology of poly (urethane urea) multiblock copolymers.
\newblock {\em Macromolecules}, 34(20):7066--7070, 2001.

\bibitem{kaushiva_uniaxial_2000}
BD~Kaushiva and GL~Wilkes.
\newblock Uniaxial orientation behavior and consideration of the geometric anisotropy of polyurea hard domain structure in flexible polyurethane foams.
\newblock {\em Polymer}, 41(18):6987--6991, 2000.

\bibitem{aou_characterization_2013}
Kaoru Aou, Alan~K Schrock, Valeriy~V Ginzburg, and Philip~C Price.
\newblock Characterization of polyurethane hard segment length distribution using soft hydrolysis/maldi and monte carlo simulation.
\newblock {\em Polymer}, 54(18):5005--5015, 2013.

\bibitem{as2023mechanics}
Faisal As’ad and Charbel Farhat.
\newblock A mechanics-informed deep learning framework for data-driven nonlinear viscoelasticity.
\newblock {\em Computer Methods in Applied Mechanics and Engineering}, 417:116463, 2023.

\bibitem{abdolazizi2024viscoelastic}
Kian~P Abdolazizi, Kevin Linka, and Christian~J Cyron.
\newblock Viscoelastic constitutive artificial neural networks (vcanns)--a framework for data-driven anisotropic nonlinear finite viscoelasticity.
\newblock {\em Journal of Computational Physics}, 499:112704, 2024.

\bibitem{suzuki2021data}
Jorge~L Suzuki, Tyler~G Tuttle, Sara Roccabianca, and Mohsen Zayernouri.
\newblock A data-driven memory-dependent modeling framework for anomalous rheology: Application to urinary bladder tissue.
\newblock {\em Fractal and Fractional}, 5(4):223, 2021.

\bibitem{dabiri2023fractional}
Donya Dabiri, Milad Saadat, Deepak Mangal, and Safa Jamali.
\newblock Fractional rheology-informed neural networks for data-driven identification of viscoelastic constitutive models.
\newblock {\em Rheologica Acta}, 62(10):557--568, 2023.

\bibitem{gonzalez2023model}
Ernesto Gonzalez-Saiz and Daniel Garcia-Gonzalez.
\newblock Model-driven identification framework for optimal constitutive modeling from kinematics and rheological arrangement.
\newblock {\em Computer Methods in Applied Mechanics and Engineering}, 415:116211, 2023.

\bibitem{upadhyay2024physics}
Kshitiz Upadhyay, Jan~N Fuhg, Nikolaos Bouklas, and KT~Ramesh.
\newblock Physics-informed data-driven discovery of constitutive models with application to strain-rate-sensitive soft materials.
\newblock {\em Computational Mechanics}, pages 1--30, 2024.

\bibitem{khoshnevis2024double}
Arman Khoshnevis, Ali Ahmadpour, and Ehsan Amani.
\newblock Double emulsion generation in shear-thinning fluids under electric field effects.
\newblock {\em International Journal of Mechanical Sciences}, page 109556, 2024.

\bibitem{doherty2023stabilised}
William Doherty, Timothy~N Phillips, and Zhihua Xie.
\newblock A stabilised finite element framework for viscoelastic multiphase flows using a conservative level-set method.
\newblock {\em Journal of Computational Physics}, 477:111936, 2023.

\bibitem{gruninger2024benchmarking}
Cole Gruninger, Aaron Barrett, Fuhui Fang, M~Gregory Forest, and Boyce~E Griffith.
\newblock Benchmarking the immersed boundary method for viscoelastic flows.
\newblock {\em Journal of Computational Physics}, 506:112888, 2024.

\bibitem{saltelli1999sensitivity}
Andrea Saltelli.
\newblock Sensitivity analysis: Could better methods be used?
\newblock {\em Journal of Geophysical Research: Atmospheres}, 104(D3):3789--3793, 1999.

\bibitem{morris1991factorial}
Max~D Morris.
\newblock Factorial sampling plans for preliminary computational experiments.
\newblock {\em Technometrics}, 33(2):161--174, 1991.

\bibitem{sobol1990sensitivity}
Il'ya~Meerovich Sobol'.
\newblock On sensitivity estimation for nonlinear mathematical models.
\newblock {\em Matematicheskoe modelirovanie}, 2(1):112--118, 1990.

\bibitem{saltelli1999quantitative}
Andrea Saltelli, Stefano Tarantola, and KP-S Chan.
\newblock A quantitative model-independent method for global sensitivity analysis of model output.
\newblock {\em Technometrics}, 41(1):39--56, 1999.

\bibitem{sasikumar2023sensitivity}
Aravind Sasikumar, Albert Turon, Ivan~R C{\'o}zar, Oriol Vallmaj{\'o}, Jorge~Camacho Casero, Matthias De~Lozzo, and Said Abdel-Monsef.
\newblock Sensitivity analysis methodology to identify the critical material properties that affect the open hole strength of composites.
\newblock {\em Journal of Composite Materials}, 57(10):1791--1805, 2023.

\bibitem{rao2023frequency}
Yanni Rao, Wei Wan, Ning Wei, and Kui Wang.
\newblock Frequency-dependent dynamic moduli prediction for general bio-inspired staggered platelet reinforced composites.
\newblock {\em Polymer Composites}, 44(6):3268--3280, 2023.

\bibitem{zhou2020global}
Shuai Zhou, Yue Jia, and Chong Wang.
\newblock Global sensitivity analysis for the polymeric microcapsules in self-healing cementitious composites.
\newblock {\em Polymers}, 12(12):2990, 2020.

\bibitem{hamdia2018sensitivity}
Khader~M Hamdia, Hamid Ghasemi, Xiaoying Zhuang, Naif Alajlan, and Timon Rabczuk.
\newblock Sensitivity and uncertainty analysis for flexoelectric nanostructures.
\newblock {\em Computer Methods in Applied Mechanics and Engineering}, 337:95--109, 2018.

\bibitem{khaledi2016sensitivity}
Kavan Khaledi, Elham Mahmoudi, Maria Datcheva, D~K{\"o}nig, and Tom Schanz.
\newblock Sensitivity analysis and parameter identification of a time dependent constitutive model for rock salt.
\newblock {\em Journal of Computational and Applied Mathematics}, 293:128--138, 2016.

\bibitem{de2024thermo}
Eduardo A~Barros De~Moraes, Prakash KC, and Mohsen Zayernouri.
\newblock A thermo-electro-mechanical model for long-term reliability of aging transmission lines.
\newblock {\em arXiv preprint arXiv:2406.18860}, 2024.

\bibitem{kc2024thermo}
Prakash KC, Maryam Naghibolhosseini, and Mohsen Zayernouri.
\newblock Thermo-electro-mechanical modeling of power transmission line failures across four us states.
\newblock {\em arXiv preprint arXiv:2406.19603}, 2024.

\bibitem{kc2024multi}
Prakash KC, Maryam Naghibolhosseini, and Mohsen Zayernouri.
\newblock Multi-scenario and stochastic thermo-electro-mechanical modeling of failure in power transmission lines.
\newblock {\em arXiv preprint arXiv:2407.18857}, 2024.

\bibitem{blair_role_1947}
GW~Scott Blair.
\newblock The role of psychophysics in rheology.
\newblock {\em Journal of Colloid Science}, 2(1):21--32, 1947.

\bibitem{jaishankar_power-law_2013}
Aditya Jaishankar and Gareth~H McKinley.
\newblock Power-law rheology in the bulk and at the interface: quasi-properties and fractional constitutive equations.
\newblock {\em Proceedings of the Royal Society A: Mathematical, Physical and Engineering Sciences}, 469(2149):20120284, 2013.

\bibitem{blair_subjective_1942}
GW~Scott Blair and FMV Coppen.
\newblock The subjective conception of the firmness of soft materials.
\newblock {\em The American Journal of Psychology}, pages 215--229, 1942.

\bibitem{schiessel_generalized_1995}
H~Schiessel, R~Metzler, A~Blumen, and TF0921 Nonnenmacher.
\newblock Generalized viscoelastic models: their fractional equations with solutions.
\newblock {\em Journal of physics A: Mathematical and General}, 28(23):6567, 1995.

\bibitem{lion_thermodynamics_1997}
Alexander Lion.
\newblock On the thermodynamics of fractional damping elements.
\newblock {\em Continuum Mechanics and Thermodynamics}, 9:83--96, 1997.

\bibitem{rathinaraj_incorporating_2021}
Joshua David~John Rathinaraj, Gareth~H McKinley, and Bavand Keshavarz.
\newblock Incorporating rheological nonlinearity into fractional calculus descriptions of fractal matter and multi-scale complex fluids.
\newblock {\em Fractal and Fractional}, 5(4):174, 2021.

\bibitem{tzelepis_experimental_2023}
Demetrios~A Tzelepis, Jorge Suzuki, Yi~Feng Su, Yiyu Wang, Yong~Chae Lim, Mohsen Zayernouri, and Valeriy~V Ginzburg.
\newblock Experimental and modeling studies of ipdi-based polyurea elastomers--the role of hard segment fraction.
\newblock {\em Journal of Applied Polymer Science}, 140(10):e53592, 2023.

\bibitem{tzelepis_polyureagraphene_2023}
Demetrios~A Tzelepis, Arman Khoshnevis, Mohsen Zayernouri, and Valeriy~V Ginzburg.
\newblock Polyurea--graphene nanocomposites—the influence of hard-segment content and nanoparticle loading on mechanical properties.
\newblock {\em Polymers}, 15(22):4434, 2023.

\bibitem{kennedy_particle_1995}
James Kennedy and Russell Eberhart.
\newblock Particle swarm optimization.
\newblock In {\em Proceedings of ICNN'95-international conference on neural networks}, volume~4, pages 1942--1948. ieee, 1995.

\bibitem{kinsler2000fundamentals}
LE~Kinsler.
\newblock {\em Fundamentals of Acoustics}, volume 480.
\newblock John wiley and Sons, 2000.

\bibitem{saltelli_sensitivity_2002}
Andrea Saltelli.
\newblock Sensitivity analysis for importance assessment.
\newblock {\em Risk analysis}, 22(3):579--590, 2002.

\bibitem{saltelli_global_2008}
Andrea Saltelli, Marco Ratto, Terry Andres, Francesca Campolongo, Jessica Cariboni, Debora Gatelli, Michaela Saisana, and Stefano Tarantola.
\newblock {\em Global sensitivity analysis: the primer}.
\newblock John Wiley \& Sons, 2008.

\bibitem{zayernouri2011coherent}
M~Zayernouri and M~Metzger.
\newblock Coherent features in the sensitivity field of a planar mixing layer.
\newblock {\em Physics of Fluids}, 23(2), 2011.

\bibitem{norton_introduction_2015}
John Norton.
\newblock An introduction to sensitivity assessment of simulation models.
\newblock {\em Environmental Modelling \& Software}, 69:166--174, 2015.

\bibitem{borgonovo_sensitivity_2017}
Emanuele Borgonovo.
\newblock {\em Sensitivity analysis: an introduction for the management scientist}, volume 251.
\newblock Springer, 2017.

\bibitem{borgonovo_new_2001}
Emanuele Borgonovo and George~E Apostolakis.
\newblock A new importance measure for risk-informed decision making.
\newblock {\em Reliability Engineering \& System Safety}, 72(2):193--212, 2001.

\bibitem{keramati2024monte}
Hamed Keramati, Erik Birgersson, Sangho Kim, and Hwa~Liang Leo.
\newblock A monte carlo sensitivity analysis for a dimensionally reduced-order model of the aortic dissection.
\newblock {\em Cardiovascular Engineering and Technology}, pages 1--13, 2024.

\bibitem{saltelli_how_2010}
Andrea Saltelli and Paola Annoni.
\newblock How to avoid a perfunctory sensitivity analysis.
\newblock {\em Environmental Modelling \& Software}, 25(12):1508--1517, 2010.

\bibitem{douglas-smith_certain_2020}
Dominique Douglas-Smith, Takuya Iwanaga, Barry~FW Croke, and Anthony~J Jakeman.
\newblock Certain trends in uncertainty and sensitivity analysis: An overview of software tools and techniques.
\newblock {\em Environmental Modelling \& Software}, 124:104588, 2020.

\bibitem{puy_comprehensive_2022}
Arnald Puy, William Becker, Samuele~Lo Piano, and Andrea Saltelli.
\newblock A comprehensive comparison of total-order estimators for global sensitivity analysis.
\newblock {\em International Journal for Uncertainty Quantification}, 12(2), 2022.

\bibitem{razavi_what_2015}
Saman Razavi and Hoshin~V Gupta.
\newblock What do we mean by sensitivity analysis? the need for comprehensive characterization of “global” sensitivity in e arth and e nvironmental systems models.
\newblock {\em Water Resources Research}, 51(5):3070--3092, 2015.

\bibitem{sobol_sensitivity_1993}
IM~Sobo{\'l}.
\newblock Sensitivity estimates for nonlinear mathematical models.
\newblock {\em Math. Model. Comput. Exp.}, 1:407, 1993.

\bibitem{sobol_global_2001}
Ilya~M Sobol.
\newblock Global sensitivity indices for nonlinear mathematical models and their monte carlo estimates.
\newblock {\em Mathematics and computers in simulation}, 55(1-3):271--280, 2001.

\bibitem{mood_introduction_1950}
Alexander~McFarlane Mood.
\newblock {\em Introduction to the Theory of Statistics.}
\newblock McGraw-hill, 1950.

\bibitem{saltelli_variance_2010}
Andrea Saltelli, Paola Annoni, Ivano Azzini, Francesca Campolongo, Marco Ratto, and Stefano Tarantola.
\newblock Variance based sensitivity analysis of model output. design and estimator for the total sensitivity index.
\newblock {\em Computer physics communications}, 181(2):259--270, 2010.

\bibitem{homma_importance_1996}
Toshimitsu Homma and Andrea Saltelli.
\newblock Importance measures in global sensitivity analysis of nonlinear models.
\newblock {\em Reliability Engineering \& System Safety}, 52(1):1--17, 1996.

\bibitem{yang2012adaptive}
Xiu Yang, Minseok Choi, Guang Lin, and George~Em Karniadakis.
\newblock Adaptive anova decomposition of stochastic incompressible and compressible flows.
\newblock {\em Journal of Computational Physics}, 231(4):1587--1614, 2012.

\bibitem{zhang2014enabling}
Zheng Zhang, Xiu Yang, Ivan~V Oseledets, George~E Karniadakis, and Luca Daniel.
\newblock Enabling high-dimensional hierarchical uncertainty quantification by anova and tensor-train decomposition.
\newblock {\em IEEE Transactions on Computer-Aided Design of Integrated Circuits and Systems}, 34(1):63--76, 2014.

\bibitem{herman_salib_2017}
Jon Herman and Will Usher.
\newblock Salib: An open-source python library for sensitivity analysis.
\newblock {\em Journal of Open Source Software}, 2(9):97, 2017.

\bibitem{smith_uncertainty_2013}
Ralph~C Smith.
\newblock {\em Uncertainty quantification: theory, implementation, and applications}.
\newblock SIAM, 2013.

\bibitem{jansen_analysis_1999}
Michiel~JW Jansen.
\newblock Analysis of variance designs for model output.
\newblock {\em Computer Physics Communications}, 117(1-2):35--43, 1999.

\end{thebibliography}

\end{document}